\documentclass[11pt]{article}
\usepackage{graphicx,amsmath,amssymb,psfrag, natbib, amsfonts,enumerate,float, hyperref, ntheorem}
\usepackage[letterpaper, left=1in, top=1in, right=1in, bottom=1in,nohead,includefoot,verbose, ignoremp]{geometry}

\newcommand{\bfb}{\mbox{\boldmath $\beta$}}
\newcommand{\bfbeta}{\mbox{\boldmath $\beta$}}
\newcommand{\bfep}{\mbox{\boldmath $\epsilon$}}
\newcommand{\bfy}{\mbox{\boldmath $y$}}
\newcommand{\bfx}{\mbox{\boldmath $x$}}
\newcommand{\bfg}{\mbox{\boldmath $g$}}
\newcommand{\bone}{\mbox{\boldmath $1$}}
\newcommand{\bft}{\mbox{\boldmath $t$}}
\newcommand{\bfgambig}{\mbox{\boldmath $\gamma$}}
\newcommand{\bfgam}{\mbox{\boldmath \footnotesize {$\gamma$}}}
\newcommand{\bftheta}{\mbox{\boldmath $\theta$}}
\newcommand{\Mgam}{\mathcal{M}_{\bfgam}}
\newcommand{\Mo}{\mathcal{M}_0}
\newcommand{\bld}[1]{\mbox{\boldmath \footnotesize {$#1$}}}
\newcommand{\bldsm}[1]{\mbox{\boldmath \tiny {$#1$}}}
\newcommand{\Mt}{\mathcal{M}_T}

\newcommand\indep{\protect\mathpalette{\protect\independenT}{\perp}}
\def\independenT#1#2{\mathrel{\setbox0\hbox{$#1#2$}%
    \copy0\kern-\wd0\mkern4mu\box0}}

\newtheorem{thm}{Theorem}[section]
\newtheorem{cor}{Corollary}[section]
\newtheorem{lem}{Lemma}[section]

\newtheorem{cond}{Condition}[section]
\newtheorem*{remark*}{Remark}[section]

\title{\textbf{Block Hyper-$\bfg$ Priors in Bayesian Regression}}

\author{Agniva So$\mbox{m}^{\star}$, Christopher M. Han$\mbox{s}^{\dagger}$ and Steven N. MacEacher$\mbox{n}^{\dagger}$  \\
{\em Duke Universit$\mbox{y}^{\star}$ and The Ohio State Universit$\mbox{y}^{\dagger}$ }
}

\date{}

\begin{document}

\maketitle

\begin{abstract}
The development of prior distributions for Bayesian regression has traditionally been driven by the goal of achieving sensible model selection and parameter estimation.  The formalization of properties that characterize good performance has led to the development and popularization of thick tailed mixtures of $g$ priors such as the Zellner--Siow and hyper-$g$ priors.  The properties of a particular prior are typically illuminated under limits on the likelihood or the prior.  In this paper we introduce a new, \emph{conditional information asymptotic} that is motivated by the common data analysis setting where at least one regression coefficient is much larger than others.  We analyze existing mixtures of $g$ priors under this limit and reveal two new behaviors, \emph{Essentially Least Squares (ELS)} estimation and the \emph{Conditional Lindley's Paradox (CLP)}, and argue that these behaviors are, in general, undesirable.  As the driver behind both of these behaviors is the use of a single, latent scale parameter that is common to all coefficients, we propose a block hyper-$g$ prior, defined by first partitioning the covariates into groups and then placing independent hyper-$g$ priors on the corresponding blocks of coefficients.  We provide conditions under which {\em ELS} and the {\em CLP} are avoided by the new class of priors, and provide consistency results under traditional sample size asymptotics.
\end{abstract}


\section{Introduction}

Bayesian methods for regression address the central questions of model selection and parameter estimation. Conjugate forms have played an important role due to ease of specification and ease of update from prior to posterior.  Zellner's well-known $g$ prior \citep{zellner1986assessing} is fully conjugate for a normal theory regression model.  The conjugate form leads to quick calculation of the posterior distribution and of the marginal likelihood.  The first of these features allows one to compute estimates and make predictions for a given model while the second allows one to compare models via the Bayes factor and, with the addition of a set of prior model probabilities, to engage in model averaging.

The $g$ prior has many close cousins.  One popular formulation imposes conditionally (on an analog of $g$) independent prior distributions for the regression coefficients as in ridge regression \citep{hoerl1970ridge}.  Other formulations have been designed to handle the two main problems of estimation/prediction and model/variable selection.  For estimation/prediction, consistent estimators are desirable, and concern about prior-data conflict suggests nonlinear shrinkage.  Nonlinear shrinkage is produced by placing a scale mixture of normals \citep{west1987scale} on $g$ or on the variance of the regression coefficients.  The Bayesian lasso prior \citep{park2008bayesian,hans2009bayesian}, the orthant normal prior \citep{hans2011elastic}, the generalized double Pareto prior \citep{armagan2013generalized}, the Horseshoe prior \citep{carvalho2010horseshoe}, the normal-exponential-gamma/normal-gamma prior \citep{griffin2005alternative, griffin2010inference, griffin2011bayesian, griffin2012structuring}, and the priors of \cite{polson2010shrink} conform to a mixture representation of this form.  

The $g$ prior and the independence prior have been extended for use in model selection  by placing a prior distribution over the discrete space of regression models.  Each model consists of a subset of the covariates, with the other covariates implicitly having zero (or occasionally near-zero) regression coefficients.  The prior specific to each model is a $g$ (or independence) prior, often with a hyper-prior placing a distribution over $g$.  \cite{mitchell1988bayesian} provided a forerunner of these priors, the ``spike and slab'' prior.  Additional priors along these lines have been investigated in \cite{george1993variable}, \cite{george1997approaches}, \cite{foster1994risk}, \cite{kass1995reference}, \cite{george2000calibration}, \cite{fernandez2001benchmark}, \cite{johnstone2004needles} and  \cite{ishwaran2005spike}.  The computational cost of model selection is substantial, and the models are often fit with stochastic search algorithms.  Various algorithms having distinctive features include \cite{george1993variable}, \cite{berger2005posterior}, \cite{hans2007shotgun}, \cite{scott2008feature}, \cite{bott:10} and \cite{clyde2011bayesian}.  

In this work, we address model-specific estimation and prediction, model selection, and estimation and prediction under model averaging.  As such, \cite{liang2008mixtures}'s development of the hyper-$g$ prior is particularly relevant.  The hyper-$g$ prior, one of the most commonly used mixture of $g$ priors, has been designed to retain the computational efficiency of the $g$ prior while performing well for both estimation and model selection.  We revisit this performance by studying the limiting behavior of the model under a new, \emph{conditional information asymptotic}.  This type of limit is important in practice, as it provides insight into the behavior of estimation and model selection when one regression coefficient is substantially larger than other non-zero regression coefficients.  As part of the analysis, we describe \emph{Essentially Least Squares (ELS)} estimation wherein inference under a Bayesian model collapses to least squares estimation, and we identify a \emph{conditional version of Lindley's paradox (CLP)} which leads to inconsistent model selection. We show that many commonly used mixtures of $g$ priors, including the hyper-$g$ prior, suffer from both \emph{ELS} and the \emph{CLP}.  

In order to overcome the deficiencies of the scale mixture of $g$ priors, we introduce the block hyper-$g$ prior, a collection of ordinary mixture of $g$ priors applied to groups of predictors separately.  The theoretical properties of the new prior are investigated in detail under the analytically tractable blockwise orthogonal design setup, and the new prior is shown to perform well.  The new prior is suitable for situations where one can subjectively (or in an automated fashion) group predictors into blocks. The use of multiple prior variance parameters has been explored in \cite{maruyama2011fully}, \cite{rouder2012default} and \cite{min2012objective}.

In Section~\ref{sec:prelim}, we introduce notation that is used throughout the paper and describe the basic regression setup.  Section~\ref{sec:paradoxes} describes \emph{ELS} and \emph{CLP} behaviors, and contains theoretical results showing that several popular priors perform poorly on these criteria.  Section~\ref{sec:newprior} formulates the block hyper-$g$ prior and provides theoretical properties of the new prior under block orthogonal designs.  Section~\ref{sec:consistency} is dedicated to examining the consistency properties of the block hyper-$g$ prior.  We conclude the article with a brief discussion in Section~\ref{sec:disc}.  Proofs of the main results are in Appendix~\ref{app:proofs}, while proofs of other results are in Appendix \ref{appB} in the Supplementary Material.

\section{Notation and Preliminaries}\label{sec:prelim}

The basic regression problem can be described as explaining the behavior of the response vector $\bfy=(y_1, y_2, $ $ \ldots, y_n)^T$ using a known set of $p$ predictor variables $\bfx_1, \bfx_2, \ldots, \bfx_p$.  We consider the traditional setting where $n > p$.  Let $\bfgambig \in \Gamma=\{0,1\}^p$ denote the index set of the subsets of the predictor variables to be included/excluded in a model so that under a particular model $\mathcal{M}_{\bfgam}$, the $i^{th}$ element of the vector $\bfgambig$ signifies inclusion of predictor $i$ (in model $\mathcal{M}_{\bfgam}$) if $\gamma_i=1$ and exclusion if $\gamma_i=0$. Thus, $\Gamma$ describes the collection of all $2^p$ possible models and each element $\bfgambig$ represents a unique model in $\Gamma$. Let $X_{\bfgam}$ denote the $n \times p_{\gamma}$ design matrix and $\bfb_{\bfgam}$ denote the $p_{\gamma}$ length vector of regression coefficients corresponding to model $\Mgam$. The linear model $\Mgam$ can be represented as:
\[ \bfy = \bone \alpha + X_{\bfgam} \bfb_{\bfgam} + \bfep
\] 
where $\bone$ is a $n \times 1$ vector of 1's and $\alpha$ denotes the intercept which appears in every model in $\Gamma$. The vector of observed errors $\bfep=(\epsilon_1,\ldots,\epsilon_n)^\top$ is a Gaussian random vector, $\bfep \sim N(0,\sigma^2 I_n)$. The Bayesian approach places a prior on the vector of the unknown parameters $(\alpha,\bfb_{\bfgam},\sigma^2)=\bftheta_{\bfgam} \in \Theta_{\bfgam}$ corresponding to each model $\Mgam$ along with an additional prior on the model space $\Gamma$. To retain the same meaning for $\alpha$ across all models, we center the columns of $X_{\bfgam}$ so that $\bone^T X_{\bfgam}=0$ and all the covariates are orthogonal to the intercept. Traditionally, this transformation justifies specification of a common prior for $\alpha$ \citep{jeff:61,berger1998bayes}.  Recently \cite{bayarri2012criteria} provided alternative arguments for a common, flat prior on $\alpha$ based on a predictive matching criterion for the prior distribution. Without loss of generality, the  response vector $\bfy$ is also assumed to be centered at zero in our setup, i.e., $\bone^T\bfy=0$.

\subsection{Zellner's $\bfg$ prior}\label{sec:gprior}

Zellner's $g$ prior \citep{zellner1986assessing} is specified as
\begin{eqnarray}
 \pi (\alpha, \sigma^2 \mid \Mgam) &\propto& \frac{1}{\sigma^2} \nonumber  \\
 \bfb_{\bfgam} \mid g, \alpha, \sigma^2, \Mgam &\sim&  N \left( \textbf{0}, g\sigma^2 (X^T_{\bfgam} X_{\bfgam})^{-1} \right) 
 	\label{eq:zgprior}
\end{eqnarray}
which results in simple closed form expressions for the marginal likelihoods and Bayes factors. The Bayes factor for any model $\Mgam$ compared to the null model $\Mo$ is
\begin{eqnarray}\label{eq:gBF}
BF(\Mgam : \Mo)= \frac{p(\bfy \mid \Mgam)}{p(\bfy \mid \Mo)} = \frac{(1+g)^{(n-p_{\gamma}-1)/2}}{[1+g(1-R_{\bld{\gamma}}^2)]^{(n-1)/2}}
\end{eqnarray}
where $R^2_{\bfgam}$ is the coefficient of determination for the model $\Mgam$.  Under sum of squared error loss, the Bayes estimator of $\bfbeta_{\bfgam}$ is
\begin{eqnarray}
\label{eq:BayesEstimateg}
\widehat{\bfbeta}_{\bfgam} & = & \frac{g}{1+g} \widehat{\bfbeta}_{\bfgam,LS}, 
\end{eqnarray}
where $\widehat{\bfbeta}_{\bfgam,LS}$ is the least squares estimator of $\bfbeta_{\bfgam}$.

The $g$ prior requires a value for $g$.  A variety of choices have been suggested based on different considerations.  Well-known (fixed) $g$ priors include the \emph{unit information prior} \citep{kass1995reference}, the \emph{risk inflation criterion prior} \citep{foster1994risk}, the \emph{benchmark prior} \citep{fernandez2001benchmark}, the \emph{local empirical Bayes prior} \citep{hansen2001model} and the \emph{global empirical Bayes prior} \citep{george2000calibration}.  \cite{liang2008mixtures} review the fixed $g$ priors and summarize the justifications behind these specific values of $g$. 

\subsection{Mixtures of $\bfg$ priors}\label{sec:mixtures}
There are a variety of motivations for considering a ``fully Bayes'' approach where $g$ is modeled with a prior distribution, leading to a marginal prior for $\bfb_{\bfgam}$ which can be represented as a ``mixture of $g$ priors''
\[
	\pi(\bfb_{\bfgam} \mid \alpha, \sigma^2, \Mgam) = \int_0^\infty N \left( \bfb_{\bfgam} \mid \textbf{0}, g\sigma^2 (X^T_{\bfgam} X_{\bfgam})^{-1} \right)  \pi(g) dg.
\]
Careful choice of the mixing distribution can result in thick-tailed priors for $\bfb_{\bfgam}$ after marginalization of $g$.  An early example is the Zellner--Siow prior \citep{zellner1980posterior}, which can be expressed as a mixture of $g$ priors with an Inverse Gamma ($\frac{1}{2},\frac{n}{2}$) mixing density for $g$.  Mixing over $g$ also endows the Bayes estimator of $\bfb_{\bfgam}$ with data-adaptive shrinkage of the least squares estimator:
\begin{equation} 
	\widehat{\bfb}_{\bfgam} =  E(\bfb \mid \bfy, \Mgam) = E\left( \frac{g}{1+g} \mid \bfy, \Mgam \right)  \widehat{\bfb}_{\bfgam,LS}.
	\label{eq:shrink}
\end{equation}
The quantity $E(\frac{g}{1+g} \mid \bfy, \Mgam)$ is often called the \textit{shrinkage factor}.

A variety of prior distributions for $g$ have been considered in the literature. \cite{zellner1980posterior}, \cite{west2003bayesian}, \cite{maruyama2010robust}, \cite{maruyama2011fully} and \cite{bayarri2012criteria} are notable examples.  In this work, we focus primarily on the ``hyper-$g$'' prior of \cite{liang2008mixtures}.

\subsubsection{Hyper-$\bfg$ Priors}\label{sec:hyper}
The ``hyper-$g$'' prior proposed by \cite{liang2008mixtures} places a prior distribution on $g$ with density
\[ 
\pi(g) = \frac{a-2}{2} (1+g)^{-a/2} \;,\; \; g>0. 
\] 
The prior is proper if $a > 2$ and the authors suggest using $a \in (2, 4]$, with $a=3$ being the default choice.
The Bayes factor under the hyper-$g$ prior can be expressed as
\begin{eqnarray*}
BF(\Mgam : \Mo) &=& \int_0^{\infty} (1+g)^{(n-p_{\gamma}-1)/2}[1+g(1-R_{\bld{\gamma}}^2)]^{-(n-1)/2} \pi(g)dg \\
&=& \frac{a-2}{p_{\gamma}+a-2}  \; _2F_1 \left(\frac{n-1}{2},1;\frac{a+p_{\gamma}}{2};R_{\bld{\gamma}}^2\right) 
\end{eqnarray*}
where $_2F_1(\cdot)$ is the Gaussian hypergeometric function.

\section{Asymptotic Evaluations: Paradoxes, Old and New}\label{sec:paradoxes}

While the $g$ priors described in Section~\ref{sec:gprior} offer many conveniences, they are known to have several undesirable
properties commonly referred to as ``paradoxes.'' In purest form, the paradoxes associated with $g$ priors are revealed when taking a limit.  \cite{liang2008mixtures} describe two such paradoxes which arise from different limits.  The first, Lindley's Paradox, relies on a limit which weakens the prior distribution.  The second, the Information Paradox, relies on a limit where the signal in the data becomes stronger.  Both limits hold the design $X$ (and hence sample size) fixed.  These two ``old'' paradoxes can be summarized as follows.

\emph{Bartlett's Paradox/Lindley's Paradox} : Lindley's Paradox is an anomaly associated with a fixed $g$ prior when the scale factor $g$ is intentionally chosen to be large in an attempt to make the prior weakly informative.  Holding the data $(X_{\bfgam}, \bfy)$ fixed, the Bayes factor \eqref{eq:gBF} comparing any arbitrary non-null model $\Mgam$ to the null model $\Mo$ approaches zero in the limit when $g \rightarrow \infty$, irrespective of the data. The full description of the paradox contrasts this undesirable behavior with the results of a classical hypothesis test \citep{lind:57, bart:57, jeff:61, liang2008mixtures}.

\emph{Information Paradox} : The Information Paradox is associated with a strong signal in the data, as manifested by a high value of  $R_{\bfgam}^2$.  Holding $(X_{\bfgam}, \bfep)$ fixed, let $||\bfb_{\bfgam}|| \rightarrow \infty$, so that $||\widehat{\bfb}_{\bfgam,LS}|| \rightarrow \infty$ and $R_{\bfgam}^2 \rightarrow 1$.  It follows from \eqref{eq:gBF} that $BF(\Mgam:\Mo) \rightarrow (1+g)^{(n-p_{\bldsm{\gamma}}-1)/2}$, a finite constant.  Thus the Bayes factor for $\Mgam$ relative to $\Mo$ is bounded even though the likelihood evidence in favor of $\Mgam$ grows without bound \citep{zellner1986assessing, berger2001objective, liang2008mixtures}.

These undesirable properties can be avoided by using mixtures of $g$ priors with a careful choice of mixing distribution.
\cite{liang2008mixtures} provide sufficient conditions under which a prior $\pi(g)$ resolves the Information Paradox, and prove
that the hyper-$g$ prior avoids both of the above paradoxes \citep[as does the robust prior of][]{bayarri2012criteria}.  
While these ``old'' paradoxes have been studied extensively, the limits taken to produce them have further, less well known implications.  The first is initially seen with the limit in Lindley's Paradox.  The second, a new paradox, follows from a modification to \cite{liang2008mixtures}'s limit.  Qualitative descriptions of these behaviors are as follows, with formal results provided in
Section~\ref{sec:clp}.

\emph{Essentially Least Squares Estimation (ELS)} : It is well-known that under the $g$ prior, the Bayes estimator of $\bfbeta_{\bfgam}$ in \eqref{eq:BayesEstimateg} tends to $\widehat{\bfbeta}_{\bfgam,LS}$ as $g \rightarrow \infty$.  Formally, we identify \emph{ELS} behavior as $||\widehat{\bfbeta}_{\bfgam} - \widehat{\bfbeta}_{\bfgam,LS}|| / ||\widehat{\bfbeta}_{\bfgam,LS}|| \rightarrow 0$ under some appropriate limit.  In Sections~\ref{sec:asymp} and \ref{sec:clp}, we consider limits which are driven by changes to the data rather than changes to the prior and show that several common mixtures of $g$ priors exhibit {\em ELS} behavior.

\emph{Conditional Lindley's Paradox (CLP)} : The Conditional Lindley's Paradox arises when comparing two models $\mathcal{M}_{\bfgam_1}$ and $\mathcal{M}_{\bfgam_2}$ with $\mathcal{M}_{\bfgam_1} \subset \mathcal{M}_{\bfgam_2}$, where $\mathcal{M}_{\bfgam_2}$ is the ``correct'' model.  Specific asymptotics for the data (described explicitly in Section~\ref{sec:asymp}) yield $BF(\mathcal{M}_{\bfgam_2}: \mathcal{M}_{\bfgam_1}) \rightarrow 0$, compelling one to accept the smaller (incorrect) model.  

Before connecting these behaviors to existing mixtures of $g$ priors, we describe the limits driving the phenomena.

\subsection{A Conditional Information Asymptotic}\label{sec:asymp}
We consider an asymptotic analysis of a sequence of problems, where each element in the
sequence is related to the linear regression model $\bfy=\alpha \bone +X \bfb + \bfep$.  The design matrix $X$ 
is an $n \times p$ matrix with full column rank, and the columns of $X$ and the response $\bfy$ 
are centered.  Specifically, we write the linear model as 
$\bfy = \alpha \bone + X_1 \bfb_1 + X_2 \bfb_2 + \bfep$, where $X = (X_1, X_2)$, $X_1$ is an $n\times p_1$ 
matrix, $X_2$ is an $n \times p_2$ matrix and $\bfb = (\bfb_1^T, \bfb_2^T)^T$.  We construct a sequence 
$\{ \Psi_N \}_{N=1}^{\infty}$ where each element $\Psi_N$
represents the linear model with
\begin{eqnarray}\label{eq:psiseq}
\Psi_N &=& (X_{1(N)},X_{2(N)},\alpha_N,\bfb_{1(N)}, \bfb_{2(N)},\bfep_N) \nonumber \\
&=& (X_1,X_2,\alpha,\bfb_{1(N)},\bfb_{2},\bfep),
\end{eqnarray}
and $||\bfb_{1(N)}|| \rightarrow \infty$ as 
$N  \rightarrow \infty$ while $X_1, X_2, \alpha,\bfb_{2}$ and $\bfep$ are held fixed.  This is a fixed-$n$, fixed-$p_\gamma$
asymptotic, and represents a strengthening of the likelihood that is driven by one particular set of predictor variables.

We refer to this as a \textit{conditional information asymptotic}, as it can be viewed as the limit that drives 
the Information Paradox of \cite{liang2008mixtures} (i.e., $|| \bfb|| \rightarrow \infty$) with the additional condition that a portion 
of $\bfb$ remains
fixed in the analysis.  The consequences of the information limit considered in \cite{liang2008mixtures} were driven by 
$R^2 \rightarrow 1$.  The following lemma notes that the conditional information asymptotic produces the same
behavior. 

\begin{lem}\label{lem3.1}
Let $R^2 (N)$ denote the coefficient of determination for element $N$ in the sequence $\{ \Psi_N \}$ as defined in \eqref{eq:psiseq}. Then $R^2 (N) \rightarrow 1$ as $N \rightarrow \infty$.  
\end{lem}

The lemma follows immediately by noting that the error vector $\bfep$ is fixed, and hence $\widehat{\sigma}^2 = \frac{||\bld{y}- \widehat{\alpha}_{LS} \bld{1} - X \widehat{\bld{\beta}}_{LS}||^2}{n-p-1}$ also remains unchanged.  $R^2(N) = ||X\widehat{\bfb}_{LS}||^2/[ ||X\widehat{\bfb}_{LS}||^2 + (n-p-1)\widehat{\sigma}^2]$ which tends to one as $||\bfb_1|| \rightarrow \infty$.

\subsection{A Conditional Lindley's Paradox}\label{sec:clp}
In this section we investigate the behavior of several mixtures of $g$ priors under the conditional information asymptotic defined
by \eqref{eq:psiseq}. To streamline notation, we drop the subscript $\bfgambig$ from $\Mgam$ and refer to $R^2 (N)$ simply as $R^2$.  The following results apply to an arbitrary model $\Mgam$ unless otherwise mentioned. 
The first theorem reveals a behavior of the Bayes estimator ($\widehat{\bfb}$) under the hyper-$g$ prior.

\begin{thm}\label{thm3.1}
(ELS) Under the hyper-$g$ prior, for the sequence $\{ \Psi_N \}$ defined in \eqref{eq:psiseq}, $\frac{||\widehat{\bld{\beta}} - \widehat{\bld{\beta}}_{LS}|| }{ ||\widehat{\bld{\beta}}_{LS}||} \rightarrow 0$  as $N\rightarrow \infty$, provided $n \geq a + p - 1$.  
\end{thm}

The proof of the theorem is in Appendix \ref{appa1}. The theorem roughly indicates that when at least one of the coefficients in the model is large, the estimates under the hyper-$g$ prior are \emph{Essentially Least Squares}. The behavior of the hyper-$g$ prior in such a situation runs counter to the conventional wisdom that, for a low-information prior, small (near zero) coefficients should be shrunk substantially while larger coefficients should be left unchanged \citep{berger1985statistical}.  The next theorem, proved in Appendix~\ref{appa2}, shows that the hyper-$g$ prior suffers from what we call the \emph{Conditional Lindley's Paradox}.

\begin{thm} \label{thm3.2} 
(CLP)
Consider the two models $M_1$ and $M_2$ such that
\begin{eqnarray} \label{TwoModels}
 M_1 &:& \bfy = \alpha \bone + X_1 \bfb_1 + \bfep \nonumber \\ 
 M_2 &:& \bfy = \alpha \bone + X_1 \bfb_1 + X_2 \bfb_2 + \bfep 
\end{eqnarray}
where $\bfb_i$ is a vector of length $p_i > 0 \; \mbox{for} \; i=1,2$ and $p_1+p_2=p$. 
Under the hyper-$g$ prior, when $||\bfb_1|| \rightarrow \infty \; (i.e, N \rightarrow \infty)$ in the sequence $\{ \Psi_N \}$ defined in \eqref{eq:psiseq} and $n \geq a+p_1-1$, the Bayes factor $BF(M_2:M_1)$ comparing model $M_2$ to model $M_1$ goes to zero, irrespective of the data.  
\end{thm}

The import of the theorem is that when comparing a pair of nested models, if at least one of the regression coefficients common to both models is large compared to the additional coefficients in the bigger model, the Bayes factor under the hyper-$g$ prior will place too much weight on the smaller model.  The limiting case with the size of the common coefficients growing infinitely large results in choice of the small model with probability tending to 1, leading to the conclusion that the predictors $X_2$ are, with certainty, unrelated to the response.  This behavior is unsettling since no matter how different from zero the additional coefficients in the big model are, the more important predictor(s) common to both models drive model choice toward the small model.  

We note that the behavior of the Bayes factor under this limit cannot be attributed to $\sigma^2$, as the posterior of 
$\sigma^2$ is well-behaved and converges to a proper probability distribution.  

\begin{cor} \label{cor3.1}
Under the hyper-$g$ prior, the posterior distribution of $\sigma^2$ in the sequence of problems $\{\Psi_{N}\}$ defined in \eqref{eq:psiseq} converges to an $IG\left(\frac{n+1-a-p}{2},\frac{2}{(n-p-1)\widehat{\sigma}^2} \right)$ distribution when $n > a+p-1$.    
\end{cor}

A proof is in Appendix~\ref{appb1} in the supplementary materials.
The \emph{CLP} coincides with Lindley's Paradox for any fixed $g$ prior, and as such it is easily shown that fixed $g$ priors are also adversely affected by the \emph{CLP}.

\subsubsection{The CLP in Other Mixtures of $\bfg$ Priors}

Theorem~\ref{thm3.2} describes the {\em CLP} under the hyper-$g$ prior. Other mixtures of $g$ priors exhibit the same behavior.
For example, \cite{maruyama2011fully} develop a generalized $g$ prior by specifying a prior on $\bfb$ through a prior on the rotated vector $W^T\bfb$, where $W$ is defined through the singular value decomposition $X=UDW^T$. In the simple situation where $X^T X$ is a block diagonal matrix, the generalized $g$ prior suffers from the \emph{CLP}.

\begin{thm}\label{thm3.3}
Consider the models in (\ref{TwoModels}) with the assumption that $X_1 \indep X_2$.  Under the generalized $g$ prior of \cite{maruyama2011fully},  as $N \rightarrow \infty$ in the sequence $\{ \Psi_N \}$ defined in \eqref{eq:psiseq}, $BF(M_2:M_1) \rightarrow 0$ irrespective of the data.  
\end{thm}

The robust prior of \cite{bayarri2012criteria} also suffers from the \emph{CLP} whether or not $X^T X$ is block diagonal. 

\begin{thm}\label{thm3.4}
Consider the models in (\ref{TwoModels}) with $n > p_1 + 2$.  Under the robust prior of \cite{bayarri2012criteria} with the recommended hyperparameters ($a=1/2, b=1$ and $\rho=\frac{1}{p+1}$), as $N \rightarrow \infty$ in the sequence $\{ \Psi_N \}$ defined in \eqref{eq:psiseq}, $BF(M_2:M_1) \rightarrow 0$ irrespective of the data.  
\end{thm}

Proofs are in Appendices~\ref{appengeneralizedg} and \ref{appenb3} in the supplementary materials.

\section{Avoiding \emph{ELS} and \emph{CLP} Behaviors}\label{sec:newprior}

The \emph{ELS} and \emph{CLP} behaviors described in Section~\ref{sec:clp} for mixtures of $g$ priors arise as a result of the use of a single, latent scale parameter $g$ that is common to each predictor variable.  In order for the model to fit the data in the presence of one (or more) large coefficients, $g$ must be large (with high probability).  Because $g$ affects estimation of all coefficients \eqref{eq:shrink}, this has the side-effect that small coefficients are not shrunk, producing \emph{ELS} behavior.  The \emph{CLP} can be explained by an argument similar to the one that explains Lindley's Paradox: as the common parameter $g$ is driven to be larger and larger (by a portion of the data, in our case) the diminishing prior mass in the neighborhood near zero containing any small, nonzero coefficients effectively rules out these predictors.

As we show in this section, these behaviors can be avoided through the use of multiple latent mixing parameters in place of a single, common $g$.  This approach has a connection to the concept of ``local shrinkage,'' which has a rich history in the study of the related normal means problem, e.g.~\cite{strawderman1971proper} and \cite{berger1980robust}.  Recent research in this area includes \cite{scott:06}, \cite{carvalho2010horseshoe}, \cite{scott2010bayes}, \cite{polson2012half} and \cite{bhattacharya2014dirichlet}.
The use of multiple latent scale parameters in regression settings has typically focused on ridge-regression-like settings where regression coefficients are conditionally independent \emph{a priori} \citep[e.g.,][]{polson2010shrink, armagan2013generalized}.
\cite{polson2012local} consider local shrinkage in regression where the local shrinkage parameters are attached to linear combinations of the regression coefficients.  A similar setting is considered by \cite{west2003bayesian}.

Our approach is to endow collections of regression coefficients with their own, independent, mixture of $g$ priors.  Having latent
scale parameters $g_i$ that are local to collections of coefficients results in models that avoid \emph{ELS} and \emph{CLP} behavior under the conditions described in Section~\ref{sec:blockhgstudy}.  The extreme case where each predictor variable is associated with its own $g_i$ was described, but not pursued, by \cite{liang2008mixtures} as representing ``scale mixtures of independent $g$ priors.''  The approach we propose emerges as a more general version of this idea with added theoretical underpinning related to \emph{ELS} and the \emph{CLP}.


\subsection{Block $\bfg$ Priors}\label{sec:blockg}
We build a block $g$ prior distribution by partitioning the predictors into $k$ blocks, $X=(X_1, X_2, \ldots, X_k)$.  
$X_i$ is a submatrix of dimension $n\times p_i$, $i=1,2,\ldots,k$.  The subscript $\bfgambig$ is suppressed here to simplify notation.  The regression setup for the block $g$ prior is
\begin{eqnarray} \label{eq:bhg}
\bfy \mid \alpha, \bfb, \sigma^2 &\sim& N( \alpha \bone + X\bfb , \sigma^2 I), \nonumber \\
\bfb \mid \bfg, \sigma^2 &\sim& N(\boldsymbol{0},A \sigma^2), \\
\mbox{and} \ \pi(\alpha,\sigma^2) &\propto& \frac{1}{\sigma^2}, \nonumber 
\end{eqnarray}
where $A$ is a block diagonal matrix defined as
\begin{equation*}
 A = \left( \begin{array}{cccc}
    g_1 (X_1^T X_1)^{-1} & 0 &  \cdots & 0 \\
    0 & g_2 (X_2^T X_2)^{-1} & \cdots & 0 \\ 
    \vdots & \vdots &  \ddots & \vdots \\          
		0 & 0 & \cdots & g_k (X_k^T X_k)^{-1} \\
		\end{array} \right). \\
\end{equation*}		

The regression coefficients for the $k$ distinct groups of predictors are taken to be independent \textit{a priori}.  The separate scales $g_i$ allow differential shrinkage on distinct blocks, with the amount of shrinkage on a block governed almost exclusively by the block itself.  The block $g$ prior reduces to the ordinary $g$ prior when the design matrix X is block orthogonal with $k$ orthogonal blocks $X_1, X_2, \ldots, X_k$ and $\bfg=(g_1,g_2,\ldots,g_k)^T=g \bone_k$.  

The blocks allow us to capture modeling concepts in our analysis.  In applied work, we motivate blocking in various ways \citep{som:14}.  The predictors that comprise a block may be different measures of a latent construct, they may represent a group of indicators, and so on.  In such cases, it is essential that the analysis be invariant to certain reparameterizations so that the coding of predictors within a block does not affect the analysis.  The block $g$ prior has this invariance (see Appendix~\ref{appb7}).  

\subsubsection{Block Hyper-$\bfg$ Prior}\label{sec:blockhyperg}

The block hyper-$g$ prior arises as a specific mixture of block $g$ priors where the mixing distributions on the components of the vector $\bfg$ are independent hyper-$g$ priors:
\begin{equation}
	\pi(\bfg) = \prod_{i=1}^k \frac{a-2}{2}(1 + g_i)^{-a/2}, \, \ g_i > 0. \label{eq:priorg}
\end{equation}
We follow the recommendations in \cite{liang2008mixtures} regarding the choice of the hyperparameter $a$ and take $2 < a \leq 4$, with $a=3$ the default choice.  A more general form of the prior would allow each $g_i$ to have density $\pi_i$ with block-specific hyperparameters $a_i$.  To do so would, in our view, unnecessarily complicate the analysis.  While we study the behavior of the block hyper-$g$ prior in this paper, one can envision different versions of the block $g$ prior where other mixing distributions are used for the $g_i$ \citep[e.g., a block robust prior \emph{\`a la}][]{bayarri2012criteria}.

\subsection{Asymptotic Evaluations of the Block Hyper-$\bfg$ Prior}\label{sec:blockhgstudy}

We evaluate the behavior of the block hyper-$g$ prior under a limit similar to the one described in Section~\ref{sec:asymp}. Consider the model in \eqref{eq:bhg} and \eqref{eq:priorg}, with $\bfy = \alpha \bone + X_1 \bfb_1 + X_2 \bfb_2 + \cdots + X_k \bfb_k + \bfep$.  We define a sequence of problems $\Psi_N = (X_{1(N)}, \ldots, X_{k(N)},\alpha_N,$ $\bfb_{1(N)}, \ldots, \bfb_{k(N)},\bfep_N)$ where the only quantity that changes in the sequence is the group of regression parameters $\bfb_{1(N)}$:
\begin{eqnarray}\label{eq:psiseq2}
\Psi_N  = (X_{1},\ldots, X_{k}, \alpha, \bfb_{1(N)}, \bfb_2, \ldots, \bfb_{k}, \bfep)
\end{eqnarray} 
with  $ ||\bfb_{1(N)}|| \rightarrow \infty$ as $N \rightarrow \infty$.
We use the following condition in the remainder.  

\begin{cond}\label{eq:bhgorth}
The predictors and the response are centered and the design matrix is block orthogonal:
\begin{eqnarray*}
\bone \indep \bfy, X_1, X_2,\ldots,X_k  \; \mbox{and} \; X_i \indep X_j, \; \mbox{where} \; i \neq j.
\end{eqnarray*}
\end{cond}

The assumption of block orthogonality facilitates asymptotic analysis by providing simpler expressions for many posterior summaries. This condition is not essential to define the block $g$ or block hyper-$g$ prior.  Rather, it leads to the subsequent theoretical results.  Block orthogonality is commonly encountered in designed experiments and in analyses where covariates have been successively orthogonalized.   In cases where Condition \ref{eq:bhgorth} is not satisfied, a variation of the block $g$ prior \citep{som:14} can be used, for which variants of the ensuing results hold.  These results will be reported elsewhere.

Before providing the main results we summarize several aspects of the posterior distribution of the regression model defined by \eqref{eq:bhg} and \eqref{eq:priorg} subject to Condition~\ref{eq:bhgorth}.
The posterior mean of $\bfb$ given $\bfg$ is
\begin{equation}
\widehat{\bfb} = E(\bfb \mid \bfg, \bfy) = \left( \frac{g_1}{1+g_1} \widehat{\bfb}_{1,LS}^T ,\ldots, \frac{g_k}{1+g_k}\widehat{\bfb}_{k,LS}^T \right)^T, \label{eq:postmean}
\end{equation}
where $\widehat{\bfb}_{i,LS}$ denotes the component of the least squares estimator $\widehat{\bfb}_{LS}$ corresponding to block $i$.
The posterior density of $\bfg$ is
\begin{eqnarray*}
\pi(\bfg \mid \bfy) &\propto& \frac{\prod_{i=1}^k (1+g_i)^{-\frac{a+p_i}{2}} }{\|\bfy\|^{n-1} \left[ 1 - \sum_{i=1}^k \frac{g_i}{1+g_i} R_i^2 \right]^{(n-1)/2} }
\end{eqnarray*}
where $R_i^2 = \frac{\bld{y}^T P_{X_i}\bld{y}}{\bld{y}^T \bld{y}} ,\; i=1,2,\ldots,k$ and $P_{X_i}$ is the projection matrix for the column space of $X_i$. It will be useful to define $t_i =\frac{g_i}{1+g_i}$ for $i=1,\ldots,k$ so that, under a block orthogonal design, each $t_i$ represents the shrinkage factor for the $i^{th}$ block under a block $g$ prior with fixed $g$. Then 
\begin{eqnarray*}
\pi(\bft \mid \bfy) \propto \prod_{i=1}^k (1-t_i)^{\frac{a+p_i}{2}-2} (1-\sum_{i=1}^k t_i R_i^2)^{-\frac{n-1}{2}}.
\end{eqnarray*}

The following lemma, similar to Lemma \ref{lem3.1}, is the building block for the main results in this section.

\begin{lem}\label{lem4.1}
For the regression model described by \eqref{eq:bhg} and \eqref{eq:priorg} and satisfying Condition~\ref{eq:bhgorth}, as $N \rightarrow \infty$  in the sequence $\{ \Psi_N \}$ defined in \eqref{eq:psiseq2}, $R_1^2 \rightarrow 1$ and $R^2_i \rightarrow 0, \; \forall \; i\neq 1$.
\end{lem}

The first of the two main results in this section shows that \emph{ELS} behavior is avoided under the block hyper-$g$ prior.

\begin{thm} \label{thm4.1}
For the regression model described by \eqref{eq:bhg} and \eqref{eq:priorg} and satisfying Condition~\ref{eq:bhgorth}, 
\[
	E(\bfb \mid \bfy)  = \left( \begin{array}{c}
		E\left(\frac{g_1}{1+g_1}\mid \bfy \right) \widehat{\bfb}_{1,LS} \\
		\vdots \\
		E\left(\frac{g_k}{1+g_k}\mid \bfy \right) \widehat{\bfb}_{k,LS} 
		\end{array} \right).
\]
Further assume that $n \geq a+p_1-1$. Then, as $N \rightarrow \infty$  in the sequence $\{ \Psi_N \}$ defined in \eqref{eq:psiseq2}, $E\left(\frac{g_1}{1+g_1}\mid \bfy \right) \rightarrow 1$ and, for $i\neq 1$, $E\left(\frac{g_i}{1+g_i}\mid \bfy \right) \rightarrow \Delta_i$ with $ \frac{2}{a+p_i} \leq \Delta_i < 1 $.  
\end{thm}
The theorem is proved in Appendix~\ref{appa4}.  Loosely, the theorem says that coefficients in the block with at least one large coefficient show \emph{ELS} behavior, while coefficients in the other blocks (where all coefficients are relatively smaller) do not display \emph{ELS} behavior and are shrunk.  The amount of shrinkage for the (relatively) small coefficients (blocks $i \neq 1$) is driven largely by the ratio of $\bfy^T P_{X_i} \bfy$ and $\widehat{\sigma}^2$.  The lower bound $\frac{2}{a+p_i}$ occurs when $R^2_i=0$ in which case $\widehat{\bfb}_{i,LS}=0$.  The block-specific shrinkage is the key to avoiding the \emph{CLP}, as is shown in the second main result of this section.  

\begin{thm}\label{thm4.2}
Consider the two models $M_1$ and $M_2$ in (\ref{TwoModels}) with a block hyper-$g$ prior on $\bfb$ (with $k=2$) as in \eqref{eq:bhg} and \eqref{eq:priorg}. Assume that blocks $X_1$ and $X_2$ satisfy Condition~\ref{eq:bhgorth}  on the design. When $||\bfb_1|| \rightarrow \infty$ ($N \rightarrow \infty$) in the sequence $\{\Psi_N\}$ defined in \eqref{eq:psiseq2}, the Bayes factor $BF(M_2:M_1)$  is bounded away from zero.
\end{thm}

This theorem, proved in Appendix~\ref{appa45}, shows that the block hyper-$g$ prior avoids the {\em CLP}.  The limiting posterior distribution does not declare the regression coefficients for the second block to be zero and does not concentrate on $M_1$.  We believe this to be appropriate behavior, as the least squares estimate of these coefficients does not change in the sequence $\{ \Psi_N \}$ and the evidence in the data is neither conclusive that the coefficients are zero or that they are not zero.

As in Section~\ref{sec:clp}, we note that the posterior distribution of $\sigma^2$ in the sequence $\{ \Psi_N \}$ under the block hyper-$g$ prior is well-behaved. The posterior for $\sigma^2$ does not converge to a standard distribution as $N \rightarrow \infty$, but Corollaries \ref{cor4.1} and \ref{cor4.2} show that the center of the limit distribution is finite and has an upper bound that can be easily calculated.  The proofs of both corollaries can be found in supplementary materials (Appendices \ref{appb3} and \ref{appb4}).  

\begin{cor}\label{cor4.1}
Consider a regression model of the form \eqref{eq:bhg} and \eqref{eq:priorg} which satisfies Condition~\ref{eq:bhgorth} and let $p=\sum_{i=1}^k p_i$. When $n > k(a-2)+p+1$, as $N \rightarrow \infty$ in the sequence $\{\Psi_N\}$ defined in \eqref{eq:psiseq2}, the sequence of posteriors of $\sigma^2$ converges to the distribution $F(\cdot)$ with density 
\[ f(\sigma^2) \propto \frac{1}{(\sigma^2)^{\frac{n+1}{2}-\frac{ka+p}{2}+k}} \exp \Big[ -\frac{(n-p-1)\widehat{\sigma}^2}{2\sigma^2} \Big] \prod_{i=2}^k \gamma \Big(\frac{a+p_i}{2}-1,\frac{(X_i\widehat{\bld{\beta}}_i)^T (X_i\widehat{\bld{\beta}}_i)}{2\sigma^2} \Big) \] 
where $\gamma(s,x)= \int_0^x t^{s-1} e^{-t} dt$ is the lower incomplete gamma function.
\end{cor}

\begin{cor}\label{cor4.2}
Consider a regression model of the form \eqref{eq:bhg} and \eqref{eq:priorg} which satisfies Condition~\ref{eq:bhgorth}.  Then in the sequence $\{ \Psi_N \}$ defined in \eqref{eq:psiseq2}, \[ \lim\limits_{N \rightarrow \infty} E(\sigma^2 \mid \bfy) \leq \frac{1}{n-1-a-p_1} \Big[ (n-p-1) \widehat{\sigma}^2 + \sum_{i=2}^k (X_i \widehat{\bfb_i})^T (X_i \widehat{\bfb_i}) \Big] \] when $n > a+p_1+1$. 
\end{cor}

Note that the bound is finite for all problems in the sequence $\{\Psi_N \}$. The upper bound for the expectation is achieved by the ordinary hyper-$g$ prior (i.e., a block hyper-$g$ prior with $k=1$). In this case
\[ \lim_{N \rightarrow\infty} E(\sigma^2 \mid \bfy) = \frac{(n-p-1) \widehat{\sigma}^2}{n-p-a-1}
\]
which is consistent with  Corollary \ref{cor3.1}.

\section{Consistency of the Block Hyper-$\bfg$ Prior}\label{sec:consistency}

In this section, we analyze the block hyper-$g$ prior with respect to three existing notions of consistency: \emph{information consistency}, \emph{model selection consistency} and \emph{prediction consistency}.  All three are considered by \cite{liang2008mixtures} with respect to the hyper-$g$ prior, and the first two among the seven ``criteria for Bayesian model choice'' posited by \cite{bayarri2012criteria}.

\subsection{Information Consistency}
This form of consistency is directly related to the \textit{Information Paradox} described in Section 3.1.  \cite{liang2008mixtures} define a Bayesian normal linear regression model under a particular prior to be \emph{information consistent} if, under an appropriate limit on the data vector $\bfy$ for a \emph{fixed} sample size $n$,  $R^2_{\bfgam} \rightarrow 1$ for model $\Mgam$ implies $BF(\Mgam : \Mo) \rightarrow \infty$, where $\Mo$ is the null model.  \cite{bayarri2012criteria} provide a formal definition of \emph{information consistency} that applies to models other than the normal linear model.  The following theorem establishes information consistency of the block hyper-$g$ prior.

\begin{thm}\label{thm5.1}
Consider a regression model of the form \eqref{eq:bhg} and \eqref{eq:priorg} which satisfies Condition~\ref{eq:bhgorth}. The block hyper-$g$ prior is ``information consistent"  when either of two sufficient conditions hold: \\
(1) $R^2 \rightarrow 1$ and  $n > k(a-2)+p+1$, where $k$ is the total number of blocks, $p=\sum_{j=1}^k p_{j}$ and $p_{j}$ is the size of block $X_{j}$.  \\
(2) For some $i=1,\ldots,k$, $R_i^2 \rightarrow 1$ and  $n \geq a+p_{i}-1$, where $R_{i}^2$ is the component of $R^2$ corresponding to the $i^{th}$ orthogonal block.
\end{thm}

The proof of the theorem is in Appendix \ref{appa5}. Note that Condition (1) provides a form of the theorem where $R^2$ approaches 1 due to growth in size of any arbitrary set of coefficients.  Under Condition (2), the coefficients growing in size all belong to a single block.  Condition (1) requires a larger sample size than does Condition (2).

\subsection{Conditions and Assumptions}\label{sec:assume}
For the remaining two consistency results we revert to the traditional asymptotic setting where parameters are held fixed and the sample size increases.  Before proceeding, we first need to fix the notion of a ``true'' model from which the data $\bfy$ are assumed to have been generated. Assume that $B_T \subseteq \{1,2,\ldots,k\}$ denotes the indices of the blocks included in the true model $\Mt$, where each block has at least one non-zero coefficient.  Then $\Mt: \bfy = \alpha_{T} \bone +  X_{T} \bfb_{T} + \bfep = \alpha_{T} \bone + \sum\limits_{i \in B_T} X_{i,T} \bfb_{i,T} + \bfep$ denotes the true data generating process. Under the model $\Mt$ there are $|B_T| = k_T$ different blocks with separate and independent hyper-$g$ priors on each block.  The following basic model assumptions and conditions are used in the results.

\begin{cond}\label{cond5.1}
The $n \times p$ design matrix $X$ grows in size with the restriction that 
\[ \lim_{n \rightarrow \infty} \frac{1}{n} X^T X = D,  \]
 for some $p \times p$ positive definite matrix $D$.
\end{cond}

The following condition is a direct consequence of Condition \ref{cond5.1} \citep{maruyama2011fully}.

\begin{cond}\label{cond5.2}
The design allows the true model $\Mt$ and any arbitrary model  $\Mgam \not \supseteq \Mt$ to be asymptotically  distinguishable:
\[  \lim_{n \rightarrow \infty} \frac{1}{n}  \bfb_{T}^T X_{T}^T (I - P_{X_{\bfgam}}) X_{T} \bfb_{T} = V_{\bfgam}  > 0.
\]
\end{cond}

These conditions are standard assumptions used to establish consistency of Bayesian procedures. \cite{fernandez2001benchmark}, \cite{liang2008mixtures}, \cite{maruyama2011fully} and \cite{bayarri2012criteria} also use these conditions (or slight variations) to demonstrate posterior model selection consistency and prediction consistency of their priors.

\subsection{Model Selection Consistency}\label{sec:modcon}
The second form of consistency we consider is posterior model selection consistency.  A Bayesian model is \emph{model selection consistent} if, when the data $\bfy$ have been generated by model $\Mt$, the posterior probability of model $\Mt$ converges in probability to 1 as the sample size $n \rightarrow \infty$ \citep{fernandez2001benchmark}. The relation between posterior model probabilities and Bayes factors ensures that this consistency criterion can be rephrased as
\begin{eqnarray} \label{eq:cons}
BF(\Mgam : \Mt) \stackrel{P}{\rightarrow} 0 \;  \mbox{ as } n \rightarrow \infty \mbox{ for any model } \Mgam \neq \Mt.
\end{eqnarray}
Criterion \eqref{eq:cons} is precisely the criterion for displaying pairwise consistency of Bayes factors in model selection, which in the fixed-$p$ case coincides with the usual model selection consistency criterion described earlier, also referred to as the strong model selection consistency property.

\begin{thm}\label{thm5.2}
Consider a regression model of the form \eqref{eq:bhg} and \eqref{eq:priorg} which satisfies Conditions~\ref{eq:bhgorth}, \ref{cond5.1} and \ref{cond5.2}. For any model $\Mgam$ such that $\Mgam \supset \Mt$ and all predictors in $\Mgam$ that are not in $\Mt$ are in blocks not in $\Mt$, $BF(\Mgam : \Mt) \stackrel{d}{\rightarrow} W_{\bfgam}$ as $n \rightarrow \infty$ for some non-degenerate random variable $W_{\bfgam}$. For all other models $M_{\bfgam}$, $BF(M_{\bfgam} : \Mt) \stackrel{P}{\rightarrow} 0$. 
\end{thm}

The theorem is proved in Appendix~\ref{appa6}.
The block hyper-$g$ prior is not model selection consistent, as the Bayes factor in \eqref{eq:cons} does not converge to zero in all situations.  This is also the case for the hyper-$g$ prior of \cite{liang2008mixtures}. The defect in these priors is that, as $n \rightarrow \infty$, the priors do not stabilize.  This defect is fixed in the hyper-$g/n$ prior of \cite{liang2008mixtures} and can be fixed in similar fashion here to produce a model selection consistent block hyper-$g/n$ prior \citep{som:14}.

\subsection{Prediction Consistency}
{\em Prediction consistency} concerns the limiting behavior (as $n \rightarrow \infty$) of the Bayes-optimal prediction $\widehat{y}^*_n$ of the true, unknown response $y^*$ for a new vector of predictors $\bfx^* \in \mathbb{R}^p$.  When the true model is not known, the Bayes-optimal prediction under squared-error loss is the Bayesian model averaged prediction
\[
	\widehat{y}^*_n = E(\alpha \mid \bfy_n) + \sum_{\bfgam \in \Gamma} \pi(\Mgam \mid \bfy_n)\, \bfx^{* T} E(\bfb \mid \bfy_n, \Mgam),
\]
where we have subscripted with $n$ the components that depend on the sample size.  \emph{Prediction consistency} is achieved when $\widehat{y}^*_n \stackrel{P}{\rightarrow} E(y^*) = \alpha_T + \bfx^{*T} \bfb_T$ as $n \rightarrow \infty$. The following lemma and its extension are used in the main result on prediction consistency of the block hyper-$g$ prior.

\begin{lem} \label{lem5.1}
Consider a regression model of the form \eqref{eq:bhg} and \eqref{eq:priorg} which satisfies Condition~\ref{eq:bhgorth}. When $\Mt$ is the true model, for any $i \in B_T$
\[ \lim_{n \rightarrow \infty} \int_{(0,1)^{k_T}} \frac{g_i}{1+g_i} \; \pi(\bfg \mid \Mt, \bfy) d\bfg = 1.
\]
\end{lem}
The proof of the lemma is in Appendix \ref{appa7}.  The lemma leads to the following result.

\begin{cor}\label{cor5.1}
Consider a regression model of the form \eqref{eq:bhg} and \eqref{eq:priorg} which satisfies Condition~\ref{eq:bhgorth}. For any model $\Mgam$ containing the true model, i.e., for any $\Mgam \supseteq \Mt$,
\[ \lim_{n \rightarrow \infty} E \left( \frac{g_i}{1+g_i} \mid \Mgam, \bfy \right) = 1,\; \mbox{if } i \in B_T.
\]
\end{cor}

While the block hyper-$g$ prior is not model selection consistent, the next theorem, proved in Appendix~\ref{appa9}, shows that it is prediction consistent under Bayesian model averaging (\textit{BMA}).

\begin{thm}\label{thm5.3}
Consider a regression model of the form \eqref{eq:bhg} and \eqref{eq:priorg} which satisfies Conditions~\ref{eq:bhgorth},  \ref{cond5.1} and \ref{cond5.2}. The predictions under BMA are consistent for this model.    
\end{thm}

\section{Discussion}\label{sec:disc}
We have identified two novel behaviors, \emph{Essentially Least Squares} estimation and the \emph{Conditional Lindley's Paradox}, that are exhibited by many common mixtures of $g$ priors in Bayesian regression.  Both behaviors stem from the use of a single latent scale parameter that is common to all regression coefficients.  We argue that \emph{ELS} behavior is, in general, undesirable, as it precludes shrinkage of small coefficients in the presence of large ones.  Similarly, we argue that priors exhibiting the \emph{CLP} should be avoided as, asymptotically, they can provide infinite evidence in support of a false hypothesis.  Our analyses are driven by a new,  \emph{conditional information asymptotic} that sheds light on a Bayesian linear model's behavior under a strengthening of the likelihood due to one component of the model.  This style of asymptotic is important in practice, as it is often the case that models are comprised of covariates with coefficients of differing magnitude.

With these considerations in mind, we developed a block hyper-$g$ prior as one possible remedy and provided conditions under which the prior does not suffer from either \emph{ELS} or the \emph{CLP}.  Traditional asymptotic analysis of the new prior revealed that, while consistent for prediction under {\em BMA}, model selection consistency is not achieved in all situations.  As stated in Section~\ref{sec:modcon}, the aspect of the prior causing this defect is the failure of the prior to stabilize as $n \rightarrow \infty$.  \cite{som:14} shows that model selection consistency can be achieved by defining a block hyper-$g/n$ prior, where the prior \eqref{eq:priorg} is replaced with
\[
	\pi(\bfg) = \prod_{i=1}^k \frac{a-2}{2n}(1 + g_i / n)^{-a/2}, \;\;\; g_i > 0.
\]
The scaling in this prior on $\bfg$ offsets the scaling in prior \eqref{eq:zgprior} on $\bfb_{\bfgam}$ much like in \cite{liang2008mixtures}.  The block hyper-$g/n$ prior shares many properties of the block hyper-$g$ prior and successfully avoids the new paradoxes in addition to providing consistency in all three aspects described in Section \ref{sec:consistency}. Properties of the block hyper-$g/n$ prior and derivations of new results will be elaborated on elsewhere.

While this paper contains theoretical developments, the new theory is connected to data analysis and raises practical issues.  One particular question is how to best select the groups or blocks of predictors in the design. To date, our data analyses under the block hyper-$g$ priors \citep{som:14} have been based on identifying predictor variables related to one another through a latent or theoretical construct. Predictors measuring the same construct are placed in the same block, as they are reasoned to likely have comparable coefficient sizes.  This suggests a single scale parameter for these related explanatory variables. In the absence of such knowledge, preliminary empirical research suggests that placing correlated predictors in the same block often leads to better performance. The existence of correlated blocks of predictors is sometimes taken to indicate the existence of previously unknown latent constructs. There is scope for theoretical investigation as to why this choice works  or if some other choice can be established as ``optimal" under specific settings. 

Finally, the results presented in Sections~\ref{sec:newprior} and \ref{sec:consistency} required block orthogonality of the predictor variables (Condition~\ref{eq:bhgorth}). We note that this condition can be relaxed to obtain similar results under a modified version of the \emph{conditional information asymptotic} under suitable adjustments to the original formulation of a block $g$ prior  \citep{som:14}.  Results along these lines will be reported elsewhere.

\appendix

\section{Appendix: Proofs of the Main Results}\label{app:proofs}

\subsection{Proof of Theorem \ref{thm3.1}} \label{appa1}

The posterior mean of the regression coefficients is
\[ \widehat{\bfb} = E \left(\frac{g}{1+g} \mid \bfy \right) \widehat{\bfb}_{LS}.
\]
For the hyper-$g$ prior, the posterior expectation of the shrinkage factor can be expressed in terms of $R^2$ \citep[see][]{liang2008mixtures}:

\[ E \left(\frac{g}{1+g} \mid \bfy \right) = \frac{2}{p+a} \frac{_2F_1(\frac{n-1}{2},2;\frac{p+a}{2}+1;R^2)} {_2F_1(\frac{n-1}{2},1;\frac{p+a}{2};R^2)} ,
\]
where $_2F_1$ is the Gaussian Hypergeometric Function. $_2F_1(a,b;c;z)$ is finite for $|z| < 1 \;\mbox{whenever} \; c>b>0$. Here, $c-b=(p+a)/2 - 1 > 0$  since $2< a \leq 4$ and $p>0$. Thus, for all values of $R^2 < 1$, both numerator and denominator are finite. We use an integral representation of the $_2 F_1$ function:

\begin{eqnarray*}
 E \left(\frac{g}{1+g} \mid \bfy \right) &=& \frac{2}{p+a} \frac{_2F_1(\frac{n-1}{2},2;\frac{p+a}{2}+1;R^2)} {_2F_1(\frac{n-1}{2},1;\frac{p+a}{2};R^2)} \\
 &=& \frac{ \int_0^1 t(1-t)^{\frac{p+a}{2}-2} (1-tR^2)^{-\frac{n-1}{2}} dt} { \int_0^1 (1-t)^{\frac{p+a}{2}-2} (1-tR^2)^{-\frac{n-1}{2}} dt}. \\
\end{eqnarray*}
Define $m=\frac{n-1}{2}$, $b=\frac{p+a}{2}-2$ and $z=R^2 (\leq 1)$ so that we have
\[ E \left(\frac{g}{1+g} \mid \bfy \right) = \frac{ \int_0^1 t(1-t)^b (1-tz)^{-m} dt} { \int_0^1 (1-t)^b (1-tz)^{-m} dt},
\]
where
\begin{eqnarray*}
\mbox{Numerator} &=& \int_0^1 t(1-t)^b \left[ \sum_{k=0}^{\infty} {m + k-1 \choose k} (tz)^k \right] dt \\
 &=& \sum_{k=0}^{\infty} \frac{(k+1) \Gamma(m+k) \Gamma(b+1)}{\Gamma(m) \Gamma(b+k+3)} z^k 
\end{eqnarray*}

\[  \mbox{and Denominator} = \int_0^1 (1-t)^b \left[ \sum_{k=0}^{\infty} {m + k-1 \choose k} (tz)^k \right] dt = \sum_{k=0}^{\infty} \frac{ \Gamma(m+k) \Gamma(b+1)}{\Gamma(m) \Gamma(b+k+2)} z^k.  \]
Thus,
\[E \left(\frac{g}{1+g} \mid \bfy \right) = \frac{\sum_{k=0}^{\infty} \frac{(k+1) \Gamma(m+k) \Gamma(b+1)}{\Gamma(m) \Gamma(b+k+3)} z^k} {\sum_{k=0}^{\infty} \frac{\Gamma(m+k) \Gamma(b+1)}{\Gamma(m) \Gamma(b+k+2)} z^k} 
=  \frac{ \sum_{k=0}^{\infty} \frac{\Gamma(m+k)}{\Gamma(b+k+2)} \frac{1}{1+\frac{b+1}{k+1}} z^k  }{ \sum_{k=0}^{\infty} \frac{\Gamma(m+k)}{\Gamma(b+k+2)}z^k }.  \]

\noindent
When $m=\frac{n-1}{2} > b+2 = \frac{p+a}{2}$ , the ratio $\frac{\Gamma(m+k)}{\Gamma(b+k+2)}$ is increasing in $k$. 

Lemma \ref{lem3.1} states that $R^2 \rightarrow 1$ as $N \rightarrow \infty$,  so the proof will be complete if we can show that $\lim\limits_{R^2 \rightarrow 1} E \left(\frac{g}{1+g} \mid \bfy \right) = 1$. 
\begin{eqnarray*}
\lim_{R^2 \rightarrow 1} E \left(\frac{g}{1+g} \mid \bfy \right) &=& \lim_{z\uparrow 1} \frac{ \sum_{k=0}^{\infty} \frac{\Gamma(m+k)}{\Gamma(b+k+2)} \frac{1}{1+\frac{b+1}{k+1}} z^k  }{ \sum_{k=0}^{\infty} \frac{\Gamma(m+k)}{\Gamma(b+k+2)}z^k } \\
\end{eqnarray*}

\noindent
\underline{Case 1}: $n > p+a+1$ 

\noindent
First note that $\frac{1}{1 + \frac{b+1}{k+1}}$ is increasing in $k$ and $\uparrow 1$ as $k \rightarrow \infty$. So for any $\eta >0, \; \exists \; N_0$ such that $\forall \; k > N_0, \frac{1}{1 + \frac{b+1}{k+1}} > 1-\eta$. 

Hence, for any $\eta > 0$, $ \lim\limits_{z \uparrow 1} \frac{ \sum_{k=0}^{\infty} \frac{\Gamma(m+k)}{\Gamma(b+k+2)} \frac{1}{1+\frac{b+1}{k+1}} z^k  }{ \sum_{k=0}^{\infty} \frac{\Gamma(m+k)}{\Gamma(b+k+2)}z^k } $
\begin{eqnarray} \label{eq:a1denom}
&=& \lim_{z \uparrow 1} \frac{ \sum_{k=0}^{N_0} \frac{\Gamma(m+k)}{\Gamma(b+k+2)} \frac{1}{1+\frac{b+1}{k+1}} z^k +  \sum_{k=N_0+1}^{\infty} \frac{\Gamma(m+k)}{\Gamma(b+k+2)} \frac{1}{1+\frac{b+1}{k+1}} z^k }{ \sum_{k=0}^{N_0} \frac{\Gamma(m+k)}{\Gamma(b+k+2)}z^k + \sum_{k=N_0+1}^{\infty} \frac{\Gamma(m+k)}{\Gamma(b+k+2)}z^k } \nonumber \\
&>& \lim_{z\uparrow 1} \frac{q_1 + (1-\eta)\sum_{k=N_0+1}^{\infty} \frac{\Gamma(m+k)}{\Gamma(b+k+2)}z^{k} }{q_2 + \sum_{k=N_0+1}^{\infty} \frac{\Gamma(m+k)}{\Gamma(b+k+2)}z^{k} } \nonumber \\
&=& (1-\eta) + \lim_{z \uparrow 1} \frac{q_1 - (1-\eta)q_2}{q_2 + \sum_{k=N_0+1}^{\infty} \frac{\Gamma(m+k)}{\Gamma(b+k+2)}z^{k}} \\
&\geq& (1-\eta) + 0 = (1-\eta)  \nonumber
\end{eqnarray}
\begin{equation*}
\implies \lim_{R^2 \rightarrow 1} E \left(\frac{g}{1+g} \mid \bfy \right) = \lim_{z\uparrow 1} \frac{ \sum_{k=0}^{\infty} \frac{\Gamma(m+k)}{\Gamma(b+k+2)} \frac{1}{1+\frac{b+1}{k+1}} z^k  }{ \sum_{k=0}^{\infty} \frac{\Gamma(m+k)}{\Gamma(b+k+2)}z^k } = 1 \hspace{2 in}
\end{equation*}
Note that $q_1$ and $q_2$ are finite numbers corresponding to the finite sums of the first $ N_0$ terms. Also $\frac{\Gamma(m+k)}{\Gamma(b+k+2)} \rightarrow \infty$ as $k \rightarrow \infty$ due to which $\sum_{k=0}^{\infty} \frac{\Gamma(m+k)}{\Gamma(b+k+2)}$ and hence the denominator in \eqref{eq:a1denom} goes to infinity causing the second term above to vanish in the limit. 

\noindent
\underline{Case 2}: $p+a-1 \leq n \leq p+a+1$ 

\noindent
Let $n=p+a-1+ 2\xi$, where $0 \leq \xi \leq 1$ and define $N_0$ as in Case 1,
\begin{eqnarray*}
\lim_{R^2 \rightarrow 1} E \left(\frac{g}{1+g} \mid \bfy \right) &=& \lim_{z \uparrow 1} \frac{ \sum_{k=0}^{\infty} \frac{\Gamma(b+k+1+\xi)}{\Gamma(b+k+1)(b+k+1)} \frac{1}{1+\frac{b+1}{k+1}} z^k}{ \sum_{k=0}^{\infty} \frac{\Gamma(b+k+1+\xi)}{\Gamma(b+k+1)(b+k+1)} z^k } 
\end{eqnarray*}

\noindent
Proceeding as in Case 1, we can show that
\begin{eqnarray*}
\lim_{R^2 \rightarrow 1} E \left(\frac{g}{1+g} \mid \bfy \right) &>& (1-\eta) + \lim_{z \uparrow 1} \frac{q_1 - (1-\eta)q_2}{q_2 + \sum_{k=N_0+1}^{\infty} \frac{\Gamma(b+k+1+\xi)}{\Gamma(b+k+1)(b+k+1)}z^{k}} \\
&\geq& (1-\eta) \;,\; \mbox{for any} \; \eta>0 
\end{eqnarray*}

\noindent
As $z \uparrow 1$, the denominator of the second term becomes $\frac{q_1 - (1-\eta)q_2}{q_2 + \sum_{k=N_0+1}^{\infty} \frac{\Gamma(b+k+1+\xi)}{\Gamma(b+k+1)(b+k+1)} }$ which tends to zero if the infinite sum $\sum_{k=N_0+1}^{\infty} \frac{\Gamma(b+k+1+\xi)}{\Gamma(b+k+1)(b+k+1)}$ diverges. It does, since $\frac{\Gamma(b+k+1+\xi)}{\Gamma(b+k+1)(b+k+1)} = O(k^{-\lambda})$, where $0 \leq \lambda \leq 1$ and $ \sum_{k=N_0+1}^{\infty} O(k^{-\lambda}) = \infty $ for any $0 \leq \lambda \leq 1$. Thus, 
\[  \lim\limits_{R^2 \rightarrow 1} E \left(\frac{g}{1+g} \mid \bfy \right) = 1. \]

A more detailed version of this proof can be found in Appendix~\ref{appb11} in the supplementary materials.

\subsection{Proof of Theorem \ref{thm3.2}}\label{appa2}

\cite{liang2008mixtures} show that
\begin{eqnarray*}
 BF(M_2:\Mo) &=&  \;_2F_1 (\frac{n-1}{2},1;\frac{a+p}{2}; R_{M_2}^2) \times \frac{a-2}{a+p-2} \\
 BF(M_1:\Mo) &=& \;_2F_1 (\frac{n-1}{2},1;\frac{a+p_1}{2}; R_{M_1}^2) \times \frac{a-2}{a+p_1-2}
\end{eqnarray*}
where $R_{M_i}^2$ is the coefficient of determination for model $M_i, i=1,2$.

The $g$ prior is invariant to linear transformation of $X$, and so we can work with an orthogonalized version of the design without loss of generality.  Specifically, we consider $Q_1 = X_1$ and $Q_2 = (I - P_{Q_1})X_2$, where $P_{Q_1}$ is the projection matrix for the column space of $Q_1$.  Then $X$ can be represented as $X=QT$ for a suitable upper triangular matrix $T$ and $X \bfb = Q \kappa$, where $\kappa=T\bfb$ also has a hyper-$g$ prior. Since $T$ is upper triangular, $||\bfb_1|| \rightarrow \infty$ is equivalent to $||\kappa_1|| \rightarrow \infty$ while $\kappa_2$ stays fixed in the sequence.  Under the block orthogonal setup, $ R_{M_1}^2 = \frac{(Q_1 \widehat{\bld{\kappa}}_1)^T (Q_1 \widehat{\bld{\kappa}}_1)}{\bld{y}^T \bld{y}}$ and $R_{M_2}^2 = R_{M_1}^2 + \frac{(Q_2 \widehat{\bld{\kappa}}_2)^T (Q_2 \widehat{\bld{\kappa}}_2)}{\bld{y}^T \bld{y}}$.  The term $(Q_2 \widehat{\bld{\kappa}}_2)^T (Q_2 \widehat{\bld{\kappa}}_2)$ is constant throughout the sequence $\{ \Psi_N \}$ and so $\frac{(Q_2 \widehat{\bld{\kappa}}_2)^T (Q_2 \widehat{\bld{\kappa}}_2)}{ \bld{y}^T \bld{y} } \rightarrow 0$ as $N \rightarrow \infty$.
\begin{eqnarray*}
BF(M_2:M_1) &=& \frac{BF(M_2:\Mo)}{BF(M_1:\Mo)} \\
&=& \frac{a+p_1-2}{a+p-2} .\; \frac{_2F_1 (\frac{n-1}{2},1;\frac{a+p}{2}; R_{M_2}^2)}{_2F_1 (\frac{n-1}{2},1;\frac{a+p_1}{2}; R_{M_1}^2)} \\
&=& \frac{\int_0^1 (1-t)^{\frac{a+p}{2}-2} (1-tR_{M_2}^2)^{-\frac{n-1}{2}} dt}{ \int_0^1 (1-t)^{\frac{a+p_1}{2}-2} (1-tR_{M_1}^2)^{-\frac{n-1}{2}}dt }
\end{eqnarray*}

\noindent
Define $b = \frac{a+p_1}{2}-2$, $m=\frac{n-1}{2}$, $R_{M_1}^2 =z$ and $R_{M_2}^2 =z+q$. When $||\beta_1|| \rightarrow \infty$, both $R_{M_2}^2$ and $R_{M_1}^2$ go to 1 which results in $z \uparrow 1$ and $q \downarrow 0$.

\[ BF(M_2:M_1) = \frac{\int_0^1 (1-t)^{b+\frac{p_2}{2}} \left[ 1-t(z+q) \right]^{-m} dt }{ \int_0^1 (1-t)^{b} \left[ 1-tz \right]^{-m} dt }
\]

\noindent
Proceeding as in Theorem \ref{thm3.1}, 
\[ \mbox{Numerator} = \sum_{k=0}^{\infty} \frac{\Gamma(m+k) \Gamma(b+1+\frac{p_2}{2})}{\Gamma(m) \Gamma(b+k+2+\frac{p_2}{2})} (z+q)^k
\]
and 
\[  \mbox{Denominator} = \sum_{k=0}^{\infty} \frac{\Gamma(m+k) \Gamma(b+1)}{\Gamma(m) \Gamma(b+k+2)} z^k
\]
Thus,
\[ BF(M_2:M_1) = \frac{\Gamma(b+1+\frac{p_2}{2})}{\Gamma(b+1)} \; \frac{ \sum_{k=0}^{\infty} \frac{\Gamma(m+k)}{\Gamma(b+k+2)} \left\{ \frac{\Gamma(b+k+2)}{\Gamma(b+k+2+ \frac{p_2}{2})} \right\}(z+q)^k }{ \sum_{k=0}^{\infty} \frac{\Gamma(m+k)}{\Gamma(b+k+2)} z^k }
\]
Hence,
\begin{eqnarray*}
\lim_{||\beta_1|| \rightarrow \infty} BF(M_2:M_1) &=& \lim_{\substack{z \rightarrow 1 \\ q \rightarrow 0}} BF(M_2:M_1) \\
&=& \lim_{z \rightarrow 1} \; \left\{ \lim_{q \rightarrow 0} BF(M_2:M_1) \right\}
\end{eqnarray*}

\noindent
The last step is justified when $ \lim\limits_{q \rightarrow 0} BF(M_2:M_1)$ exists for all $0\leq z < 1$. This holds since  $\lim \limits_{q \rightarrow 0} BF(M_2:M_1) = \frac{a+p_1-2}{a+p-2} \; \frac{_2F_1(m,1;b+2+\frac{p_2}{2};z)}{_2F_1(m,1;b+2;z)} $ which exists and is finite for all $0 \leq z < 1$.
\begin{eqnarray*}
\lim_{||\beta_1|| \rightarrow \infty} BF(M_2:M_1) &=& \lim_{z \uparrow 1} \frac{\Gamma(b+1+\frac{p_2}{2})}{\Gamma(b+1)} \; \frac{ \sum_{k=0}^{\infty} \frac{\Gamma(m+k)}{\Gamma(b+k+2)} \left\{ \frac{\Gamma(b+k+2)}{\Gamma(b+k+2 +\frac{p_2}{2})} \right\}z^k }{ \sum_{k=0}^{\infty} \frac{\Gamma(m+k)}{\Gamma(b+k+2)} z^k } 
\end{eqnarray*}

\noindent
But $\frac{\Gamma(b+k+2)}{\Gamma(b+k+2+ \frac{p_2}{2})}$ decreases to 0 as $k \rightarrow \infty$ (see Appendix \ref{appb11} for a proof). Hence given an arbitrary $\eta>0$, we can find a number $N_0$ such that $\forall \; k > N_0$, $\frac{\Gamma(b+k+2)}{\Gamma(b+k+2+ \frac{p_2}{2})} < \eta $.

\begin{eqnarray*}
&& BF(M_2:M_1) = \frac{\Gamma(b+1+\frac{p_2}{2})}{\Gamma(b+1)} \times \\
 &&  \frac{ \sum_{k=0}^{N_0} \frac{\Gamma(m+k)}{\Gamma(b+k+2)} \left\{ \frac{\Gamma(b+k+2)}{\Gamma(b+k+2+ \frac{p_2}{2})} \right\}z^k + \sum_{k=N_0+1}^{\infty} \frac{\Gamma(m+k)}{\Gamma(b+k+2)} \left\{ \frac{\Gamma(b+k+2)}{\Gamma(b+k+2+ \frac{p_2}{2})} \right\}z^k }{ \sum_{k=0}^{N_0} \frac{\Gamma(m+k)}{\Gamma(b+k+2)} z^k + \sum_{k=N_0+1}^{\infty} \frac{\Gamma(m+k)}{\Gamma(b+k+2)} z^k } \\
&<& \frac{q_1 + \eta \sum_{k=N_0+1}^{\infty} \frac{\Gamma(m+k)}{\Gamma(b+k+2)} z^k }{q_2 + \sum_{k=N_0+1}^{\infty} \frac{\Gamma(m+k)}{\Gamma(b+k+2)} z^k} \times \frac{\Gamma(b+1+\frac{p_2}{2})}{\Gamma(b+1)} \\
&=& \frac{\Gamma(b+1+\frac{p_2}{2})}{\Gamma(b+1)} \; . \frac{q_1+\eta T}{q_2 + T} \;,\; \mbox{with } T = \sum_{k=N_0+1}^{\infty} \frac{\Gamma(m+k)}{\Gamma(b+k+2)} z^k  \\
&=& \frac{\Gamma(b+1+\frac{p_2}{2})}{\Gamma(b+1)} \left[ \frac{\eta(q_1+T)}{q_2+T} + \frac{(1-\eta)q_1}{q_2+T} \right]
\end{eqnarray*}

\noindent
We later show that $||\beta_1|| \rightarrow \infty$ (or $z \uparrow 1$) implies that $T \rightarrow \infty$ when $n \geq a+p_1-1$.
\begin{eqnarray*}
\implies \lim_{||\beta_1|| \rightarrow \infty} BF(M_2:M_1) &\leq& \lim_{T \rightarrow \infty} \frac{\Gamma(b+1+\frac{p_2}{2})}{\Gamma(b+1)} \left[ \frac{\eta(q_1+T)}{q_2+T} + \frac{(1-\eta)q_1}{q_2+T} \right] \\
&=& \eta \; \frac{\Gamma(b+1+\frac{p_2}{2})}{\Gamma(b+1)} \\
\mbox{Hence }  \lim_{||\beta_1|| \rightarrow \infty} BF(M_2:M_1) &=& 0
\end{eqnarray*}

\noindent
We now prove that, for $n \geq a+p_1-1$, $T \rightarrow \infty$ when $||\beta_1|| \rightarrow \infty$.

\noindent
\underline{Case 1}: $n > a+p_1+1$

\noindent
Then $m>b+2$ and so $\frac{\Gamma(m+k)}{\Gamma(b+k+2)} \uparrow \infty$ as $k \rightarrow \infty$. 
\[ \lim_{z \uparrow 1} T = \sum_{k=N_0+1}^{\infty} \frac{\Gamma(m+k)}{\Gamma(b+k+2)} = \infty
\]

\noindent
\underline{Case 2}: $a+p_1-1 \leq n \leq a+p_1+1$

\noindent
Let $n=a+p_1-1+2\xi$ where $0 \leq \xi \leq 1$. Then $m=\frac{a+p_1}{2} -1 + \xi$ and $\frac{\Gamma(m+k)}{\Gamma(b+k+2)} = \frac{\Gamma(\frac{a+p_1}{2}-1+\xi+k)}{\Gamma(\frac{a+p_1}{2}+k)}$. \\
For $0 \leq \xi \leq 1$, $\frac{\Gamma(\frac{a+p_1}{2}-1+\xi+k)}{\Gamma(\frac{a+p_1}{2}+k)} =\frac{\Gamma(\frac{a+p_1}{2}-1+k+\xi)}{\Gamma(\frac{a+p_1}{2}-1+k+1)} = O(k^{-\lambda})$ for some $0 \leq\lambda \leq 1$. But $\sum_{k=N_0+1}^{\infty} O(k^{-\lambda}) = \infty$ for such values of $\lambda$ implying that $T \rightarrow \infty$ as $z \uparrow 1$.

A more detailed version of this proof can be found in Appendix~\ref{appDetail3.2} in the supplementary materials.

\subsection{A Preliminary Lemma}

The following lemma, proved in Appendix \ref{appa3}, is useful to derive the next set of theoretical results.

\begin{lem} \label{lem4.2}
If $f_1(t_m)$ and $f_2(t_m)$ denote properly normalized pdfs on $(0,1)$ with
\begin{eqnarray*}
f_1(t_m) &\propto& \int_{(0,1)^{k-1}} \left[ \prod_{i=1}^k (1-t_i)^{\frac{a+p_i}{2}-2} \right] (1-\sum_{i=1}^k t_i R_i^2)^{-\frac{n-1}{2}} d\bft_{-m} \\
\mbox{ and } f_2(t_m) &\propto& \int_{(0,1)^{k-1}} \left[ \prod_{i=1}^k (1-t_i)^{\frac{a+p_i}{2}-2} \right] (1- t_j R_j^2)^{-\frac{n-1}{2}} d\bft_{-m}
\end{eqnarray*}
for some $ m , j \in \{1,2,...,k \} $, where $j$ may or may not equal $m$ and $\bft_{-m}=\{ t_i : i \neq m\}$, then $E_{f_1} (t_{m}) \geq E_{f_2} (t_m)$. Strict inequality holds when $R_m^2 > 0$ and $R_i^2 > 0$, for at least one $i \neq m$.
\end{lem}

\subsection{ Proof of Theorem \ref{thm4.1}} \label{appa4}

The first part of the proof is trivial and follows directly from \eqref{eq:postmean}. 

In the block orthogonal setup
\[ \pi(\bfg \mid \bfy) \propto \frac{\prod_{j=1}^k (1+g_j)^{-\frac{a+p_j}{2}} }{\left[ 1 -  \sum_{j=1}^k \frac{g_j}{g_j+1} R_j^2 \right]^{(n-1)/2} } \]

\noindent
So for any $i=1,2,..,k$,
\begin{eqnarray*}
E \left(\frac{g_i}{1+g_i} \mid \bfy \right) &=& \frac{ \int_{(0,1)^k} t_i \prod_{j=1}^k (1-t_j)^{\frac{a+p_j}{2}-2} (1-\sum_{j=1}^k t_j R_j^2)^{-\frac{n-1}{2}} d\bft}{ \int_{(0,1)^k} \prod_{j=1}^k (1-t_j)^{\frac{a+p_j}{2}-2} (1-\sum_{j=1}^k t_j R_j^2)^{-\frac{n-1}{2}} d\bft } \\
&\geq& \frac{ \int_{(0,1)^k}  t_i  \prod_{j=1}^k (1-t_j)^{\frac{a+p_j}{2}-2} (1- t_i R_i^2)^{-\frac{n-1}{2}} d\bft}{ \int_{(0,1)^k} \prod_{j=1}^k (1-t_j)^{\frac{a+p_j}{2}-2} (1- t_i R_i^2)^{-\frac{n-1}{2}} d\bft } \\
&& \hspace{1.4 in} (\mbox{by Lemma }\ref{lem4.2}) \\
&=& \frac{2}{a+p_i} \; \frac{_2F_1(\frac{n-1}{2},2;\frac{a+p_i}{2}+1;R_i^2)}{_2F_1(\frac{n-1}{2},1;\frac{a+p_i}{2};R_i^2)} 
\end{eqnarray*} 

\[  \mbox{Hence } \lim_{N \rightarrow \infty} E \left(\frac{g_i}{1+g_i} \mid \bfy \right) \geq  \lim_{N \rightarrow \infty}  \frac{2}{a+p_i} \; \frac{_2F_1(\frac{n-1}{2},2;\frac{a+p_i}{2}+1;R_i^2)}{_2F_1(\frac{n-1}{2},1;\frac{a+p_i}{2};R_i^2)} \]

\noindent
As $N \rightarrow \infty$, $R^2_1 \rightarrow 1$ so that for $i=1$,
\begin{eqnarray*}
\lim_{N \rightarrow \infty} E \left(\frac{g_1}{1+g_1} \mid \bfy \right) &\geq&  \lim_{z \rightarrow 1}  \frac{2}{a+p_1} \; \frac{_2F_1(\frac{n-1}{2},2;\frac{a+p_1}{2}+1;z)}{_2F_1(\frac{n-1}{2},1;\frac{a+p_1}{2};z)}  \\
&=& 1 \; ,\; \mbox{ when } n \geq a+p_1-1\mbox{ (see Theorem \ref{thm3.1}) }  
\end{eqnarray*}
But $E \left(\frac{g_1}{1+g_1} \mid \bfy \right) \leq 1$, implying that $E \left(\frac{g_1}{1+g_1} \mid \bfy \right) \rightarrow 1$ in the limit. 

For $i > 1$ and  $m \neq i$, $ E \left(\frac{g_i}{1+g_i} \mid \bfy \right)$
\begin{eqnarray*} 
&=&  \frac{ \int_{(0,1)^k} t_i \prod_{j=1}^k (1-t_j)^{\frac{a+p_j}{2}-2} (1-\sum_{j=1}^k t_j R_j^2)^{-\frac{n-1}{2}} d\bft}{ \int_{(0,1)^k} \prod_{j=1}^k (1-t_j)^{\frac{a+p_j}{2}-2} (1-\sum_{j=1}^k t_j R_j^2)^{-\frac{n-1}{2}} d\bft }  \\
&\geq& \frac{ \int_{(0,1)^k} t_i \prod_{j=1}^k (1-t_j)^{\frac{a+p_j}{2}-2} (1- t_m R_m^2)^{-\frac{n-1}{2}} d\bft}{ \int_{(0,1)^k} \prod_{j =1}^{k} (1-t_j)^{\frac{a+p_j}{2}-2} (1- t_m R_m^2)^{-\frac{n-1}{2}} d\bft } \; \; (\mbox{by Lemma }\ref{lem4.2}) \\
&=& \frac{\mbox{Beta}(2,\frac{a+p_i}{2}-1)}{\mbox{Beta}(1,\frac{a+p_i}{2}-1)} \times \frac{_2F_1(\frac{n-1}{2},1;\frac{a+p_m}{2};R^2_m)}{_2F_1(\frac{n-1}{2},1;\frac{a+p_m}{2};R^2_m)} = \frac{2}{a+p_i} 
\end{eqnarray*}

Thus, $\lim\limits_{N \rightarrow \infty} E \left(\frac{g_i}{1+g_i} \mid \bfy \right) \geq \frac{2}{a+p_i}$. 

\begin{eqnarray*}
E \left(\frac{g_i}{1+g_i} \mid \bfy \right) &=& \frac{ \int_{(0,1)^k} t_i \prod_{j=1}^k (1-t_j)^{\frac{a+p_j}{2}-2} (1-\sum_{j=1}^k t_j R_j^2)^{-\frac{n-1}{2}} d\bft}{ \int_{(0,1)^k} \prod_{j=1}^k (1-t_j)^{\frac{a+p_j}{2}-2} (1-\sum_{j=1}^k t_j R_j^2)^{-\frac{n-1}{2}} d\bft }  \\
&\leq& \frac{ \int_{(0,1)^k} t_i \prod_{j=1}^k (1-t_j)^{\frac{a+p_j}{2}-2} (1-\sum_{j \neq i} R_j^2 - t_i R_i^2)^{-\frac{n-1}{2}} d\bft}{ \int_{(0,1)^k} \prod_{j=1}^k (1-t_j)^{\frac{a+p_j}{2}-2} (1-\sum_{j \neq i} R_j^2 - t_i R_i^2)^{-\frac{n-1}{2}} d\bft } \\
&& \hspace{1 in} (\mbox{using a variation of Lemma } \ref{lem4.2}) 
\end{eqnarray*}

\begin{eqnarray*}
&=& \frac{ (1-\sum_{j \neq i} R_j^2)^{-(n-1)/2} }{ (1-\sum_{j \neq i} R_j^2)^{-(n-1)/2} } \;  \frac{ \int_0^1 t_i (1-t_i)^{\frac{a+p_j}{2}-2} (1- t_i \frac{R_i^2}{1-\sum_{j \neq i} R_j^2})^{-\frac{n-1}{2}} dt_i}{ \int_0^1 (1-t_i)^{\frac{a+p_j}{2}-2} (1 - t_i \frac{R_i^2}{1-\sum_{j \neq i} R_j^2})^{-\frac{n-1}{2}} dt_i } \\
&=& \frac{2}{a+p_i} \frac{_2F_1(\frac{n-1}{2},2;\frac{a+p_i}{2}+1;\kappa_i)}{_2F_1(\frac{n-1}{2},1;\frac{a+p_i}{2};\kappa_i)}  < 1
\end{eqnarray*}
\noindent
where $\kappa_i = \frac{R_i^2}{1-\sum_{j \neq i} R_j^2}$.

\noindent
For this sequence of problems, $0<\kappa_i<1$ is fixed for all $i \neq 1$, because
\[ \kappa_i = \frac{R_i^2}{1-\sum_{j \neq i} R_j^2} =\frac{\bfy^T P_{X_i} \bfy}{(n-p-1)\hat{\sigma}^2+\bfy^T P_{X_i} \bfy} 
\]
and so
\[ \lim_{N \rightarrow \infty} E \left(\frac{g_i}{1+g_i} \mid \bfy \right) < 1, \; \mbox{for} \; i\neq 1.
\]

A more detailed version of this proof can be found in Appendix~\ref{appDetail4.1} in the supplementary materials.

\subsection{Proof of Theorem \ref{thm4.2}}\label{appa45}

The Bayes factors $BF(M_i:\Mo)$ comparing the models $M_i, i=1,2$ to the null (intercept only) model are
\begin{eqnarray*}
BF(M_2:\Mo) &=& \left(\frac{a-2}{2} \right)^2 \int_{0}^1 \int_{0}^1 \prod_{i=1}^2 (1-t_i)^{\frac{a+p_i}{2}-2} (1-\sum_{i=1}^2 t_i R_i^2)^{-\frac{n-1}{2}} dt_1 dt_2 \\
\mbox{and } BF(M_1:\Mo) &=& \frac{a-2}{2} \int_{0}^1 (1-t_1)^{\frac{a+p_1}{2}-2} (1- t_1 R_1^2)^{-\frac{n-1}{2}} dt_1 
\end{eqnarray*}
\begin{eqnarray*}
\mbox{Thus, } BF(M_2:M_1) &=& \frac{BF(M_2:\Mo)}{BF(M_1:\Mo)}   \\
&=& \frac{a-2}{2} \frac{ \int_{0}^1 \int_{0}^1 \prod_{i=1}^2 (1-t_i)^{\frac{a+p_i}{2}-2} (1-\sum_{i=1}^2 t_i R_i^2)^{-\frac{n-1}{2}}  dt_1 dt_2 }{\int_{0}^1 (1-t_1)^{\frac{a+p_1}{2}-2} (1- t_1 R_1^2)^{-\frac{n-1}{2}} dt_1 } \hspace{1 cm} \\
&\geq& \frac{a-2}{2} \frac{\int_{0}^1 \int_{0}^1 \prod_{i=1}^2 \left[ (1-t_i)^{\frac{a+p_i}{2}-2} (1-t_i R_i^2)^{-\frac{n-1}{2}} \right] dt_1 dt_2}{ \int_0^1 (1-t_1)^{\frac{a+p_1}{2}-2} (1- t_1 R_1^2)^{-\frac{n-1}{2}} dt_1 } \\
&=& \frac{a-2}{2} \int_0^1 (1-t_2)^{\frac{a+p_2}{2}-2} (1- t_2 R_2^2)^{-\frac{n-1}{2}} dt_2
\end{eqnarray*}

As $||\beta_1|| \rightarrow \infty$, $R_1^2 \rightarrow 1$ and $R_2^2 \rightarrow 0$, and so
\[ \lim_{||\beta_1|| \rightarrow \infty} BF(M_2:M_1) \geq \frac{a-2}{2} \int_0^1 (1-t_2)^{\frac{a+p_2}{2}-2} dt_2 = \frac{a-2}{a+p_2-2}
\]

A more detailed version of this proof can be found in Appendix~\ref{appb45} in the supplementary materials.

\subsection{Proof of Theorem \ref{thm5.1}}\label{appa5}

Information consistency under model $\Mgam$ is equivalent to $BF(\Mgam:\Mo) \rightarrow \infty$ as  $R^2_{\bfgam}\rightarrow 1$. Dropping the subscript $\bfgam$ for convenience, we first establish the sufficiency of Condition (2). We know $R^2=\sum_{j=1}^k R^2_j$ and Condition (2) of the theorem enforces $R^2_i \rightarrow 1$ for a given block $i$, implying  $R^2_{j} \rightarrow 0 \;\; \forall \; j \neq i$. Then
\begin{eqnarray*}
 BF(\Mgam:\Mo) &=& \left( \frac{a-2}{2} \right)^{k} \int_{(0,1)^{k}} \left[ \prod_{j=1}^{k} (1-t_j)^{\frac{a+p_j}{2}-2} \right] \left(1-\sum_{j=1}^{k} t_j R_{j}^2 \right)^{-\frac{n-1}{2}} d\bft \\
&\geq& \left( \frac{a-2}{2} \right)^{k} \prod_{j=1}^{k}  \left[ \int_{0}^1 (1-t_j)^{\frac{a+p_j}{2}-2} \left(1- t_j R_{j}^2 \right)^{-\frac{n-1}{2}} dt_j \right]  \\
\implies \lim\limits_{R^2_{i} \rightarrow 1} BF(\Mgam:\Mo) &\geq& \lim\limits_{R^2_i \rightarrow 1} \left( \frac{a-2}{2} \right)^{k} \int_{0}^1 (1-t_i)^{\frac{a+p_i}{2}-2} \left(1- t_i R_{i}^2 \right)^{-\frac{n-1}{2}} dt_i \\
& & \hspace{1.2 in} \times \prod_{j\neq i}  \left[ \int_0^1 (1-t_j)^{\frac{a+p_{j}}{2}-2} dt_j \right] 
\end{eqnarray*}
 
\noindent
The first term on the RHS goes to $\infty$ when $n \geq a+p_{i}-1$ while the rest of the terms converge to nonzero constants. Hence the Bayes factor also diverges in the limit. 

Under Condition (1),  the block structure does not play any role in driving consistency. For an arbitrary $0 < \eta < 1$, 
\begin{eqnarray*}
 BF(\Mgam:\Mo) &=& \left( \frac{a-2}{2} \right)^{k} \int_{(0,1)^k} \left[ \prod_{j=1}^{k} (1-t_j)^{\frac{a+p_j}{2}-2} \right] \left(1-\sum_{j=1}^{k} t_j R_{j}^2 \right)^{-\frac{n-1}{2}} d\bft \\
&>&  \left( \frac{a-2}{2} \right)^k \int_{(1-\eta,1)^k} \left[ \prod_{j=1}^{k} (1-t_j)^{\frac{a+p_j}{2}-2} \right] (1- (1- \eta) R^2 )^{-\frac{n-1}{2}} d\bft 
\end{eqnarray*}
Hence, by the Monotone Convergence Theorem, 
\begin{eqnarray*}
\lim\limits_{R^2 \rightarrow 1} BF(\Mgam:\Mo)  &\geq& \eta^{-\frac{n-1}{2}} \left( \frac{a-2}{2} \right)^{k} \prod_{j=1}^{k} \left[  \int_{1-\eta}^1 (1-t_j )^{\frac{a+p_j}{2}-2} dt_j \right]  \\ 
& = & \eta^{\frac{(a-2) k +p - n +1}{2} } \prod_{j=1}^{k} \left[ \frac{a-2}{a+p_j-2} \right]
\end{eqnarray*}
where $p = \sum_{j=1}^{k} p_{j}$. When $n > k(a-2)+p+1$, the exponent of $\eta$ is negative, indicating that $\lim\limits_{R^2 \rightarrow 1} BF(\Mgam:\Mo)  = \infty$.

\subsection{ Proof of Theorem \ref{thm5.2}} \label{appa6}

Let $R_{i,T}^2$ and $p_{i,T}$ represent the component of $R^2$ and the number of predictors in the $i^{th}$ block of the true model $\Mt$ while $R_{i,\gamma}^2$ and $p_{i,\gamma}$ denote the corresponding entities for model $M_{\bfgam}$. Further assume $B_{\gamma}$ denotes the set of indices of the blocks within $M_{\bfgam}$, and let $k_{\gamma}=|B_{\gamma}|$. Recall that $B_T$ and $k_T$ are the block indices and the number of blocks  respectively in model $\Mt$. 

We shall use the following two lemmas in the proof of this theorem. The results from Lemma \ref{lema1} are slightly generalized versions of Lemmas B.2 and B.3 from \cite{maruyama2011fully} and can be proved in a similar way. Conditions \ref{cond5.1} and \ref{cond5.2} are used to prove Lemma \ref{lema1}.
 
\begin{lem}\label{lema1}
Let $R_{i,\gamma}^2$ and $R_{i,T}^2$ denote the $i^{th}$ component of $R^2$ under an arbitrary model $\Mgam$ and the true model $\Mt$ respectively. \\
(i) Then for $i \in B_T$,
\begin{eqnarray*}
R_{i,\gamma}^2 &\stackrel{P}{\rightarrow}& \frac{\bfb_{i,T}^T D_{i,T} \bfb_{i,T} - V_{i,\gamma}}{\sigma^2 + \sum_{j \in B_T} \bfb_{j,T}^T D_{j,T} \bfb_{j,T} +  \alpha_T^2  } \nonumber \\
R_{i,T}^2 &\stackrel{P}{\rightarrow}& \frac{\bfb_{i,T}^T D_{i,T} \bfb_{i,T} }{\sigma^2 + \sum_{j \in B_T} \bfb_{j,T}^T D_{j,T} \bfb_{j,T} + \alpha_T^2 } 
\end{eqnarray*}
where $D_{i,\gamma}= \lim\limits_{n \rightarrow \infty} \frac{1}{n} X_{i,\gamma}^T  X_{i,\gamma}$ is positive definite for all models $\Mgam$ (and all blocks) and
\begin{equation*}
V_{i,\gamma}= \lim\limits_{n \rightarrow \infty} \frac{1}{n} \bfb_{i,T}^T X_{i,T}^T (I - P_{X_{i,\gamma}}) X_{i,T} \bfb_{i,T} = \Bigg\{ \begin{array}{cc}
     0 & X_{i,\gamma} \supseteq X_{i,T} \\
   > 0 &  X_{i,\gamma} \not \supset X_{i,T} \\
       \end{array}  
\end{equation*}
(ii) When $i \not \in B_T$, $ R_{i,\gamma}^2 \stackrel{P}{\rightarrow} 0 $. \\
(iii) For $i \in B_{\gamma}\backslash B_T$, $n  R_{i,\gamma}^2 \stackrel{d}{\rightarrow} c \chi^2_{p_{i,\gamma}} = O_p (1)$. (c is a constant) \\
(iv) For any model $M_{\gamma} \supset \Mt$, $\left( \frac{1- \sum_{j \in B_T} R^2_{j,T}}{1-  \sum_{j \in B_{\gamma}} R^2_{j,\gamma}} \right)^n = \left( \frac{1-R^2_T}{1-R^2_{\gamma}} \right)^n$ is bounded from above in probability.
\end{lem}

\begin{lem}\label{lema2}
Consider the function $h(\bft)$  defined on $(0,1)^{|I|}$ as $h(\bft)= \sum\limits_{i \in I} b_i \log(1-t_i) - m\log(1-\sum\limits_{i\in I} t_i r_i)$ where $I$ is a set of indices, each $b_i>0$, $r_i\geq0$, $m > \sum\limits_{i \in I} b_i$ and $\sum\limits_{i \in I} r_i < 1$. Then the (unique) maximum of $h(\bft)$ is attained at the point $\bft=\bft^*$ in the interior of the set $(0,1)^{|I|}$ with $t_i^* = 1- \frac{b_i (1-r)}{r_i (m-b)}$ for all $i \in I$, where $b= \sum\limits_{i \in I} b_i$  and $r = \sum\limits_{i \in I} r_i$.\\
If we denote the Hessian matrix as $H(\bft)=(( H_{ij}(\bft) )) = (( \frac{\partial^2 h(\bft)}{\partial t_i t_j} ))$, then 

\begin{eqnarray*}
H_{ij}(\bft^*) &=& \frac{(m-b)^2 r_i r_j}{m (1-r)^2 } \; \;, \; \; i \neq j \\
\mbox{and } H_{ii}(\bft^*) &=& -\frac{(m-b)^2 r_i^2}{(1-r)^2}  \left[ \frac{1}{b_i} -  \frac{1}{m} \right]
\end{eqnarray*}
\end{lem}

\noindent
The proof of Lemma \ref{lema2} is skipped for brevity. It is not difficult to check that the partial derivatives of $h(\bft)$ attain a value of zero at $\bft^*$. The Hessian matrix is non-positive definite at $\bft=\bft^*$.

Returning to the proof of the theorem, first consider the case $\Mt \neq \Mo$.

\noindent
The Bayes factor comparing model $\Mgam$ to the true model $\Mt$ can be written as
\begin{eqnarray}
&& BF(\Mgam : \Mt) = \frac{BF(\Mgam:\Mo)}{BF(\Mt:\Mo)} \nonumber \\
&=& \left( \frac{a-2}{2} \right)^{k_{\gamma}-k_T} \frac{\int \prod\limits_{i \in B_{\gamma}} (1-t_i)^{\frac{a+p_{i,\gamma}}{2}-2} (1-\sum\limits_{i \in B_{\gamma}} t_i R_{i,\gamma}^2)^{-\frac{n-1}{2}} d\bft }{ \int \prod\limits_{i \in B_T} (1-t_i)^{\frac{a+p_{i,T}}{2}-2} (1-\sum\limits_{i \in B_T} t_i R_{i,T}^2)^{-\frac{n-1}{2}} d\bft } \nonumber \\
&=& \left( \frac{a-2}{2} \right)^{k_{\gamma}-k_T} \frac{\int \exp( h_{\gamma}(\bft) ) d\bft}{\int \exp( h_{T}(\bft) ) d\bft} \label{eq:BFlap}
\end{eqnarray}
where $h_{j}(\bft)= \sum\limits_{i \in B_j} (\frac{a+p_i}{2}-2)  \log(1-t_{i}) -  \frac{n-1}{2} \log(1-\sum\limits_{i \in B_j} t_i R_{i,j}^2)$; $j \in \{ \gamma, T\}$. We apply Lemma \ref{lema2} to these two functions which necessitates:  

\noindent
(A) $\frac{a+p_{i,\gamma}}{2} > 2$ for all $i \in B_{\gamma}$ and $\frac{a+p_{i,T}}{2} > 2$ for all $i \in B_{T}$.\\
(B) $n > p+1+ (a-4) \times \min(k_{\gamma},k_T)$.

(B) is satisfied since $n \geq p+2$  and $a \leq 4$. For (A) to hold, we must have $a+p_{i} > 4$ for all blocks in both models. This is true for $a > 3$, and we proceed with this portion of the proof.  The proof for the $2 < a \leq 3$ case is in Appendix~\ref{appb5} in the supplementary materials.


Using the multivariate generalization of the Laplace approximation, \eqref{eq:BFlap} is approximated up to an $O(\frac{1}{n})$ term as
\begin{eqnarray*}
BF(\Mgam:\Mt) &\approx& \left( \frac{a-2}{2} \right)^{k_{\gamma}-k_T}  \frac{|H_T (\hat{\bft}_{T})|^{1/2}}{|H_{\gamma} (\hat{\bft}_{\gamma})|^{1/2}} \times   \\
&&  \frac{\exp \left[ \sum\limits_{i \in B_{\gamma}} b_{i,\gamma} \log(1- \hat{t}_{i,\gamma}) - m \log(1-\sum\limits_{i \in B_{\gamma}} \hat{t}_{i,\gamma} R_{i,\gamma}^2 ) \right]}{\exp\left[ \sum\limits_{i \in B_T} b_{i,T} \log(1- \hat{t}_{i,T}) - m \log(1-\sum\limits_{i \in B_{T}} \hat{t}_{i,T} R_{i,T}^2 ) \right]}
\end{eqnarray*}
where $b_{i,j} = \frac{a+p_{i,j}}{2}-2$, $m = \frac{n-1}{2}$ and $H_j (\hat{\bft}_j)$ is the Hessian matrix of $h_j(\bft)$ evaluated at the maximizer $\hat{\bft}_j$ of $h_j(\bft)$; $j \in \{ \gamma, T \}$.

Using Lemmas \ref{lema1} and \ref{lema2}, $|H_{T} (\hat{\bft}_{T})| = O(m^{2 k_T})$ and $|H_{\gamma} (\hat{\bft}_{\gamma})| = O(m^{2 q_{\gamma}})$ for some $0 \leq q_{\gamma} \leq k_T$. This follows from $( H_{\gamma} (\hat{\bft}_{\gamma}) )_{ii} = O(m^2 R^2_{i,\gamma})$ and so $q_{\gamma}=k_{\gamma}-L_{\gamma}$, where $L_{\gamma}$ is the number of components of $R^2_{i,\gamma}$ going to zero in probability.  

For large $m = \frac{n-1}{2}$ (or equivalently for large $n$), 
\begin{eqnarray*}
BF(\Mgam:\Mt) &\approx&  O(m^{k_T - q_{\gamma}}) \exp \Bigg[ -m \log \left( \frac{m(1-R_{\gamma}^2)}{m - b_{\gamma}} \right) + m \log \left( \frac{m(1-R_{T}^2)}{m - b_{T}} \right) \\
&+& \hspace{-.3 cm}   \sum_{i \in B_{\gamma}} b_{i,\gamma} \log \left( \frac{b_{i,\gamma} (1-R^2_{\gamma})}{R^2_{i,\gamma}(m-b_{\gamma})} \right) -   \sum_{i \in B_{T}} b_{i,T} \log \left( \frac{b_{i,T} (1-R^2_{T})}{R^2_{i,T}(m-b_{T})} \right) \Bigg] \left( \frac{a-2}{2} \right)^{k_{\gamma}-k_T}  \\
&=& \hspace{-.3 cm} O(m^{k_T - q_{\gamma}}) \exp \Big[ m \log\Big( \frac{1-R^2_T}{1-R^2_{\gamma}} \Big) + (b_T - b_{\gamma}) \log m - \sum_{i \in B_{\gamma}} b_{i,\gamma} \log(R_{i,\gamma}^2) + O(1) \Big]
\end{eqnarray*}
where $R^2_j = \sum_{i \in B_j} R^2_{i,j}$ and $b_j = \sum_{i \in B_j} b_{i,j}$ for $j \in \{ \gamma, T \}$.

\noindent
\underline{Case 1}: $\Mgam \not \supset \Mt$

From Lemma \ref{lema1} (i) and (ii), $R_{\gamma}^2 < R_T^2$ and so $ \log\left( \frac{1-R^2_T}{1-R^2_{\gamma}} \right) < 0$. Again $\sum_{i \in B_{\gamma}} b_{i,\gamma} \log(R_{i,\gamma}^2) = \sum_{i \in J^c} b_{i,\gamma} \log(R_{i,\gamma}^2) + \sum_{i \in J} b_{i,\gamma} \log(R_{i,\gamma}^2)  = C +  \sum_{i \in J} b_{i,\gamma} $ $ \log(m R_{i,\gamma}^2) - \sum_{i \in J} b_{i,\gamma} \log m$, where $J \subseteq B_{\gamma}$ is the set of indices such that $R^2_{i,\gamma} \rightarrow 0$ for $i \in J$. Then
\begin{eqnarray*}
\lim_{m \rightarrow \infty} BF(\Mgam: \Mt) = \lim_{m \rightarrow \infty} O(m^s) \cdot  O( f^m) = 0
\end{eqnarray*}
where $0<f<1$ and $s$ is some real number which might be positive or negative depending on the block structures of models $\Mgam$ and $\Mt$. The limit is zero  since the second term goes to zero at an exponential rate and the first term is either bounded in probability ($s=0$) or goes to zero ($s<0$) or to infinity ($s>0$) at a polynomial rate.

\noindent
\underline{Case 2}: $\Mgam \supset \Mt$ \\
\underline{Case 2A}: $\Mgam$ has the same block structure as $\Mt$, i.e., $B_{\gamma}=B_T$, but has more predictors in at least one block.

In this case $R_{\gamma}^2 \geq R_T^2$ but due to Lemma \ref{lema1} (iv), $\left( \frac{1-R^2_T}{1-R^2_{\gamma}} \right)^m$ is bounded in probability. Also Lemma \ref{lema1} (i) confirms that none of the $R_{i,\gamma}^2$ converge to zero in the limit and so $q_{\gamma} =k_{\gamma}= k_T$. 

Note that $b_{\gamma} =  \sum_{i \in B_{\gamma}} [ \frac{a+p_{i,\gamma}}{2}-2 ] = \frac{p_{\gamma}}{2} - (4-a)\frac{k_{\gamma}}{2}$ and $b_T= \frac{p_{T}}{2} - (4-a)\frac{k_T}{2}$. For $\Mgam$, $b_T - b_{\gamma} = \frac{p_T - p_{\gamma}}{2} < 0$  and $k_T = k_{\gamma}$. Hence, 
\begin{eqnarray*}
\lim_{m \rightarrow \infty} BF(\Mgam: \Mt) = \lim_{m \rightarrow \infty} O(1)\exp \left[ O(1) + \left( \frac{p_T - p_{\gamma}}{2} \right) \log m \right] = 0
\end{eqnarray*}

\noindent
\underline{Case 2B}: $\Mgam$ has more blocks than $\Mt$, i.e., $B_{\gamma} \supset B_T$, and in addition has more predictors in at least one of the blocks common to both models.

As before $\left( \frac{1-R^2_T}{1-R^2_{\gamma}} \right)^m$ is bounded in probability and all $R_{i,\gamma}^2 \rightarrow 0$ for $i \in B_{\gamma}\backslash B_T$ which implies that $q_{\gamma}=k_{\gamma}-( k_{\gamma} - k_T ) = k_T$. Since $b_T - b_{\gamma}=  \sum_{i \in B_{T}} \frac{p_{i,T}}{2} -  (4-a)\frac{k_{T}}{2}  -  \sum_{i \in B_{\gamma}}  \frac{p_{i,\gamma}}{2} + (4-a)\frac{k_{\gamma}}{2}$ , we have
\begin{eqnarray*}
 && m \log\left( \frac{1-R^2_T}{1-R^2_{\gamma}} \right) + (b_T - b_{\gamma}) \log m - \sum_{i \in B_{\gamma}} b_{i,\gamma} \log(R_{i,\gamma}^2) \\
&=& O(1) + \left[ \frac{1}{2} \sum_{i \in B_{T}} (p_{i,T}- p_{i,\gamma}) + \frac{1}{2} \sum_{i \in B_{\gamma} \backslash B_T} (0 - p_{i,\gamma}) - (4-a)\frac{k_T - k_{\gamma}}{2} \right] \log m \\
&&  - \sum_{i \in B_{\gamma} \backslash B_T} \hspace{-.2 cm} b_{i,\gamma} \log(m R_{i,\gamma}^2) + \sum_{i \in B_{\gamma} \backslash B_T} \hspace{-.2 cm} b_{i,\gamma} \log m  \hspace{.3 in} \Big(m R^2_{i,\gamma} = O(1)  \mbox{ by Lemma \ref{lema1} (iii)}\Big) \\
&=& O(1) + \log m \Big[  \frac{1}{2} \sum_{i \in B_{T}} (p_{i,T}- p_{i,\gamma}) - \frac{1}{2} \sum_{i \in B_{\gamma} \backslash B_T} p_{i,\gamma} + (a-4)\frac{k_T - k_{\gamma}}{2} \\
&& \hspace{2 in} + \sum_{i \in B_{\gamma} \backslash B_T} \left( \frac{a+p_{i,\gamma}}{2} - 2 \right) \Big] 
\end{eqnarray*}

\begin{eqnarray*}
&=& O(1) + \log m  \left[  \frac{1}{2} \sum_{i \in B_{T}} (p_{i,T}- p_{i,\gamma}) \right] \\
\therefore \lim\limits_{m \rightarrow \infty} BF(\Mgam:\Mt) &=& \lim_{m \rightarrow \infty} O(1)\exp \left[ O(1) + \frac{\sum_{i \in B_{T}} (p_{i,T}- p_{i,\gamma}) }{2} \log m \right] = 0
\end{eqnarray*}
since $p_{i,\gamma} \geq p_{i,T}$ $\forall \; i \in B_T$ with strict inequality for at least one $i$ ensures that the above sum is strictly negative.

\noindent
\underline{Case 2C}: $\Mgam$ has more blocks than $\Mt$, but has the same set of predictors in all blocks common to both.

The Bayes factor in this case is
\begin{eqnarray*}
 BF(\Mgam:\Mt) &\approx& O(1)\exp \left[ O(1) + \frac{\sum_{i \in B_{T}} (p_{i,T}- p_{i,\gamma}) }{2} \log m \right] \\
&=& O(1) \exp [ O(1) + 0 ] = O(1), \; \mbox{for all } m
\end{eqnarray*}
since $p_{i,\gamma} = p_{i,T}$ $\forall \; i \in B_T$. The $O(1)$ term  is  a combination of  constants and random variables. As a consequence, the Bayes factor cannot equal 0 with probability 1. This is the only case where the Bayes factor of any arbitrary model $\Mgam$ compared to the true model $\Mt$ does not go to zero with increasing sample size, violating the principle of model selection consistency.

The case when $\Mt = \Mo$ is identical to Case 2C with $B_T = \phi$. In this case, $P\Big( \lim\limits_{n \rightarrow \infty} BF(\Mgam: \Mt) = 0 \; \Big) = 0$. Model selection consistency does not hold.

\subsection{ Proof of Theorem \ref{thm5.3}}  \label{appa9}

Using the same notation as in Theorem \ref{thm5.2}, we observe that block orthogonality of the design gives
\begin{eqnarray*}
 E(\bfb \mid \bfy, \Mgam) &=&  E \left[ E(\bfb \mid \bfy, \bfg, \Mgam) \right] \\
&=& \int \left( \begin{array}{c}
\frac{g_{i_1}}{1+g_{i_1}} \widehat{\bfb}_{i_1,\bld{\gamma},LS}\\
\cdots \\
\cdots \\
\frac{ g_{i_{k_{\bldsm{\gamma}}}} }{ 1+g_{i_{k_{\bldsm{\gamma}}}} } \widehat{\bfb}_{i_{k_{\bldsm{\gamma}}}, \bld{\gamma}, LS} 
\end{array} \right) \pi(\bfg \mid \bfy, \Mgam) d\bfg \\
\mbox{and } E(\alpha \mid \bfy) &=& \widehat{\alpha}_{LS}
\end{eqnarray*}
assuming $B_{\bfgam} =\{ i_1, i_2, \ldots, i_{k_{\bld{\gamma}}} \}$. 

When $\Mt=\Mo$, $\widehat{\bfb}_{\bfgam, LS} \stackrel{P}{\rightarrow} 0$ and $\widehat{\alpha}_{LS} \stackrel{P}{\rightarrow} \alpha_{T}$ for every $\Mgam$ since least squares estimators are consistent. Thus the model averaged prediction $\widehat{y}^*_n$ converges to $E(y^*) = \alpha_T$.

Denote the set of all models belonging to Case 2C of Theorem \ref{thm5.2} together with $\Mt$ by $\Omega$. We have shown that, as $n \rightarrow \infty$, $\pi(\Mgam \mid \bfy) \rightarrow 0$ for any model $\Mgam \not \in \Omega$. Thus, $\lim \limits_{n \rightarrow \infty} \sum\limits_{\bfgam: \Mgam \in \Omega} \pi(\Mgam \mid \bfy) = 1$. 

When $\Mt \neq \Mo$, the least squares estimates are consistent for $\Mgam \in \Omega$, so that $\widehat{\bfb}_{i, \bfgam, LS} \stackrel{P}{\rightarrow} \bfb_{i,T}$ for $i \in B_T$ and $\widehat{\bfb}_{i, \bfgam, LS} \stackrel{P}{\rightarrow} 0$ for $i \not\in B_T$. Hence 
\begin{eqnarray*}
\lim\limits_{n \rightarrow \infty} \widehat{y}_n^*  &=&   \alpha_{T} +  \sum_{\bfgam:\Mgam \in \Omega} \lim_{n \rightarrow \infty} \pi(\Mgam \mid \bfy)  \bfx^{* T} \times \\
 &&    \left( \begin{array}{c}
\bfb_{i_1,T} \lim\limits_{n \rightarrow \infty} \int \frac{g_{i_1}}{1+g_{i_1}} \pi(\bfg \mid \bfy, \Mgam) d\bfg \\
\cdots \\
\cdots \\
\bfb_{i_{k_{\bldsm{\gamma}}}, T} \lim\limits_{n \rightarrow \infty} \int \frac{ g_{i_{k_{\bldsm{\gamma}}}} }{ 1+g_{i_{k_{\bldsm{\gamma}}}} } \pi(\bfg \mid \bfy, \Mgam) d\bfg  \end{array} \right). 
\end{eqnarray*}

Use Lemma \ref{lem5.1} and Corollary \ref{cor5.1} to get $\lim \limits_{n \rightarrow \infty} E \left( \frac{g_i}{1+g_i} \mid \Mgam, \bfy \right) = 1$ $\forall$ $i \in B_{T}$ when $\Mgam \in \Omega$, while  $0 \leq  E \left( \frac{g_i}{1+g_i} \mid \Mgam, \bfy \right) \leq 1$ for any other $i$. So, 
\begin{eqnarray*}
\lim_{n \rightarrow \infty} \widehat{y}_n^* &=& \alpha_{T} + \bfx^{* T} \bfb_{T} \lim_{n \rightarrow \infty} \sum_{\bfgam:\Mgam \in \Omega} \pi(\Mgam \mid \bfy)  = E(y^*)
\end{eqnarray*}
indicating that the block hyper-$g$ prior is prediction consistent under {\em BMA}.

\bibliographystyle{apa}
\bibliography{BlockHyperG}

\section{Supplementary Materials}\label{appB}

\subsection{Proof of Corollary \ref{cor3.1}} \label{appb1}

It can be verified that 
\begin{eqnarray*}
\pi(\sigma^2,g \mid \bfy) &\propto&  (1+g)^{-\frac{a+p}{2}} \frac{1}{\sigma^{n+1}} \exp\left[ -\frac{1}{2\sigma^2} \bfy^T (I - \frac{g}{1+g} P_{X})  \bfy  \right] \\
\mbox{So } \pi(\sigma^2 \mid \bfy) &\propto& \int_0^{\infty} (1+g)^{-\frac{a+p}{2}} \frac{1}{\sigma^{n+1}} \exp\left[ -\frac{||\bfy||^2}{2\sigma^2} (1-\frac{g}{1+g} R^2) \right] dg \\
&\propto& \frac{1}{\sigma^{n+1}}\int_0^1 (1-t)^{\frac{a+p}{2}-2} \exp \left[-\frac{||\bfy||^2}{2\sigma^2} \{ (1-R^2) + R^2(1-t)\}  \right] dt \\
&\propto& \frac{1}{\sigma^{n+1}} \exp \left[ -\frac{||\bfy||^2 (1-R^2)}{2\sigma^2} \right] \int_0^1 (1-t)^{\frac{a+p}{2}-2} \exp \left[ -\frac{||\bfy||^2 R^2}{2\sigma^2} (1-t) \right] dt
\end{eqnarray*}
 
Now $||\bfy||^2 (1-R^2) = (n-p-1) \widehat{\sigma}^2$, which is fixed for all $N$, so that
\vspace{-.3 cm}
\begin{eqnarray*}
\pi(\sigma^2 \mid \bfy) &\propto& \frac{1}{\sigma^{n+1}} \exp \left[ -\frac{(n-p-1)\widehat{\sigma}^2}{2\sigma^2} \right] \int_0^{\frac{||\bld{y}||^2 R^2}{2\sigma^2}} x^{\frac{a+p}{2} -2} e^{-x} dx \times (\sigma^2)^{\frac{a+p}{2}-1}\\
&\propto& \frac{1}{(\sigma^2)^{(n+1)/2-(a+p)/2+1}} \exp \left[ -\frac{(n-p-1)\widehat{\sigma}^2}{2\sigma^2} \right] \int_0^{\frac{||\bld{y}||^2 R^2}{2\sigma^2}} x^{\frac{a+p}{2} -2} e^{-x} dx 
\end{eqnarray*}

As $N \rightarrow \infty, ||\bfy|| \rightarrow \infty \; \mbox{and} \; R^2 \rightarrow 1 \; \mbox{so that} \; \frac{||\bld{y}||^2 R^2}{2\sigma^2} \rightarrow \infty$. (Lemma \ref{lem3.1})
\[ \implies \lim_{N \rightarrow \infty} \pi(\sigma^2 \mid \bfy) \propto \frac{1}{(\sigma^2)^{(n+1)/2-(a+p)/2+1}} \exp \left[ -\frac{(n-p-1)\widehat{\sigma}^2}{2\sigma^2} \right] \int_0^{\infty} x^{\frac{a+p}{2} -2} e^{-x} dx     \]
\[ \propto \frac{1}{(\sigma^2)^{(n+1)/2-(a+p)/2+1}} \exp \left[ -\frac{(n-p-1)\widehat{\sigma}^2}{2\sigma^2} \right]  \]

Thus, the distribution of $\lim\limits_{N \rightarrow \infty} \pi(\sigma^2 \mid \bfy)$ is Inverse Gamma with shape parameter= $\frac{n+1-a-p}{2}$ and scale parameter = $\frac{2}{(n-p-1)\widehat{\sigma}^2}$. The mean of this distribution is $\frac{(n-p-1)\widehat{\sigma}^2}{(n-p-a-1)}$.

\subsection{Proof of Theorem \ref{thm3.3}}\label{appengeneralizedg}

Using the expressions derived in  \cite{maruyama2011fully} and defining the two models $M_{1}$ and $M_{2}$ as in Theorem \ref{thm3.2}, it follows that
\begin{eqnarray*}
BF(M_{2}:M_{1}) &=& C \; \frac{(1-R_1^2)^{(n-p_1-2)/2-a}}{(1-R_2^2)^{(n-p-2)/2-a}} \frac{(1-Q_2^2)^{-p/2-a-1}}{(1-Q_1^2)^{-p_1/2-a-1}}
\end{eqnarray*}
where $C$ is a constant, $R_1^2 = \lambda_1^2+\ldots+\lambda_{p_1}^2$, $R_2^2 = \tau_1^2+\ldots+\tau_{p}^2$, $Q_1^2 = \sum_{i=1}^{p_1} (1-\frac{1}{\nu_{i1}})\lambda_i^2$ and $Q_2^2 = \sum_{i=1}^{p} (1-\frac{1}{\nu_{i2}}) \tau_i^2$. Here $\lambda_i = \mbox{corr}(\boldsymbol{u_i},\bfy)$, the correlation between the response $\bfy$ and $\boldsymbol{u_i}$, the $i^{th}$ principal component of $X_1$ in model $M_1$ and $\tau_i$ is the corresponding entity for model $M_2$ with $X=(X_1, X_2)$. $\nu_{i1}$ and $\nu_{i2}$ are arbitrary constants appearing in the prior covariance matrix with the restriction that $\nu_{ij} \geq \nu_{i(j+1)} \geq 1$ for all $i, j$.

For the asymptotic considered here and because of the block orthogonality assumption ($X_1 \indep X_2$), $R_1^2 \rightarrow 1$ and $R_2^2 - R_1^2 \rightarrow 0$ as in Theorem \ref{thm3.2}.  Since $\nu_{ij} \geq 1$ $\forall \; i,j$, as $N \rightarrow \infty$, $Q_1^2 \rightarrow \eta_1$ and $Q_2^2 \rightarrow \eta_2$  for some $\eta_1, \eta_2$ satisfying $1 > \eta_1,\eta_2 \geq 0$.

\begin{eqnarray*}
\lim_{||\bfb_1|| \rightarrow \infty} BF(M_{2}:M_{1}) & = & \lim_{\substack{z \rightarrow 1 \\ q \rightarrow 0} } C \; \frac{(1-z)^{(n-p_1-2)/2-a}}{(1-z-q)^{(n-p-2)/2-a}} \frac{(1-\eta_2)^{-p/2-a-1}}{(1-\eta_1)^{-p_1/2-a-1}} \\
&=& C^* \times \lim_{\substack{z \rightarrow 1 \\ q \rightarrow 0} }  \frac{(1-z)^{(n-p_1-2)/2-a}}{(1-z-q)^{(n-p-2)/2-a}} \\
&=& C^* \times \lim_{z \rightarrow 1} (1-z)^{(p-p_1)/2}=0. 
\end{eqnarray*}

\subsection{Proof of Theorem \ref{thm3.4}}\label{appenb3}

As before, consider the  setup as in Theorem \ref{thm3.2} and assume that $X_1 \indep X_2$ without  loss of generality. If the design matrix is not block orthogonal, we can use the technique in the proof of  Theorem \ref{thm3.2}  since the robust prior is also a specific mixture of a $g$ prior and transform the design to be block orthognal. Using expressions derived in \cite{bayarri2012criteria}
\begin{eqnarray*}
BF(M_{2}:M_{1})= C \left( \frac{Q_{10}}{Q_{20}} \right)^{(n-1)/2} \frac{ _2F_1 \left[ (p+1)/2; (n-1)/2; (p+3)/2; \frac{(1-Q_{20}^{-1})(p+1)}{1+n} \right]}{ _2F_1 \left[ (p_1+1)/2; (n-1)/2; (p_1+3)/2; \frac{(1-Q_{10}^{-1})(p_1+1)}{1+n} \right]}
\end{eqnarray*} 
where $Q_{10}=1-R_1^2$ and $Q_{20}=1-R^2_2$.

Since $ _2F_1(a;b;c;z)=\; _2F_1(b;a;c;z) $, $BF(M_{2}:M_{1})$
\begin{eqnarray*}
&=& C \left( \frac{Q_{10}}{Q_{20}} \right)^{-(n-1)/2} \frac{ _2F_1 \left[  (n-1)/2; (p+1)/2; (p+3)/2; \frac{(1-Q_{20}^{-1})(p+1)}{1+n} \right]}{ _2F_1 \left[ (n-1)/2; (p_1+1)/2;  (p_1+3)/2; \frac{(1-Q_{10}^{-1})(p_1+1)}{1+n} \right]} \\
&=&  C \left( \frac{1-z}{1-z-q} \right)^{(n-1)/2} \frac{ _2F_1 \left[  (n-1)/2; (p+1)/2; (p+3)/2; -\frac{(z+q)(p+1)}{(1+n)(1-z-q)} \right]}{ _2F_1 \left[ (n-1)/2; (p_1+1)/2;  (p_1+3)/2; -\frac{z (p_1+1)}{(1+n)(1-z)} \right]} \\
&=&  C^* \left( \frac{1-z}{1-z-q} \right)^{(n-1)/2} \frac{\int_{0}^1 t^{(p-1)/2} \left( 1+\frac{t(z+q)(p+1)}{(n+1)(1-z-q)} \right)^{-(n-1)/2} dt }{\int_{0}^1 t^{(p_1-1)/2} \left( 1+\frac{t z(p_1+1)}{(n+1)(1-z)} \right)^{-(n-1)/2} dt}
\end{eqnarray*} 
where $z = R_1^2$, $q=R_2^2-R_1^2$, and $C$ and $C^*$ are constants.  Here $N \rightarrow \infty$ implies $z \uparrow 1$ and $q \downarrow 0$. This is because $R^2_{2} = R^2_{1} + \frac{(X_2 \hat{\beta}_2)^T (X_2 \hat{\beta}_2)}{y^T y}$ due to the block orthogonal design.

Hence  $\lim\limits_{||\bfb_1|| \rightarrow \infty} BF(M_{2}:M_{1})$
\begin{eqnarray}
&=& \lim_{\substack{z \rightarrow 1 \\ q \rightarrow 0}} C^* \left( \frac{1-z}{1-z-q} \right)^{(n-1)/2}  \frac{\int_{0}^1 t^{(p-1)/2} \left( 1+\frac{t(z+q)(p+1)}{(n+1)(1-z-q)} \right)^{-(n-1)/2} dt }{\int_{0}^1 t^{(p_1-1)/2} \left( 1+\frac{t z(p_1+1)}{(n+1)(1-z)} \right)^{-(n-1)/2} dt}  \nonumber \\
&=& \lim_{\substack{z \rightarrow 1 \\ q \rightarrow 0}} C^* \frac{\int_{0}^1 t^{(p-1)/2} \left( 1-z-q +\frac{t(z+q)(p+1)}{n+1} \right)^{-(n-1)/2} dt }{\int_{0}^1 t^{(p_1-1)/2} \left( 1-z +\frac{t z(p_1+1)}{n+1} \right)^{-(n-1)/2} dt}   \nonumber \\
&=& \lim_{z \rightarrow 1}  C^* \frac{\int_{0}^1 t^{(p-1)/2} \left( 1-z +\frac{t z (p+1)}{n+1} \right)^{-(n-1)/2} dt }{\int_{0}^1 t^{(p_1-1)/2} \left( 1-z +\frac{t z(p_1+1)}{n+1} \right)^{-(n-1)/2} dt} \nonumber \\
&&  \hspace{1.1 in} (\mbox{by the Monotone Convergence Theorem}) \nonumber \\
&=&    \lim_{z \rightarrow 1}  K \frac{\int_{0}^1 t^{(p-1)/2} \left( B+ t \right)^{-(n-1)/2}  dt }{\int_{0}^1 t^{(p_1-1)/2} \left( B^* + t \right)^{-(n-1)/2}    dt} \label{eq:robustlim}
\end{eqnarray} 
where $K$ is a constant, $B=B(z)=\frac{(n+1)(1-z)}{(p+1)z}$ and $B^*=B^*(z)=\frac{(n+1)(1-z)}{(p_1+1)z}$. Clearly both $B$ and $B^*$ go to zero as $z \rightarrow 1$ and $\frac{B}{B^*}= \frac{p_1+1}{p+1}$  for all $z$.

The limit in \eqref{eq:robustlim} is zero when $n > p_1 + 2$ (see Appendix \ref{appb2} for a proof).

\subsection{Proof of Invariance of Block $\bfg$ Priors to Reparameterizations by Blockwise Affine Transformations} \label{appb7}

Suppose we have two separate regression problems

\noindent
\underline{Problem 1}:
\begin{eqnarray*} 
\bfy \mid \alpha, \bfb, \sigma^2  &\sim& N( \alpha \bone + X_1 \beta_1 + \ldots + X_k \beta_k , \sigma^2 I)  \\
\bfb \mid \bfg, \sigma^2 &\sim& N(0,A \sigma^2) \\
\pi(\alpha,\sigma^2)  &\propto& \frac{1}{\sigma^2}
\end{eqnarray*}
with  $ A = \left( \begin{array}{cccc}
    g_1 (X_1^T X_1)^{-1} & 0 &  \cdots & 0 \\
    0 & g_2 (X_2^T X_2)^{-1} & \cdots & 0 \\ 
    \vdots & \vdots &  \ddots & \vdots \\
     0 & 0 & \cdots & g_k (X_k^T X_k)^{-1} \\
     \end{array} \right)$\\
		
\noindent
and \underline{Problem 2}:
\begin{eqnarray*} 
\bfy \mid \alpha, \bfgam, \sigma^2 &\sim& N( \alpha \bone + Z_1 \gamma_1 + \ldots + Z_k \gamma_k , \sigma^2 I)  \\
\bfgam \mid \bfg, \sigma^2 &\sim& N(0, B \sigma^2) \\
\pi(\alpha,\sigma^2)  &\propto& \frac{1}{\sigma^2}
\end{eqnarray*}
with  $ B = \left( \begin{array}{cccc}
    g_1 (Z_1^T Z_1)^{-1} & 0 &  \cdots & 0 \\
    0 & g_2 (Z_2^T Z_2)^{-1} & \cdots & 0 \\ 
    \vdots & \vdots &  \ddots & \vdots \\
		0 & 0 & \cdots & g_k (Z_k^T Z_k)^{-1} \\
		\end{array} \right)$\\

\noindent
We wish to show that if  Problem 2 is a reparameterization of Problem 1 by a linear map acting within blocks, then inference is the same for both problems.

If the predictors in Problem 2 are a within-block linear transformation of the predictors in Problem 1, then there exist non-singular matrices $P_i$ such that $Z_i = X_i P_i$ for each $i=1,2,...,k$.  If we define  $ P = \left( \begin{array}{cccc}
    P_1 & 0 &  \cdots & 0 \\
    0 & P_2 & \cdots & 0 \\ 
    \vdots & \vdots &  \ddots & \vdots \\
		0 & 0 & \cdots & P_k \\
		\end{array} \right),$
it is clear that $Z=(Z_1,Z_2,...,Z_k)$ and $X=(X_1,X_2,...,X_k)$ are related by $Z=XP$. So $\bfy=\alpha \bone+Z\bfgam + \epsilon$ can be rewritten as $\bfy=\alpha \bone+XP\bfgam + \epsilon$. For each $i$, 
\[ (Z_i^T Z_i)^{-1} = [(X_i P_i)^T (X_i P_i)]^{-1} = [P_i^T X_i^T X_i P_i]^{-1} = P_i^{-1} (X^T_i X_i)^{-1} (P_i^T)^{-1}
\]

\noindent
It follows now that $B = P^{-1}A(P^T)^{-1}$ so that
\begin{eqnarray*}
\bfgam \mid \bfg, \sigma^2 &\sim& N(0, B \sigma^2) = N(0, P^{-1}A(P^T)^{-1} \sigma^2) \\
\implies P\bfgam \mid \bfg, \sigma^2 &\sim& N(0, P \left[ P^{-1}A(P^T)^{-1} \sigma^2 \right] P^T) = N(0,A\sigma^2) \\
&\stackrel{d}{=}& \bfb \mid \bfg, \sigma^2
\end{eqnarray*}

\noindent
The equivalence of these two priors signifies that inferences from Problem 1 and from Problem 2 will be exactly the same. Hence the block $g$ prior is invariant under a blockwise affine transformation of the problem in any general design.

\subsection{Proof of Lemma \ref{lem4.2}} \label{appa3}

The proof relies on a result in, for example, \cite{lehmann2005testing} on stochastic ordering of random variables.

\textbf{\underline{UR: A Useful Result}} \label{UR} (\textbf{Stochastic ordering of densities/random variables}) 

Let $p_{\theta}(x)$ be a family of densities on the real line with monotone likelihood ratio in x. Then 
\begin{enumerate}
\item For any $\theta <\theta^{\prime}$, the cumulative distribution functions of X under $\theta$ and $\theta^{\prime}$ satisfy $F_{\theta^{\prime}}(x) \leq F_{\theta}(x)$ for all $x$.
\item If $\psi$ is a non-decreasing function of $x$, then $E_{\theta}(\psi(X))$ is a non-decreasing function of $\theta$.   
\end{enumerate}

\noindent
If $C_1$  and $C_2$ denote the normalizing constants for the two densities then
\begin{eqnarray*}
f_1(t_m) &=& \frac{1}{C_1}  \int_{(0,1)^{k-1}} \left[ \prod_{i=1}^k (1-t_i)^{\frac{a+p_i}{2}-2} \right] (1-\sum_{i=1}^k t_i R_i^2)^{-\frac{n-1}{2}} d\bft_{-m} \\
f_2(t_m) &=& \frac{1}{C_2}  \int_{(0,1)^{k-1}} \left[ \prod_{i=1}^k (1-t_i)^{\frac{a+p_i}{2}-2} \right] (1- t_j R_j^2)^{-\frac{n-1}{2}} d\bft_{-m} \\
\mbox{So,} \; \frac{f_1(t_m)}{f_2(t_m)} &=& \frac{C_2}{C_1} \frac{ \int \prod_{i=1}^k (1-t_i)^{\frac{a+p_i}{2}-2} (1-\sum_{i=1}^k t_i R_i^2)^{-\frac{n-1}{2}} d\bft_{-m} }{ \int \prod_{i=1}^k (1-t_i)^{\frac{a+p_i}{2}-2} (1- t_j R_j^2)^{-\frac{n-1}{2}} d\bft_{-m}  } 
\end{eqnarray*}

\noindent
\underline{Case 1}: $m \neq j$

\[ \frac{f_1(t_m)}{f_2(t_m)} = \frac{C_2}{C_1} \frac{(1-t_m)^{\frac{a+p_m}{2}-2}}{(1-t_m)^{\frac{a+p_m}{2}-2}}  \frac{ \int \prod_{i \neq m} (1-t_i)^{\frac{a+p_i}{2}-2} (1-\sum_{i=1}^k t_i R_i^2)^{-\frac{n-1}{2}} d\bft_{-m} }{ \int \prod_{i \neq m} (1-t_i)^{\frac{a+p_i}{2}-2} (1- t_j R_j^2)^{-\frac{n-1}{2}} d\bft_{-m}  }  \] 
\[  = \frac{C_2}{C_1} \frac{ \int \prod_{i \neq m} (1-t_i)^{\frac{a+p_i}{2}-2} (1-\sum_{i \neq m}^k t_i R_i^2 - t_m R^2_m)^{-\frac{n-1}{2}} d\bft_{-m} }{ \int \prod_{i \neq m} (1-t_i)^{\frac{a+p_i}{2}-2} (1- t_j R_j^2)^{-\frac{n-1}{2}} d\bft_{-m}  }
\]

\noindent
Note that $(1- \sum_{i=1}^k t_i R_i^2)^{-(n-1)/2}$ is a non-decreasing function of $t_m$, and so $\frac{f_1(t_m)}{f_2(t_m)}$ is a non-decreasing function of $t_m$.

Applying the Useful Result (UR) stated earlier, $f_1$ is stochastically larger than $f_2$ and, using part (2) of the UR with the strictly increasing function $\psi(t_m)=t_m$, we have $E_{f_1}(t_m) \geq E_{f_2}(t_m)$. In fact, the ratio $\frac{f_1(t_m)}{f_2(t_m)}$ is strictly increasing and we have strict inequality in the result when $R_m^2 > 0$.
\; \\

\noindent
\underline{Case 2}: $m=j$

\begin{eqnarray*}
\frac{f_1(t_j)}{f_2(t_j)} &=& \frac{C_2}{C_1} \frac{(1-t_j)^{\frac{a+p_j}{2}-2}}{(1-t_j)^{\frac{a+p_j}{2}-2}}  \frac{ \int \prod_{i \neq j} (1-t_i)^{\frac{a+p_i}{2}-2} (1-\sum_{i=1}^k t_i R_i^2)^{-\frac{n-1}{2}} d\bft_{-j} }{ \int \prod_{i \neq j} (1-t_i)^{\frac{a+p_i}{2}-2} (1- t_j R_j^2)^{-\frac{n-1}{2}} d\bft_{-j}  }   \\
&=&  \frac{C_2}{C_1} \frac{1}{(1- t_j R_j^2)^{-\frac{n-1}{2}}} \frac{ \int \prod_{i \neq j} (1-t_i)^{\frac{a+p_i}{2}-2} (1-\sum_{i=1}^k t_i R_i^2)^{-\frac{n-1}{2}} d\bft_{-j} }{ \int \prod_{i \neq j} (1-t_i)^{\frac{a+p_i}{2}-2} d\bft_{-j} }  \\
&=& \frac{C_2}{C_1} \frac{ \int \prod_{i \neq j} (1-t_i)^{\frac{a+p_i}{2}-2} (1-\sum_{i \neq j} t_i \frac{R_i^2}{1-t_j R_j^2})^{-\frac{n-1}{2}} d\bft_{-j} }{ \int \prod_{i \neq j} (1-t_i)^{\frac{a+p_i}{2}-2} d\bft_{-j} }  
\end{eqnarray*}

\noindent
The function $(1-\sum_{i \neq j} t_i \frac{R_i^2}{1-t_j R_j^2})^{-\frac{n-1}{2}}$ is non-decreasing in $t_j$ and so $\frac{f_1(t_j)}{f_2(t_j)}$ is a non-decreasing function of $t_j$. Using UR, we conclude that $E_{f_1}(t_j) \geq E_{f_2}(t_j)$. Strict inequality holds in this case when the ratio $\frac{f_1(t_m)}{f_2(t_m)}$ is strictly increasing and this happens when $R^2_m > 0$ and at least one $R^2_i > 0$, $i \neq m$.

\subsection{Proof of Corollary \ref{cor4.1}} \label{appb3}

It is easy to show that
\[ \pi(\sigma^2, \bfg \mid \bfy) \propto \prod_{i=1}^k (1+g_i)^{-\frac{a+p_i}{2}} \frac{1}{\sigma^{n+1}} \exp\left[ -\frac{1}{2\sigma^2} \bfy^T(I - \sum_{i=1}^k \frac{g_i}{1+g_i} P_{X_i} )\bfy \right] \]

 So, $\pi(\sigma^2 \mid \bfy)$
\begin{eqnarray*}
&\propto&  \frac{1}{\sigma^{n+1}} \int_{(0,\infty)^k} \prod_{i=1}^k (1+g_i)^{-\frac{a+p_i}{2}} \exp\left[ -\frac{||\bfy||^2}{2\sigma^2} (1- \sum_{i=1}^k \frac{g_i}{1+g_i} R_i^2) \right] d\bfg \\
 &\propto&  \frac{1}{\sigma^{n+1}} \int_{(0,1)^k} \prod_{i=1}^k (1 - t_i)^{\frac{a+p_i}{2}-2} \exp\left[ -\frac{||\bfy||^2}{2\sigma^2} (1- \sum_{i=1}^k t_i R_i^2) \right] d\bft \\
 &\propto&  \frac{1}{\sigma^{n+1}} \int_{(0,1)^k} \prod_{i=1}^k (1 - t_i)^{\frac{a+p_i}{2}-2} \exp\left[ -\frac{||\bfy||^2}{2\sigma^2} (1- \sum_{i=1}^k R_i^2 + \sum_{i=1}^k (1-t_i)R^2_i) \right] d\bft \\
&\propto& \frac{1}{\sigma^{n+1}}  \exp \left[ -\frac{||\bfy||^2 (1- \sum_{i=1}^k R_i^2)}{2\sigma^2} \right] \times \prod_{i=1}^k  \left[ \int_0^1 (1-t_i)^{\frac{a+p_i}{2}-2} \exp\left( -\frac{||\bfy||^2}{2\sigma^2} R_i^2(1-t_i) \right)  dt_i \right] \\
 &\propto& \frac{1}{\sigma^{n+1}}  \exp \left[ -\frac{(n-p-1)\widehat{\sigma}^2}{2\sigma^2} \right] \prod_{i=1}^k \left[ \int_0^{\frac{||\bld{y}||^2 R_i^2}{2\sigma^2}} x_i^{\frac{a+p_i}{2}-2} e^{-x_i} dx_i \right] \times \prod_{i=1}^k (\sigma^2)^{\frac{a+p_i}{2}-1} 
\end{eqnarray*}

\noindent
Now as $N \rightarrow \infty, ||\bfy|| \rightarrow \infty \; \mbox{and} \; R_1^2 \rightarrow 1$ while $R_i^2 \rightarrow 0$ $\forall \; i \neq 1 $ (Lemma \ref{lem4.1}). 
The expression $\frac{||\bld{y}||^2 R_i^2}{2\sigma^2} = \frac{\bld{y}^T P_{X_i} \bld{y}}{2\sigma^2} =\frac{(X_i \widehat{\bld{\beta}}_i)^T (X_i\widehat{\bld{\beta}}_i)}{2\sigma^2}$  is constant for $i \neq 1$ and goes to $\infty$ for $i=1$. So only $\frac{||\bld{y}||^2 R_1^2}{2\sigma^2} \rightarrow \infty$, while the other integrals are over a finite unchanging domain.
\begin{eqnarray*}
\lim_{N \rightarrow \infty} \pi(\sigma^2 \mid \bfy) &\propto&  \frac{1}{(\sigma^2)^{\frac{n+1}{2}+k -\frac{ka+\sum_{i=1}^k p_i}{2}}} \exp \left[ -\frac{(n-p-1)\widehat{\sigma}^2}{2\sigma^2} \right] \\
 && \hspace{.8 in} \times \lim_{N \rightarrow \infty}  \prod_{i=2}^k \left[ \int_0^{\frac{||\bld{y}||^2 R_i^2}{2\sigma^2}} x_i^{\frac{a+p_i}{2}-2} e^{-x_i} dx_i \right] 
\end{eqnarray*}

\[ \propto  \frac{1}{(\sigma^2)^{\frac{n+1}{2}+k -\frac{ka+p}{2}}} \exp \left[ -\frac{(n-p-1)\widehat{\sigma}^2}{2\sigma^2} \right] \prod_{i=2}^k \gamma \left(\frac{a+p_i}{2}-1, \frac{(X_i \widehat{\bld{\beta}}_i)^T (X_i\widehat{\bld{\beta}}_i)}{2\sigma^2}  \right) \]

\noindent 
where $\gamma(\cdot,\cdot)$ is the lower incomplete gamma function. 

The normalizing constant in the density $\pi (\sigma^2 \mid \bfy)$ exists for any $N$ since 
\begin{eqnarray*}
&& \frac{1}{(\sigma^2)^{[n-1-k(a-2)-p]/2 \;+ 1} } \exp \left[ -\frac{(n-p-1)\widehat{\sigma}^2}{2\sigma^2} \right]   \prod_{i=2}^k \left[ \int_0^{\frac{||\bld{y}||^2 R_i^2}{2\sigma^2}} x_i^{\frac{a+p_i}{2}-2} e^{-x_i} dx_i \right] \\
&& \leq \frac{1}{(\sigma^2)^{[n-1-k(a-2)-p]/2 \;+ 1} } \exp \left[ -\frac{(n-p-1)\widehat{\sigma}^2}{2\sigma^2} \right]   \prod_{i=2}^k \left[ \int_0^{\infty} x_i^{\frac{a+p_i}{2}-2} e^{-x_i} dx_i \right]
\end{eqnarray*}
which is integrable as a function of $\sigma^2$ over $(0, \infty)$ if $n > k(a-2)+p+1$. The integral (over $\sigma^2$) for the expression above is also clearly bounded away from zero as long as $||\bfy||^2 R_i^2$ does not equal  zero or converge to zero for any $i=2,\ldots,k$. But for any such $i$, $||\bfy||^2 R_i^2 = \bfy^T P_{X_i} \bfy$ is strictly greater than zero with probability one in every element of the sequence $\{ \Psi_N \}$  (and also in the limit). This guarantees that the normalizing constant in $\pi (\sigma^2 \mid \bfy)$  is finite and non-zero for all $N$ and validates the existence of a proper limiting distribution for the sequence of posteriors of $\sigma^2$.

\subsection{Proof of Corollary \ref{cor4.2}}\label{appb4}

Proceeding as in Corollary \ref{cor4.1},
\begin{eqnarray*}
\pi(\sigma^2 \mid \bfy) &\propto&  \frac{1}{\sigma^{n+1}} \int_{(0,\infty)^k} \prod_{i=1}^k (1+g_i)^{-\frac{a+p_i}{2}} \exp\left[ -\frac{||\bfy||^2}{2\sigma^2} (1- \sum_{i=1}^k \frac{g_i}{1+g_i} R_i^2) \right] d\bfg 
\end{eqnarray*}
\begin{eqnarray*}
&&  \mbox{So,  } E(\sigma^2 \mid \bfy) \\
 &=& \frac{ \int_{(0,1)^k} \int_0^{\infty} \frac{1}{\sigma^{n-1}} \exp \left[ -\frac{||\bld{y}||^2}{2\sigma^2} (1-\sum_{i=1}^k t_i R^2_i) \right] \prod_{i=1}^k (1-t_i)^{\frac{a+p_i}{2}-2} d\sigma^2 d\bft   }{ \int_{(0,1)^k} \int_0^{\infty} \frac{1}{\sigma^{n+1}} \exp \left[ -\frac{||\bld{y}||^2}{2\sigma^2} (1-\sum_{i=1}^k t_i R^2_i) \right] \prod_{i=1}^k (1-t_i)^{\frac{a+p_i}{2}-2} d\sigma^2  d\bft } \\
&=& \frac{ \int_{(0,1)^k}  \prod_{i=1}^k (1-t_i)^{\frac{a+p_i}{2}-2} \int_0^{\infty} \frac{1}{(\sigma^2)^{\frac{n-3}{2}+1}} \exp \left[ -\frac{||\bld{y}||^2}{2\sigma^2} (1-\sum_{i=1}^k t_i R^2_i) \right]  d\sigma^2 d\bft   }{ \int_{(0,1)^k}  \prod_{i=1}^k (1-t_i)^{\frac{a+p_i}{2}-2} \int_0^{\infty} \ \frac{1}{(\sigma^2)^{\frac{n-1}{2}+1}} \exp \left[ -\frac{||\bld{y}||^2}{2\sigma^2} (1-\sum_{i=1}^k t_i R^2_i) \right] d\sigma^2  d\bft } \\
&=& \frac{||\bfy||^2}{n-3} \frac{\int_{(0,1)^k} \prod_{i=1}^k (1-t_i)^{\frac{a+p_i}{2}-2} (1-\sum_{i=1}^k t_i R^2_i)^{-\frac{n-3}{2}} d\bft }{ \int_{(0,1)^k} \prod_{i=1}^k (1-t_i)^{\frac{a+p_i}{2}-2} (1-\sum_{i=1}^k t_i R^2_i)^{-\frac{n-1}{2}} d\bft} 
\end{eqnarray*}
This leads to
\begin{eqnarray} \label{eq:cor4.2ineq} 
 \hspace{1 in}  E(\sigma^2 \mid \bfy)  &\leq&  \frac{||\bfy||^2}{n-3} \frac{\int_{(0,1)^k} \prod_{i=1}^k (1-t_i)^{\frac{a+p_i}{2}-2} (1- t_1 R^2_1)^{-\frac{n-3}{2}} d\bft }{ \int_{(0,1)^k} \prod_{i=1}^k (1-t_i)^{\frac{a+p_i}{2}-2} (1- t_1 R^2_1)^{-\frac{n-1}{2}} d\bft} \\  
&=& \frac{||\bfy||^2}{n-3} \frac{\int_{0}^1 (1-t_1)^{\frac{a+p_1}{2}-2} (1-t_1 R^2_1)^{-\frac{n-3}{2}} dt_1 }{ \int_0^1 (1-t_1)^{\frac{a+p_1}{2}-2} (1-t_1 R^2_1)^{-\frac{n-1}{2}} dt_1} \nonumber \\
&=& \frac{||\bfy||^2}{n-3} \; \; \frac{_2F_1 \left(\frac{n-3}{2},1;\frac{a+p_1}{2};R_1^2 \right)}{_2F_1 \left(\frac{n-1}{2},1;\frac{a+p_1}{2};R_1^2 \right)} \nonumber
\end{eqnarray}

\noindent
To show that \eqref{eq:cor4.2ineq} holds, define the pdfs $f_1$ and $f_2$  as
\begin{eqnarray*}
f_1(t_m) &\propto& \int_{(0,1)^{k-1}} \left[ \prod_{i=1}^k (1-t_i)^{\frac{a+p_i}{2}-2} \right] (1-\sum_{i=1}^k t_i R_i^2)^{-\frac{n-1}{2}} d\bft_{-m} \\
\mbox{ and } f_2(t_m) &\propto& \int_{(0,1)^{k-1}} \left[ \prod_{i=1}^k (1-t_i)^{\frac{a+p_i}{2}-2} \right] (1- t_1 R_1^2)^{-\frac{n-1}{2}} d\bft_{-m}
\end{eqnarray*}
and use Lemma \ref{lem4.2} to obtain the result
\[  E_{f_1}(1- \sum_{i=1}^k t_i R^2_i) = 1 - \sum_{i=1}^k R^2_i  E_{f_1} (t_i) \leq 1 - \sum_{i=1}^k R^2_i  E_{f_2} (t_i) \leq 1 -  R^2_1  E_{f_2} (t_1) =   E_{f_2} ( 1 -   t_1 R^2_1 )
\]

Note that for the ordinary hyper-$g$ prior, the inequality is replaced by the equality
\[  E(\sigma^2 \mid \bfy) = \frac{||\bfy||^2}{n-3} \; \; \frac{_2F_1 \left(\frac{n-3}{2},1;\frac{a+p}{2};R^2 \right)}{_2F_1 \left(\frac{n-1}{2},1;\frac{a+p}{2};R^2 \right)} \]

\noindent
We shall use the following identity for Gaussian hypergeometric functions to simplify the RHS of the inequality (equality under the ordinary hyper-$g$):
\begin{eqnarray*}
\lim\limits_{z \rightarrow 1} (1-z)^{a+b-c}\; _2F_1\left( a,b;c;z \right) = \frac{\Gamma(a+b-c)\Gamma(c)}{\Gamma(a)\Gamma(b)}
\end{eqnarray*}
when $a+b-c > 0$.\\
As $N \rightarrow \infty, ||y|| \rightarrow \infty \; \mbox{ and } \; R_1^2 \rightarrow 1$. Thus,
\begin{eqnarray*}
&& \lim\limits_{N\rightarrow \infty} \frac{||\bfy||^2}{n-3} \; \; \frac{_2F_1 \left(\frac{n-3}{2},1;\frac{a+p_1}{2};R_1^2 \right)}{_2F_1 \left(\frac{n-1}{2},1;\frac{a+p_1}{2};R_1^2 \right)} \\
&=& \lim_{N\rightarrow \infty} \frac{||\bfy||^2 (1-R_1^2)}{n-3} \frac{ (1-R^2_1)^{\frac{n-3}{2}+1-\frac{a+p_1}{2}}\; _2F_1\left(\frac{n-3}{2},1;\frac{a+p_1}{2};R_1^2 \right)}{ (1-R^2_1)^{\frac{n-1}{2}+1-\frac{a+p_1}{2}}\; _2F_1\left(\frac{n-1}{2},1;\frac{a+p_1}{2};R_1^2 \right)} \\
&=&  \frac{(n-3)/2}{(n-1-a-p_1)/2}  \times \lim_{N\rightarrow \infty} \frac{||\bfy||^2 (1-R_1^2)}{n-3} \;,\mbox{ provided }  n > a+p_1+1 
\end{eqnarray*}
\begin{eqnarray*}
\mbox{ So, } \lim_{N \rightarrow \infty} E(\sigma^2 \mid \bfy) &\leq&  \lim_{N \rightarrow \infty} \frac{||\bfy||^2(1-R_1^2)}{n-1-a-p_1} \\
&=& \frac{1}{n-1-a-p_1} \left[ (n-p-1) \widehat{\sigma}^2 + \sum_{i=2}^k (X_i \widehat{\bfb_i})^T (X_i \widehat{\bfb_i}) \right] 
\end{eqnarray*}

\subsection{ Proof of Lemma \ref{lem5.1}}  \label{appa7}

We follow the same notation as in Theorem \ref{thm5.2}. For any $i \in B_T$,  
\begin{eqnarray*}
&& \int_{(0,1)^{k_T}} \frac{g_i}{1+g_i} \; \pi(\bfg \mid \Mt, \bfy) d\bfg \\
 &=& \frac{ \int_{(0,1)^{k_T}} t_i \prod_{j=1}^{k_T} (1-t_j)^{\frac{a+p_{j,T}}{2}-2}(1- \sum_{j=1}^{k_T} t_j R_{j,T}^2)^{-\frac{n-1}{2}} d\bft }{ \int_{(0,1)^{k_T}} \prod_{j=1}^{k_T} (1-t_j)^{\frac{a+p_{j,T}}{2}-2}(1- \sum_{j=1}^{k_T} t_j R_{j,T}^2)^{-\frac{n-1}{2}} d\bft } \\
&\geq& \frac{ \int_{(0,1)^{k_T}} t_i \prod_{j=1}^{k_T} (1-t_j)^{\frac{a+p_{j,T}}{2}-2}(1- t_i R_{i,T}^2)^{-\frac{n-1}{2}} d\bft}{ \int_{(0,1)^{k_T}} \prod_{j=1}^{k_T} (1-t_j)^ {\frac{a+p_{j,T}}{2}-2}(1- t_i R_{i,T}^2)^{-\frac{n-1}{2}} d \bft} \;  \; \mbox{(by Lemma }\ref{lem4.2})\\
&=& \frac{\int_0^1 t_i (1-t_i)^{\frac{a+p_{i,T}}{2}-2} (1-t_i R_{i,T}^2)^{-\frac{n-1}{2}}  dt_i }{\int_0^1 (1-t_i)^{\frac{a+p_{i,T}}{2}-2} (1-t_i R_{i,T}^2)^{-\frac{n-1}{2}}  dt_i} 
\end{eqnarray*}

For a specific $i$, define $m=\frac{n-1}{2}$, $b=\frac{a+p_{i,T}}{2}-2$ and $z=R_{i,T}^2$, where $0<z<1$ for all $n$ (since the predictor $X_{i,T}$ is part of the true model). Then, for that index $i$
\begin{eqnarray*} 
 E \left( \frac{g_i}{1+g_i} \mid \Mt, \bfy \right) &\geq& \frac{\int_0^1 t (1-t)^b (1-tz)^{-m} dt }{\int_0^1 (1-t)^{b} (1-tz)^{-m} dt} \\
&=& \frac{\sum_{k=0}^{\infty} \frac{\Gamma(m+k)}{\Gamma(b+k+2)} \frac{1}{1+ (b+1)/(k+1)} z^k }{ \sum_{k=0}^{\infty} \frac{\Gamma(m+k)}{\Gamma(b+k+2)} z^k } \; \;  \mbox{(see Theorem } \ref{thm3.1}) \\
\end{eqnarray*}

\noindent
Given an arbitrary $\eta > 0$ , $\exists$ $N_0$ (which does not depend on $m$) such that $\forall$ $k > N_0$, $\frac{1}{1+ \frac{b+1}{k+1}} > 1 - \eta$. So
\begin{eqnarray*} 
E \left( \frac{g_i}{1+g_i} \mid \Mt, \bfy \right) &>& \frac{ \sum_{k=0}^{N_0} \frac{\Gamma(m+k)(k+1)}{\Gamma(b+k+3)}z^k + (1-\eta) \sum_{k=N_0+1}^{\infty} \frac{\Gamma(m+k)}{\Gamma(b+k+2)}z^k  }{ \sum_{k=0}^{N_0} \frac{\Gamma(m+k)}{\Gamma(b+k+2)}z^k + \sum_{k=N_0+1}^{\infty} \frac{\Gamma(m+k)}{\Gamma(b+k+2)}z^k } \\
&=& \frac{q_1 + (1-\eta)T}{q_2 + T} = (1-\eta) + \frac{q_1 - (1- \eta)q_2}{q_2 + T} \\
\end{eqnarray*}

\noindent
To prove the lemma we have to show that for any $i \in B_T, \; \lim\limits_{m \rightarrow \infty} E \left(  \frac{g_i}{1+g_i} \mid \Mt, \bfy  \right)$ $= 1$, which is equivalent to $\lim\limits_{m \rightarrow \infty} \frac{q_1 - (1- \eta)q_2}{q_2 + T} = 0$.

\[ \frac{q_1 - (1- \eta)q_2}{q_2 + T} = \frac{[q_1 - (1- \eta)q_2]/\Gamma(m+N_0+1)}{ \frac{q_2}{\Gamma(m+N_0+1)} + \frac{z^{N_0+1}}{\Gamma(b+N_0+3)} + \sum_{k=N_0+2}^{\infty} \frac{\Gamma(m+k) z^k}{\Gamma(b+k+2) \Gamma(m+N_0+1)}  }
\]

For any $0<z<1$, $\sum_{k=N_0+2}^{\infty} \frac{\Gamma(m+k) z^k}{\Gamma(b+k+2)\Gamma(m+N_0+1)} < \infty$ for all finite $m$. Let $\rho = \lim\limits_{n \rightarrow \infty} R_{i,T}^2$, which must satisfy $0 < \rho < 1$ since the collection of predictors $X_{i,T}$ is part of the true model (see Lemma \ref{lema1} in Appendix \ref{appa6}). With a large enough $m$ it is possible to make $\frac{q_1 - (1- \eta)q_2}{\Gamma(m+N_0+1)}$ arbitrarily small while the denominator above is a finite number exceeding  $\frac{(\rho-\Phi)^{N_0+1}}{\Gamma(b+N_0+3)}$, where $\Phi$ is some small positive number (smaller than $\rho$).

Thus, we can find a large enough number $M$ so that given any arbitrary $\eta>0$ and $\delta>0$, whenever $m \geq M$
\begin{eqnarray*}
E \left( \frac{g_i}{1+g_i} \mid \Mt, \bfy \right) &>& 1- \eta - \delta \\
\implies \lim_{n \rightarrow \infty} E \left( \frac{g_i}{1+g_i} \mid \Mt, \bfy \right) &=& 1 ,\; \mbox{for any } i \in B_T
\end{eqnarray*}

\subsection{Finding the Limit of \eqref{eq:robustlim} in Appendix \ref{appenb3}} \label{appb2}

\begin{eqnarray*}
\mbox{ Let } P =  \frac{I_1}{I_2} = \frac{\int_{0}^1 t^{(p-1)/2} \left( B+ t \right)^{-(n-1)/2}  dt }{\int_{0}^1 t^{(p_1-1)/2} \left( B^* + t \right)^{-(n-1)/2} dt } 
\end{eqnarray*}
For any fixed $B^* > 0$ ($z < 1$), $I_2$ is finite and so $f_{B^*} (t) = \frac{1}{C_{B^*}} t^{(p_1-1)/2} (B^*+t)^{-(n-1)/2}$ is a proper density function with $C_{B^*} = I_2$ being the normalizing constant.

\begin{eqnarray*}
P &=& \frac{1}{C_{B^*}} \int_0^1 t^{(p-1)/2} \left[  \frac{p_1+1}{p+1} B^* + t \right]^{-(n-1)/2} dt \\
&=& E_{f_{B^*}} \left( t^{(p-p_1)/2} \left[ \frac{B^* + t}{\frac{p_1+1}{p+1}B^* + t} \right]^{(n-1)/2} \right) \\
&=& P_{f_{B^*}} (t < \epsilon) E_{f_{B^*}} \left( t^{(p-p_1)/2} \left[ \frac{B^* + t}{\frac{p_1+1}{p+1}B^* + t} \right]^{(n-1)/2} \mid t < \epsilon \right) \\
&&  +  P_{f_{B^*}} (t \geq \epsilon) E_{f_{B^*}} \left( t^{(p-p_1)/2} \left[ \frac{B^* + t}{\frac{p_1+1}{p+1}B^* + t} \right]^{(n-1)/2} \mid t \geq \epsilon \right) \; \; (\mbox{for a small } \epsilon > 0)\\
&<&  P_{f_{B^*}} (t < \epsilon)  \left[ \epsilon^{(p-p_1)/2} \left(\frac{p+1}{p_1+1} \right)^{(n-1)/2} \right] +  P_{f_{B^*}} (t \geq \epsilon) \left[ \left( \frac{B^* + \epsilon}{\frac{p_1+1}{p+1}B^* + \epsilon}\right)^{(n-1)/2}  \right]
\end{eqnarray*}
The last inequality follows from the fact that $t^{(p-p_1)/2}$ is increasing in $t$  while  $ \frac{B^* + t}{\frac{p_1+1}{p+1}B^* + t} $ is decreasing in $t$ on the interval $(0,1)$. It can also be shown that for any integer $k$ so that $k > \frac{p}{p_1}$, $ \frac{B^* + \epsilon}{\frac{p_1+1}{p+1}B^* + \epsilon} < k$ and hence,
\[ P \;< \; P_{f_{B^*}} (t < \epsilon)  \left[ \epsilon^{(p-p_1)/2} \left(\frac{p+1}{p_1+1} \right)^{(n-1)/2} \right] +  P_{f_{B^*}} (t \geq \epsilon) k^{(n-1)/2}
\]
We will show that $ P_{f_{B^*}} (t < \epsilon) \rightarrow 1$ as $B^* \rightarrow 0$ ($N \rightarrow \infty$) if $n > p_1+2$. For any arbitrary $0 < \epsilon < 1$,
\[ \lim_{N \rightarrow \infty} P \;  \leq \;  \epsilon^{(p-p_1)/2} \left(\frac{p+1}{p_1+1} \right)^{(n-1)/2} \; +  \; 0 = \epsilon^{(p-p_1)/2} \left(\frac{p+1}{p_1+1} \right)^{(n-1)/2}
\]

\noindent
It is possible to choose $\epsilon$ arbitrarily small, making the upper bound for the limit of $P$ arbitrarily small, and so $P \rightarrow 0$ as claimed in \eqref{eq:robustlim}.

To show that $\lim\limits_{N \rightarrow \infty}  P_{f_{B^*}} (t < \epsilon) =1$ for any $0 < \epsilon <1$, consider
\begin{eqnarray*}
\frac{ P_{f_{B^*}} (t < \epsilon) }{ P_{f_{B^*}} (t \geq \epsilon) } &=& \frac{ \frac{1}{C_{B^*}} \int_0^{\epsilon} f_{B^*}(t) dt }{ \frac{1}{C_{B^*}} \int^1_{\epsilon} f_{B^*}(t) dt }  \\
&\geq& \frac{ \int_{\epsilon^{j+1}}^{\epsilon^j} t^{(p_1-1)/2} (B^*+t)^{-(n-1)/2}dt }{ \int_{\epsilon}^1 t^{(p_1-1)/2}(B^*+t)^{-(n-1)/2}dt } \;, \; \mbox{for some } j \geq 2 \\
&\geq& \frac{ (\epsilon^j - \epsilon^{j+1}) \inf\limits_{\epsilon^{j+1}<t<\epsilon^j} t^{(p_1-1)/2} (B^*+t)^{-(n-1)/2} }{ (1-\epsilon) \sup\limits_{\epsilon < t < 1} t^{(p_1-1)/2} (B^*+t)^{-(n-1)/2} } 
\end{eqnarray*}
For $B^* < \frac{\epsilon^{j+1}(n-p_1)}{p_1-1}$, the above term equals \; $\epsilon^j \frac{(\epsilon^j)^{\frac{p_1-1}{2}} (B^* + \epsilon^j)^{-(n-1)/2}  }{ \epsilon^{\frac{p_1-1}{2}} (B^* + \epsilon)^{-(n-1)/2} }$\;,  since the function $ t^{(p_1-1)/2} (B^*+t)^{-(n-1)/2}$ is decreasing on $(\epsilon^{j+1},1)$. So,
\begin{eqnarray*}
 \lim_{N \rightarrow \infty} \frac{ P_{f_{B^*}} (t < \epsilon) }{ 1 - P_{f_{B^*}} (t < \epsilon) } &\geq& \epsilon^j  \; \epsilon^{(j-1)\frac{p_1-1}{2}} \epsilon^{-(j-1)(n-1)/2}  \; \; \; (\mbox{since } B^* \downarrow 0)\\
&=& \epsilon^{\frac{1}{2} [2j + (j-1)(p_1-n)] }
\end{eqnarray*}
This relation holds for any $j \geq 2$ and hence also for the limit as $j \rightarrow \infty$.  Since $0<\epsilon<1$, $\lim\limits_{j \rightarrow \infty} \epsilon^{\frac{1}{2} [2j + (j-1)(p_1-n)] } = \epsilon^{(n-p_1)/2} \lim\limits_{j \rightarrow \infty} \epsilon^{\frac{j}{2} [2 + p_1-n] } = \infty$ when $n>p_1+2$. This shows that $\lim\limits_{N \rightarrow \infty}  P_{f_{B^*}}(t < \epsilon) = 1$, for any $0 < \epsilon < 1$ if $n > p_1+2$, completing the argument for the convergence of \eqref{eq:robustlim} to zero.

\subsection{Proof of Theorem \ref{thm5.2} when $2 < a  \leq 3$} \label{appb5}

We rely on Conditions \ref{cond5.1} and \ref{cond5.2} in this proof through the use of Lemma \ref{lema1}. When $\Mgam \neq \Mo$, as in the proof of Theorem \ref{thm5.2}  for $3< a \leq 4$,
\begin{eqnarray*}
BF(\Mgam : \Mt) &=& \left( \frac{a-2}{2} \right)^{k_{\gamma}-k_T} \frac{\int \prod\limits_{i \in B_{\gamma}} (1-t_i)^{\frac{a+p_{i,\gamma}}{2}-2} (1-\sum\limits_{i \in B_{\gamma}} t_i R_{i,\gamma}^2)^{-\frac{n-1}{2}} d\bft }{ \int \prod\limits_{i \in B_T} (1-t_i)^{\frac{a+p_{i,T}}{2}-2} (1-\sum\limits_{i \in B_T} t_i R_{i,T}^2)^{-\frac{n-1}{2}} d\bft }  \\
&=& \left( \frac{a-2}{2} \right)^{k_{\gamma}-k_T} \frac{\int \prod\limits_{i \in J_{\gamma}} (1-t_i)^{-\frac{1}{2}}  h_{\gamma}(\bft)  d\bft}{\int  \prod\limits_{i \in J_T} (1-t_i)^{-\frac{1}{2}}   h_{T}(\bft)  d\bft} 
\end{eqnarray*}
where  $h_{j} (\bft) = \prod\limits_{i \in J_j} \left[ (1-t_i)^{\frac{a+p_{i,j}+1}{2}-2} \right] \prod\limits_{i \in B_j \backslash J_j}  \left[ (1-t_i)^{\frac{a+p_{i,j}}{2}-2 } \right] (1-\sum\limits_{i \in B_j} t_i R_{i,j}^2)^{-\frac{n-1}{2}}$ = $\prod\limits_{i \in B_j} \left[ (1-t_i)^{\frac{a+p_{i,j}^*}{2}-2} \right]$ $(1-\sum\limits_{i \in B_j} t_i R_{i,j}^2)^{-\frac{n-1}{2}}$,  with $p^*_{i,j}=p_{i,j}+1$ when $i \in J_j$ and $p^*_{i,j}=p_{i,j}$ otherwise. Here both $J_{\gamma} \subseteq B_{\gamma}$ and $J_T \subseteq B_T$ are the block indices corresponding to $a+p_{i,j} \leq 4$ which happens only when $2 < a \leq 3$ and the related $p_{i,j} = 1$; $j \in \{\gamma, T\}$.

It is clear that $b_{i,j}^*=\frac{a+p_{i,j}^*}{2}-2 > 0$ $\forall \; i$, as in the situation described in Appendix \ref{appa6}. Using a variation of the Laplace approximation, 
\begin{eqnarray}
 \lim\limits_{m \rightarrow \infty} BF(\Mgam:\Mt)  &=&  \left( \frac{a-2}{2} \right)^{k_{\gamma}-k_T}  \lim\limits_{m \rightarrow \infty}  \frac{|H_T (\hat{\bft}_{T})|^{1/2}}{|H_{\gamma} (\hat{\bft}_{\gamma})|^{1/2}} \frac{ \prod\limits_{i \in J_{\gamma}} (1- \hat{t}_{i,\gamma})^{-1/2} }{ \prod\limits_{i \in J_T} (1- \hat{t}_{i,T})^{-1/2} }  \nonumber \\
&\times&    \frac{\exp[ \sum\limits_{i \in B_{\gamma}} b_{i,\gamma}^* \log(1- \hat{t}_{i,\gamma}) - m \log(1-\sum\limits_{i \in B_{\gamma}} \hat{t}_{i,\gamma} R_{i,\gamma}^2 ) ]}{\exp[ \sum\limits_{i \in B_T} b^*_{i,T} \log(1- \hat{t}_{i,T}) - m \log(1-\sum\limits_{i \in B_{T}} \hat{t}_{i,T} R_{i,T}^2 ) ]}   \label{eq:BFal3}
\end{eqnarray}
with the same definitions for $H_j (\bft_j)$, $\hat{t}_{i,j}$ and $m$  as in  Appendix \ref{appa6}.

If we can show that \eqref{eq:BFal3} behaves exactly like the Laplace approximation to \eqref{eq:BFlap} for any model $\Mgam \neq \Mt$, the proof is complete. Observe that the difference between the approximations to  \eqref{eq:BFlap} and \eqref{eq:BFal3} comes just from the extra term $\frac{ \prod\limits_{i \in J_{\gamma}} (1- \hat{t}_{i,\gamma})^{-1/2} }{\prod\limits_{i \in J_T} (1- \hat{t}_{i,T})^{-1/2}} $ which equals

\begin{eqnarray} \label{eq:extra} 
 \frac{  \prod\limits_{i \in J_{\gamma}} [\frac{b_{i,\gamma} (1-R^2_{\gamma})}{R^2_{i,\gamma}(m-b_{\gamma})} ]^{-1/2} }{\prod\limits_{i \in J_T}  [ \frac{b_{i,T} (1-R^2_T)}{R^2_{i,T}(m-b_T)} ]^{-1/2} } =  \left[ \frac{ \prod\limits_{i \in J_{\gamma}} O(m R^2_{i,\gamma}) }{ \prod\limits_{i \in J_T} O(m R^2_{i,T}) } \right]^{\frac{1}{2}}
\end{eqnarray}

\noindent
\underline{Case 1}: $\Mgam \not \supset \Mt$

\noindent
In this case, the remaining term in the Laplace approximation to the integral goes to zero in probability at an exponential rate (see Appendix \ref{appa6}) while the extra part will either be bounded, go to infinity or go to zero in probability at a polynomial rate. For all models $\Mgam$  belonging to Case 1, \eqref{eq:BFal3} equals zero.\\

\noindent
\underline{Case 2}: $\Mgam \supset \Mt$

\noindent
\underline{Case 2A}:  For Case 2A of Theorem \ref{thm5.2}, $B_{\gamma} = B_T$ and $\Mt \subset M_{\gamma}$ which means that whenever $p_{i,\gamma}=1$, we must have $p_{i,T} = 1$ whereas it is possible to have $p_{i,T}=1$ along with $p_{i,\gamma} > 1$. This indicates that $|J_T| \geq |J_{\gamma}|$ for this group of models. Also note that all $R^2_{i,j} \stackrel{P}{\rightarrow} C$ for some non-zero constant $C$ (by Lemma \ref{lema1} (i) from Appendix \ref{appa6}). Hence \eqref{eq:extra} reduces to $O(m^{\frac{|J_{\gamma}| - |J_T|}{2} })$ while the other part of the integral is of the order $O(m^{\frac{p_T^*-p_{\gamma}^*}{2}}) = O(m^{\frac{p_T-p_{\gamma}+|J_T|-|J_{\gamma}|}{2}})$. This implies that \eqref{eq:BFal3} equals  $\lim\limits_{m \rightarrow \infty} O(m^{p_T - p_{\gamma}}) = 0$, as in Case 2A from Appendix \ref{appa6}.

\noindent
\underline{Case 2B}:   For Case 2B of Theorem \ref{thm5.2}, $B_{\gamma} \supset B_T$ and $\Mt \subset M_{\gamma}$. As before for $i \in B_T \bigcap  B_{\gamma} = B_T$, whenever $p_{i,\gamma}=1$, we must have $p_{i,T} = 1$ whereas it is possible to have $p_{i,T}=1$ along with $p_{i,\gamma} > 1$. This translates to the inequality $|J_{\gamma} \bigcap B_T| \leq |J_T|$. Whereas for $i \in B_{\gamma}\backslash B_T$, it might happen that $p_{i,\gamma}=1$, but the corresponding $p_{i,T}=0$ since the block does not exist. However, $R_{i,\gamma}^2 \rightarrow 0$ and $nR_{i,\gamma}^2 \stackrel{d}{\rightarrow} c \chi^2_{p_{i,\gamma}}$  for $i \in B_{\gamma} \backslash B_T$  (by Lemma \ref{lema1} (ii) and (iii)) so that $m R_{i,\gamma}^2 = O(1)$ for all $i \in J_{\gamma} \backslash B_T$. Thus \eqref{eq:extra} becomes $O(m^{ \frac{|J_{\gamma} \bigcap B_T| - |J_T|}{2}} )$  and the other part of the integral is $O\Big(m^{\frac{1}{2} \sum\limits_{i \in B_T} [ p_{i,T}^*-p_{i,\gamma}^* ]} \Big) = O\Big(m^{\frac{1}{2} \sum\limits_{i \in B_T} [ p_{i,T}-p_{i,\gamma} ] + \frac{1}{2}[ \; |J_T| - |J_{\gamma} \bigcap B_T| \;]} \Big)$. \eqref{eq:BFal3} is zero since it reduces to $\lim\limits_{m \rightarrow \infty} O\Big(m^{\frac{1}{2} \sum\limits_{i \in B_T} [ p_{i,T}-p_{i,\gamma} ]} \Big)=0$ and we have the same result as in Case 2B from Appendix \ref{appa6}.

\noindent
\underline{Case 2C}:   For Case 2C of Theorem \ref{thm5.2}, we know that there is a one to one equivalence between $p_{i,T}=1$ and $p_{i,\gamma}=1$ for $i \in B_T$ since the predictors are exactly the same in all blocks common to $\Mgam$ and $\Mt$. This means that $|J_{\gamma} \bigcap B_T| = |J_T|$, while as before, $m R_{i,\gamma}^2 = O(1)$ for all $i \in J_{\gamma} \backslash B_T$. Now \eqref{eq:extra} is of the order $O(m^{ \frac{|J_{\gamma} \bigcap B_T| - |J_T|}{2}} ) = O(1)$  and the other term in the Laplace approximation is also $O(1)$, leading to the same conclusion  as in Case 2C from Appendix \ref{appa6} that the limit of the Bayes factor does not equal zero with probability 1.

The situation when $\Mgam = \Mo$ is similar to Case 2C here and in Theorem \ref{thm5.2} with $B_T = \phi$. By the same reasoning, the Bayes factor does not converge to zero in probability.

\subsection{Detailed Proof of Theorem \ref{thm3.1}} \label{appb11}

The posterior mean of the regression coefficients is
\[ \widehat{\bfb} = E \left(\frac{g}{1+g} \mid \bfy \right) \widehat{\bfb}_{LS}
\]
where $\widehat{\bfb}$ denotes the posterior mean of the regression coefficient $\bfb$  and $\widehat{\bfb}_{LS}$ denotes the estimate of $\bfb$ under least squares.

For the hyper-$g$ prior, the posterior expectation of the shrinkage factor can be expressed in terms of $R^2$ \citep[see][]{liang2008mixtures}

\[ E \left(\frac{g}{1+g} \mid \bfy \right) = \frac{2}{p+a} \frac{_2F_1(\frac{n-1}{2},2;\frac{p+a}{2}+1;R^2)} {_2F_1(\frac{n-1}{2},1;\frac{p+a}{2};R^2)} ,
\]
where $_2F_1$ is the Gaussian Hypergeometric Function. $_2F_1(a,b;c;z)$ is finite for $|z| < 1 \;\mbox{whenever} \; c>b>0$. Here, $c-b=(p+a)/2 - 1 > 0$  since $2< a \leq 4$ and $p>0$. Thus, for all values of $R^2 < 1$, both numerator and denominator are finite. We use an integral representation of the $_2 F_1$ function.
\begin{eqnarray*}
 E \left(\frac{g}{1+g} \mid \bfy \right) &=& \frac{2}{p+a} \frac{_2F_1(\frac{n-1}{2},2;\frac{p+a}{2}+1;R^2)} {_2F_1(\frac{n-1}{2},1;\frac{p+a}{2};R^2)} \\
 &=& \frac{ \int_0^1 t(1-t)^{\frac{p+a}{2}-2} (1-tR^2)^{-\frac{n-1}{2}} dt} { \int_0^1 (1-t)^{\frac{p+a}{2}-2} (1-tR^2)^{-\frac{n-1}{2}} dt} \\
\end{eqnarray*}
Define $m=\frac{n-1}{2}$, $b=\frac{p+a}{2}-2$ and $z=R^2 (\leq 1)$ so that we have
\[ E \left(\frac{g}{1+g} \mid \bfy \right) = \frac{ \int_0^1 t(1-t)^b (1-tz)^{-m} dt} { \int_0^1 (1-t)^b (1-tz)^{-m} dt} 
\]

In our problem we have $m>b> -\frac{1}{2}$ since $a > 2$ and $p \geq 1$ which satisfies the requirement of $b > -1$ required later in the proof.
\begin{eqnarray*}
\mbox{Numerator} &=& \int_0^1 t(1-t)^b \left[ \sum_{k=0}^{\infty} {m + k-1 \choose k} (tz)^k \right] dt \\
 &=& \int_0^1 \sum_{k=0}^{\infty} {m + k-1 \choose k} z^k t^{k+1} (1-t)^b dt \\
 &=& \sum_{k=0}^{\infty}   {m + k-1 \choose k} z^k \int_0^1 t^{k+1} (1-t)^b dt  \\
&& \hspace{1 in} (\mbox{the infinite sum converges for all } |z|<1) 
 \end{eqnarray*}

\begin{eqnarray*}
 &=& \sum_{k=0}^{\infty} \frac{(k+1) \Gamma(m+k) \Gamma(b+1)}{\Gamma(m) \Gamma(b+k+3)} z^k  \; \;  \; (\mbox{need } b > -1)
\end{eqnarray*}

Similarly we can show by interchanging the positions of the infinite sum and the integral that

\begin{eqnarray*}
 \mbox{Denominator} &=& \int_0^1 (1-t)^b \left[ \sum_{k=0}^{\infty} {m + k-1 \choose k} (tz)^k \right] dt  \\
 &=& \sum_{k=0}^{\infty}   {m + k-1 \choose k} z^k \int_0^1 t^k (1-t)^b dt = \sum_{k=0}^{\infty} \frac{ \Gamma(m+k) \Gamma(b+1)}{\Gamma(m) \Gamma(b+k+2)} z^k  .
 \end{eqnarray*}
Thus,
\[E \left(\frac{g}{1+g} \mid \bfy \right) = \frac{\sum_{k=0}^{\infty} \frac{(k+1) \Gamma(m+k) \Gamma(b+1)}{\Gamma(m) \Gamma(b+k+3)} z^k} {\sum_{k=0}^{\infty} \frac{\Gamma(m+k) \Gamma(b+1)}{\Gamma(m) \Gamma(b+k+2)} z^k} 
=  \frac{ \sum_{k=0}^{\infty} \frac{\Gamma(m+k)}{\Gamma(b+k+2)} \frac{1}{1+\frac{b+1}{k+1}} z^k  }{ \sum_{k=0}^{\infty} \frac{\Gamma(m+k)}{\Gamma(b+k+2)}z^k }  \]

\noindent
When $m=\frac{n-1}{2} > b+2 = \frac{p+a}{2}$ , we show that $\frac{\Gamma(m+k)}{\Gamma(b+k+2)}$ is increasing in $k$. Consider the function $D(k) = \log \Gamma (m+k) - \log \Gamma (b+k+2)$ which has the derivative, where $\Psi(\cdot)$ is the digamma function
\begin{eqnarray*}
D^{\prime} (k) &=& \Psi(m+k) - \Psi(b+k+2) \\
&=& \int_0^{\infty} \left[ \frac{e^{-t}}{t} - \frac{e^{-(m+k)t}}{1-e^{-t}} \right] dt - \int_0^{\infty} \left[ \frac{e^{-t}}{t} - \frac{e^{-(b+k+2)t}}{1-e^{-t}} \right] dt \\
&=&  \int_0^{\infty} \left[ \frac{e^{-(b+k+2)t} - e^{-(m+k)t}}{1-e^{-t}} \right] dt  > 0 
\end{eqnarray*}
$\forall \; k \in \mathbb{N}$,  whenever  $m > b+2$.

This implies that $D(k)$ is increasing in $k$ and so is $\exp(D(k))$. The algebra above makes use of the following standard integral representation of the digamma function
\[ \Psi(x) =   \int_0^{\infty} \left[ \frac{e^{-t}}{t} - \frac{e^{-tx}}{1-e^{-t}} \right] dt \; \;,\; \mbox{ when } x > 0 
\]

Lemma \ref{lem3.1} states that $R^2 \rightarrow 1$ as $N \rightarrow \infty$ so that $  \lim\limits_{N \rightarrow \infty} \widehat{\bfb} = \lim\limits_{R^2 \rightarrow 1} E \left(\frac{g}{1+g} \mid \bfy \right) \widehat{\bfb}_{LS} $. 

\noindent
The proof will be complete if we can show that $\lim\limits_{R^2 \rightarrow 1} E \left(\frac{g}{1+g} \mid \bfy \right) = 1$. Throughout, we make use of the expression
\begin{eqnarray*}
\lim_{R^2 \rightarrow 1} E \left(\frac{g}{1+g} \mid \bfy \right) &=& \lim_{z\uparrow 1} \frac{ \sum_{k=0}^{\infty} \frac{\Gamma(m+k)}{\Gamma(b+k+2)} \frac{1}{1+\frac{b+1}{k+1}} z^k  }{ \sum_{k=0}^{\infty} \frac{\Gamma(m+k)}{\Gamma(b+k+2)}z^k } \\
\end{eqnarray*}

\newpage
\noindent
\underline{Case 1}: $n > p+a+1$ 

\noindent
First note that $\frac{1}{1 + \frac{b+1}{k+1}}$ is increasing in k and $\uparrow 1$ as $k \rightarrow \infty$. So for any $\eta >0, \; \exists \; N_0$ such that $\forall \; k > N_0, \frac{1}{1 + \frac{b+1}{k+1}} > 1-\eta$. 

Hence, $ \lim\limits_{z \uparrow 1} \frac{ \sum_{k=0}^{\infty} \frac{\Gamma(m+k)}{\Gamma(b+k+2)} \frac{1}{1+\frac{b+1}{k+1}} z^k  }{ \sum_{k=0}^{\infty} \frac{\Gamma(m+k)}{\Gamma(b+k+2)}z^k } $
\begin{eqnarray*}
&=& \lim_{z \uparrow 1} \frac{ \sum_{k=0}^{N_0} \frac{\Gamma(m+k)}{\Gamma(b+k+2)} \frac{1}{1+\frac{b+1}{k+1}} z^k +  \sum_{k=N_0+1}^{\infty} \frac{\Gamma(m+k)}{\Gamma(b+k+2)} \frac{1}{1+\frac{b+1}{k+1}} z^k }{ \sum_{k=0}^{N_0} \frac{\Gamma(m+k)}{\Gamma(b+k+2)}z^k + \sum_{k=N_0+1}^{\infty} \frac{\Gamma(m+k)}{\Gamma(b+k+2)}z^k } \\
&>& \lim_{z\uparrow 1} \frac{q_1 + (1-\eta)\sum_{k=N_0+1}^{\infty} \frac{\Gamma(m+k)}{\Gamma(b+k+2)}z^{k} }{q_2 + \sum_{k=N_0+1}^{\infty} \frac{\Gamma(m+k)}{\Gamma(b+k+2)}z^{k} } \\
&=& (1-\eta) + \lim_{z \uparrow 1} \frac{q_1 - (1-\eta)q_2}{q_2 + \sum_{k=N_0+1}^{\infty} \frac{\Gamma(m+k)}{\Gamma(b+k+2)}z^{k}} \\
&\geq& (1-\eta) + 0 = (1-\eta) 
\end{eqnarray*}
\begin{equation*}
\implies \lim_{R^2 \rightarrow 1} E \left(\frac{g}{1+g} \mid \bfy \right) = \lim_{z\uparrow 1} \frac{ \sum_{k=0}^{\infty} \frac{\Gamma(m+k)}{\Gamma(b+k+2)} \frac{1}{1+\frac{b+1}{k+1}} z^k  }{ \sum_{k=0}^{\infty} \frac{\Gamma(m+k)}{\Gamma(b+k+2)}z^k } = 1 \hspace{2 in}
\end{equation*}
Note that $q_1$ and $q_2$ are finite numbers corresponding to the finite sums of the first $ N_0$ terms. Also $\frac{\Gamma(m+k)}{\Gamma(b+k+2)} \rightarrow \infty$ as $k \rightarrow \infty$ due to which $\sum_{k=0}^{\infty} \frac{\Gamma(m+k)}{\Gamma(b+k+2)}$ and hence the denominator goes to infinity causing the second term above to vanish in the limit. 

\noindent
\underline{Case 2}: $p+a-1 \leq n \leq p+a+1$ 

\noindent
Let $n=p+a-1+ 2\xi$, where $0 \leq \xi \leq 1$ and define $N_0$ as in Case 1,

\begin{eqnarray*}
E \left(\frac{g}{1+g} \mid \bfy \right)  &=& \frac{ \sum_{k=0}^{\infty} \frac{\Gamma(b+k+1+\xi)}{\Gamma(b+k+2)} \frac{1}{1+\frac{b+1}{k+1}} z^k  }{ \sum_{k=0}^{\infty} \frac{\Gamma(b+k+1+\xi)}{\Gamma(b+k+2)}z^k } \\
\lim_{R^2 \rightarrow 1} E \left(\frac{g}{1+g} \mid \bfy \right) &=& \lim_{z \uparrow 1} \frac{ \sum_{k=0}^{\infty} \frac{\Gamma(b+k+1+\xi)}{\Gamma(b+k+1)(b+k+1)} \frac{1}{1+\frac{b+1}{k+1}} z^k}{ \sum_{k=0}^{\infty} \frac{\Gamma(b+k+1+\xi)}{\Gamma(b+k+1)(b+k+1)} z^k } \\
\end{eqnarray*}

\noindent
Proceeding as in Case 1, we can show that

\begin{eqnarray*}
\lim_{R^2 \rightarrow 1} E \left(\frac{g}{1+g} \mid \bfy \right) &>& (1-\eta) + \lim_{z \uparrow 1} \frac{q_1 - (1-\eta)q_2}{q_2 + \sum_{k=N_0+1}^{\infty} \frac{\Gamma(b+k+1+\xi)}{\Gamma(b+k+1)(b+k+1)}z^{k}} \\
&\geq& (1-\eta) \;,\; \mbox{for any} \; \eta>0 
\end{eqnarray*}

\noindent
As $z \uparrow 1$, the denominator of the second term becomes $\frac{q_1 - (1-\eta)q_2}{q_2 + \sum_{k=N_0+1}^{\infty} \frac{\Gamma(b+k+1+\xi)}{\Gamma(b+k+1)(b+k+1)} }$ which tends to zero if the infinite sum $\sum_{k=N_0+1}^{\infty} \frac{\Gamma(b+k+1+\xi)}{\Gamma(b+k+1)(b+k+1)}$ diverges. It does, since $\frac{\Gamma(b+k+1+\xi)}{\Gamma(b+k+1)(b+k+1)} = O(k^{-\lambda})$, where $0 \leq \lambda \leq 1$ and $ \sum_{k=N_0+1}^{\infty} O(k^{-\lambda}) = \infty $ for any $0 \leq \lambda \leq 1$. Thus,

\[ \lim_{R^2 \rightarrow 1} E \left(\frac{g}{1+g} \mid \bfy \right) = 1. \]

\noindent
\underline{Case 3}: $n < p+a-1$ (proving necessity of the constraint)

\begin{eqnarray*}
\lim_{R^2 \rightarrow 1} E \left(\frac{g}{1+g} \mid \bfy \right) &=& \lim_{R^2 \rightarrow 1} \frac{ \int_0^1 t(1-t)^{\frac{p+a}{2}-2} (1-tR^2)^{-\frac{n-1}{2}} dt} { \int_0^1 (1-t)^{\frac{p+a}{2}-2} (1-tR^2)^{-\frac{n-1}{2}} dt} \\
&=& \frac{ \int_0^1 t(1-t)^{\frac{p+a}{2}-2-\frac{n-1}{2}} dt} { \int_0^1 (1-t)^{\frac{p+a}{2}-2-\frac{n-1}{2}} dt} \\
&=&  \frac{ \mbox{Beta} (2,\frac{p+a-n-1}{2}) }{ \mbox{Beta} (1,\frac{p+a-n-1}{2})}  =    \frac{2}{p+a-n+1} 
\end{eqnarray*}
which is strictly less than 1 $\forall \; n < p+a-1$ and has a minimum value of $\frac{2}{p+a}$ when $n=1$.

\subsection{Detailed Proof of Theorem \ref{thm3.2}}\label{appDetail3.2}

\cite{liang2008mixtures} show that
\begin{eqnarray*}
 BF(M_2:\Mo) &=&  \;_2F_1 (\frac{n-1}{2},1;\frac{a+p}{2}; R_{M_2}^2) \times \frac{a-2}{a+p-2} \\
 BF(M_1:\Mo) &=& \;_2F_1 (\frac{n-1}{2},1;\frac{a+p_1}{2}; R_{M_1}^2) \times \frac{a-2}{a+p_1-2}
\end{eqnarray*}
where $R_{M_i}^2$ is the coefficient of determination for model $M_i, i=1,2$.

The $g$ prior is invariant to linear transformation of $X$, and so we can work with an orthogonalized version of the design without loss of generality.  Specifically, we consider $Q_1 = X_1$ and $Q_2 = (I - P_{Q_1})X_2$, where $P_{Q_1}$ is the projection matrix for the column space of $Q_1$.  Then $X$ can be represented as $X=QT$ for a suitable upper triangular matrix $T$ and $X \bfb = Q \kappa$, where $\kappa=T\bfb$ also has a hyper-$g$ prior. Since $T$ is upper triangular, $||\bfb_1|| \rightarrow \infty$ is equivalent to $||\kappa_1|| \rightarrow \infty$ while $\kappa_2$ stays fixed in the sequence.  Under the block orthogonal setup, $ R_{M_1}^2 = \frac{(Q_1 \widehat{\bld{\kappa}}_1)^T (Q_1 \widehat{\bld{\kappa}}_1)}{\bld{y}^T \bld{y}}$ and $R_{M_2}^2 = R_{M_1}^2 + \frac{(Q_2 \widehat{\bld{\kappa}}_2)^T (Q_2 \widehat{\bld{\kappa}}_2)}{\bld{y}^T \bld{y}}$.  The term $(Q_2 \widehat{\bld{\kappa}}_2)^T (Q_2 \widehat{\bld{\kappa}}_2)$ is constant throughout the sequence $\{ \Psi_N \}$ and so $\frac{(Q_2 \widehat{\bld{\kappa}}_2)^T (Q_2 \widehat{\bld{\kappa}}_2)}{ \bld{y}^T \bld{y} } \rightarrow 0$ as $N \rightarrow \infty$.

\begin{eqnarray*}
BF(M_2:M_1) &=& \frac{BF(M_2:\Mo)}{BF(M_1:\Mo)} \\
&=& \frac{a+p_1-2}{a+p-2} .\; \frac{_2F_1 (\frac{n-1}{2},1;\frac{a+p}{2}; R_{M_2}^2)}{_2F_1 (\frac{n-1}{2},1;\frac{a+p_1}{2}; R_{M_1}^2)} \\
&=& \frac{\int_0^1 (1-t)^{\frac{a+p}{2}-2} (1-tR_{M_2}^2)^{-\frac{n-1}{2}} dt}{ \int_0^1 (1-t)^{\frac{a+p_1}{2}-2} (1-tR_{M_1}^2)^{-\frac{n-1}{2}}dt }
\end{eqnarray*}

\noindent
Define $b = \frac{a+p_1}{2}-2$, $m=\frac{n-1}{2}$, $R_{M_1}^2 =z$ and $R_{M_2}^2 =z+q$. When $||\beta_1|| \rightarrow \infty$, both $R_{M_2}^2$ and $R_{M_1}^2$ go to 1 which results in $z \uparrow 1$ and $q \downarrow 0$.

\[ BF(M_2:M_1) = \frac{\int_0^1 (1-t)^{b+\frac{p_2}{2}} \left[ 1-t(z+q) \right]^{-m} dt }{ \int_0^1 (1-t)^{b} \left[ 1-tz \right]^{-m} dt }
\]

\noindent
Proceeding as in Theorem \ref{thm3.1}, 
\[ \mbox{Numerator} = \sum_{k=0}^{\infty} \frac{\Gamma(m+k) \Gamma(b+1+\frac{p_2}{2})}{\Gamma(m) \Gamma(b+k+2+\frac{p_2}{2})} (z+q)^k
\]
and 
\[  \mbox{Denominator} = \sum_{k=0}^{\infty} \frac{\Gamma(m+k) \Gamma(b+1)}{\Gamma(m) \Gamma(b+k+2)} z^k
\]

Thus,
\[ BF(M_2:M_1) = \frac{\Gamma(b+1+\frac{p_2}{2})}{\Gamma(b+1)} \; \frac{ \sum_{k=0}^{\infty} \frac{\Gamma(m+k)}{\Gamma(b+k+2)} \left\{ \frac{\Gamma(b+k+2)}{\Gamma(b+k+2+ \frac{p_2}{2})} \right\}(z+q)^k }{ \sum_{k=0}^{\infty} \frac{\Gamma(m+k)}{\Gamma(b+k+2)} z^k }
\]

\begin{eqnarray*}
\mbox{Hence } \lim_{||\beta_1|| \rightarrow \infty} BF(M_2:M_1) &=& \lim_{\substack{z \rightarrow 1 \\ q \rightarrow 0}} BF(M_2:M_1) \\
&=& \lim_{z \rightarrow 1} \; \left\{ \lim_{q \rightarrow 0} BF(M_2:M_1) \right\}
\end{eqnarray*}

\noindent
The last step is justified when $ \lim\limits_{q \rightarrow 0} BF(M_2:M_1)$ exists for all $0\leq z < 1$. This holds since  $\lim \limits_{q \rightarrow 0} BF(M_2:M_1) = \frac{a+p_1-2}{a+p-2} \; \frac{_2F_1(m,1;b+2+\frac{p_2}{2};z)}{_2F_1(m,1;b+2;z)} $ which exists and is finite for all $0 \leq z < 1$.
\begin{eqnarray*}
\lim_{||\beta_1|| \rightarrow \infty} BF(M_2:M_1) &=& \lim_{z \uparrow 1} \frac{\Gamma(b+1+\frac{p_2}{2})}{\Gamma(b+1)} \; \frac{ \sum_{k=0}^{\infty} \frac{\Gamma(m+k)}{\Gamma(b+k+2)} \left\{ \frac{\Gamma(b+k+2)}{\Gamma(b+k+2 +\frac{p_2}{2})} \right\}z^k }{ \sum_{k=0}^{\infty} \frac{\Gamma(m+k)}{\Gamma(b+k+2)} z^k } 
\end{eqnarray*}

\noindent
But $\frac{\Gamma(b+k+2)}{\Gamma(b+k+2+ \frac{p_2}{2})}$ decreases to 0 as $k \rightarrow \infty$ (see Appendix \ref{appb11} for a proof). Hence given an arbitrary $\eta>0$, we can find a number $N_0$ such that $\forall \; k > N_0$, $\frac{\Gamma(b+k+2)}{\Gamma(b+k+2+ \frac{p_2}{2})} < \eta $. 

\begin{eqnarray*}
&& BF(M_2:M_1) =  \frac{\Gamma(b+1+\frac{p_2}{2})}{\Gamma(b+1)} \; \times \\
 && \frac{ \sum_{k=0}^{N_0} \frac{\Gamma(m+k)}{\Gamma(b+k+2)} \left\{ \frac{\Gamma(b+k+2)}{\Gamma(b+k+2+ \frac{p_2}{2})} \right\}z^k + \sum_{k=N_0+1}^{\infty} \frac{\Gamma(m+k)}{\Gamma(b+k+2)} \left\{ \frac{\Gamma(b+k+2)}{\Gamma(b+k+2+ \frac{p_2}{2})} \right\}z^k }{ \sum_{k=0}^{N_0} \frac{\Gamma(m+k)}{\Gamma(b+k+2)} z^k + \sum_{k=N_0+1}^{\infty} \frac{\Gamma(m+k)}{\Gamma(b+k+2)} z^k } \\
&<& \frac{q_1 + \eta \sum_{k=N_0+1}^{\infty} \frac{\Gamma(m+k)}{\Gamma(b+k+2)} z^k }{q_2 + \sum_{k=N_0+1}^{\infty} \frac{\Gamma(m+k)}{\Gamma(b+k+2)} z^k} \times \frac{\Gamma(b+1+\frac{p_2}{2})}{\Gamma(b+1)} 
 \end{eqnarray*}

\begin{eqnarray*}
&=& \frac{\Gamma(b+1+\frac{p_2}{2})}{\Gamma(b+1)} \; . \frac{q_1+\eta T}{q_2 + T} \;,\; \mbox{with } T = \sum_{k=N_0+1}^{\infty} \frac{\Gamma(m+k)}{\Gamma(b+k+2)} z^k  \\
&=& \frac{\Gamma(b+1+\frac{p_2}{2})}{\Gamma(b+1)} \left[ \frac{\eta(q_1+T)}{q_2+T} + \frac{(1-\eta)q_1}{q_2+T} \right]
\end{eqnarray*}

\noindent
We later show that $||\beta_1|| \rightarrow \infty$ (or $z \uparrow 1$) implies that $T \rightarrow \infty$ when $n \geq a+p_1-1$.

\begin{eqnarray*}
\implies \lim_{||\beta_1|| \rightarrow \infty} BF(M_2:M_1) &\leq& \lim_{T \rightarrow \infty} \frac{\Gamma(b+1+\frac{p_2}{2})}{\Gamma(b+1)} \left[ \frac{\eta(q_1+T)}{q_2+T} + \frac{(1-\eta)q_1}{q_2+T} \right] \\
&=& \eta \; \frac{\Gamma(b+1+\frac{p_2}{2})}{\Gamma(b+1)} 
\end{eqnarray*}

\begin{eqnarray*}
\mbox{Hence } \lim_{||\beta_1|| \rightarrow \infty} BF(M_2:M_1) &\leq& \eta \; \frac{\Gamma(b+1+\frac{p_2}{2})}{\Gamma(b+1)} \;,\; \mbox{for any arbitrary choice of } \eta >0, \\ 
\mbox{and } \lim_{||\beta_1|| \rightarrow \infty} BF(M_2:M_1) &=& 0
\end{eqnarray*}

\noindent
We now prove that, for $n \geq a+p_1-1$, $T \rightarrow \infty$ when $||\beta_1|| \rightarrow \infty$.

\noindent
\underline{Case 1}: $n > a+p_1+1$

\noindent
Then $m>b+2$ and so $\frac{\Gamma(m+k)}{\Gamma(b+k+2)} \uparrow \infty$ as $k \rightarrow \infty$. 
\[ \lim_{z \uparrow 1} T = \sum_{k=N_0+1}^{\infty} \frac{\Gamma(m+k)}{\Gamma(b+k+2)}
\]

\noindent
So $T \rightarrow \infty$ for such values of $n$. \\

\noindent
\underline{Case 2}: $a+p_1-1 \leq n \leq a+p_1+1$

\noindent
Let $n=a+p_1-1+2\xi$ where $0 \leq \xi \leq 1$. Then $m=\frac{a+p_1}{2} -1 + \xi$ and $\frac{\Gamma(m+k)}{\Gamma(b+k+2)} = \frac{\Gamma(\frac{a+p_1}{2}-1+\xi+k)}{\Gamma(\frac{a+p_1}{2}+k)}$.  
For $0 \leq \xi \leq 1$, $\frac{\Gamma(\frac{a+p_1}{2}-1+\xi+k)}{\Gamma(\frac{a+p_1}{2}+k)} =\frac{\Gamma(\frac{a+p_1}{2}-1+k+\xi)}{\Gamma(\frac{a+p_1}{2}-1+k+1)} = O(k^{-\lambda})$ for some $0 \leq\lambda \leq 1$. But $\sum_{k=N_0+1}^{\infty} O(k^{-\lambda}) = \infty$ for such values of $\lambda$ implying that $T \rightarrow \infty$ as $z \uparrow 1$. \\

\noindent
\underline{Case 3}: $n < a+p_1-1$ (proving necessity of the constraint)
\begin{eqnarray*}
\lim_{||\beta_1 \rightarrow \infty||} BF(M_2:M_1) &=& \lim_{\substack{R_{M_1}^2 \rightarrow 1 \\ R_{M_2}^2 \rightarrow 1}} \frac{\int_0^1 (1-t)^{\frac{a+p}{2}-2} (1-tR_{M_2}^2)^{-\frac{n-1}{2}} dt}{ \int_0^1 (1-t)^{\frac{a+p_1}{2}-2} (1-tR_{M_1}^2)^{-\frac{n-1}{2}}dt } \\
&=&  \frac{\int_0^1 (1-t)^{\frac{a+p}{2}-2 - \frac{n-1}{2} } dt}{ \int_0^1 (1-t)^{\frac{a+p_1}{2}-2 - \frac{n-1}{2}} dt } \\
&=& \frac{a+p_1-n-1}{a+p-n-1} < 1.
\end{eqnarray*}

\noindent
When $n < a+p_1-1$, the Bayes Factor $BF(M_2:M_1)$ is strictly less than 1 and still favors the smaller model $M_1$ in the limit, but not with overwhelming evidence as in the other two cases.

\subsection{Detailed  Proof of Theorem \ref{thm4.1}} \label{appDetail4.1}

In the block orthogonal setup
\[ \pi (\bfg \mid \bfy) \propto \frac{\prod_{j=1}^k (1+g_j)^{-\frac{a+p_j}{2}} }{\left[ 1 -  \sum_{j=1}^k \frac{g_j}{g_j+1} R_j^2 \right]^{(n-1)/2} } \]

\noindent
So for any $i=1,2,..,k$,
\begin{eqnarray*}
E \left(\frac{g_i}{1+g_i} \mid \bfy \right) &=& \frac{ \int_{(0,1)^k} t_i \prod_{j=1}^k (1-t_j)^{\frac{a+p_j}{2}-2} (1-\sum_{j=1}^k t_j R_j^2)^{-\frac{n-1}{2}} d\bft}{ \int_{(0,1)^k} \prod_{j=1}^k (1-t_j)^{\frac{a+p_j}{2}-2} (1-\sum_{j=1}^k t_j R_j^2)^{-\frac{n-1}{2}} d\bft } \\
&\geq& \frac{ \int_{(0,1)^k}  t_i  \prod_{j=1}^k (1-t_j)^{\frac{a+p_j}{2}-2} (1- t_i R_i^2)^{-\frac{n-1}{2}} d\bft}{ \int_{(0,1)^k} \prod_{j=1}^k (1-t_j)^{\frac{a+p_j}{2}-2} (1- t_i R_i^2)^{-\frac{n-1}{2}} d\bft } \\
&& \hspace{1.4 in} (\mbox{by Lemma }\ref{lem4.2}) \\
&=& \frac{2}{a+p_i} \; \frac{_2F_1(\frac{n-1}{2},2;\frac{a+p_i}{2}+1;R_i^2)}{_2F_1(\frac{n-1}{2},1;\frac{a+p_i}{2};R_i^2)} \\
\mbox{Hence } \lim_{N \rightarrow \infty} E \left(\frac{g_i}{1+g_i} \mid \bfy \right) &\geq&  \lim_{N \rightarrow \infty}  \frac{2}{a+p_i} \; \frac{_2F_1(\frac{n-1}{2},2;\frac{a+p_i}{2}+1;R_i^2)}{_2F_1(\frac{n-1}{2},1;\frac{a+p_i}{2};R_i^2)} .\\
\end{eqnarray*} 

\noindent
As $N \rightarrow \infty$, $R^2_1 \rightarrow 1$ so that for $i=1$,
\begin{eqnarray*}
\lim_{N \rightarrow \infty} E \left(\frac{g_1}{1+g_1} \mid \bfy \right) &\geq&  \lim_{z \rightarrow 1}  \frac{2}{a+p_1} \; \frac{_2F_1(\frac{n-1}{2},2;\frac{a+p_1}{2}+1;z)}{_2F_1(\frac{n-1}{2},1;\frac{a+p_1}{2};z)}  \\
&=& 1   
\end{eqnarray*}
when  $n \geq a+p_1-1 $ (see Theorem \ref{thm3.1}). But $E \left(\frac{g_1}{1+g_1} \mid \bfy \right) \leq 1$, implying that $E \left(\frac{g_1}{1+g_1} \mid \bfy \right) \rightarrow 1$ in the limit. 

For $i > 1$ and  $m \neq i$, $ E \left(\frac{g_i}{1+g_i} \mid \bfy \right) $
\begin{eqnarray*}
  &=&  \frac{ \int_{(0,1)^k} t_i \prod_{j=1}^k (1-t_j)^{\frac{a+p_j}{2}-2} (1-\sum_{j=1}^k t_j R_j^2)^{-\frac{n-1}{2}} d\bft}{ \int_{(0,1)^k} \prod_{j=1}^k (1-t_j)^{\frac{a+p_j}{2}-2} (1-\sum_{j=1}^k t_j R_j^2)^{-\frac{n-1}{2}} d\bft }  \\
&\geq& \frac{ \int_{(0,1)^k} t_i \prod_{j=1}^k (1-t_j)^{\frac{a+p_j}{2}-2} (1- t_m R_m^2)^{-\frac{n-1}{2}} d\bft}{ \int_{(0,1)^k} \prod_{j =1}^{k} (1-t_j)^{\frac{a+p_j}{2}-2} (1- t_m R_m^2)^{-\frac{n-1}{2}} d\bft } \; \; (\mbox{by Lemma }\ref{lem4.2}) \\
&=& \frac{\mbox{Beta}(2,\frac{a+p_i}{2}-1)}{\mbox{Beta}(1,\frac{a+p_i}{2}-1)} \times \frac{_2F_1(\frac{n-1}{2},1;\frac{a+p_m}{2};R^2_m)}{_2F_1(\frac{n-1}{2},1;\frac{a+p_m}{2};R^2_m)} = \frac{2}{a+p_i} 
\end{eqnarray*}
Thus, $\lim\limits_{N \rightarrow \infty} E \left(\frac{g_i}{1+g_i} \mid \bfy \right) \geq \frac{2}{a+p_i}$. Equality in the relation above is attained when $R_i^2 = 0$ or $R_i^2 \rightarrow 0$ and $\lim\limits_{N \rightarrow \infty} \sum\limits_{j \neq i} R^2_j < 1$, i.e, when no linear combination of the predictors explains all of the variation in the response.

\begin{eqnarray*}
E \left(\frac{g_i}{1+g_i} \mid \bfy \right) &=& \frac{ \int_{(0,1)^k} t_i \prod_{j=1}^k (1-t_j)^{\frac{a+p_j}{2}-2} (1-\sum_{j=1}^k t_j R_j^2)^{-\frac{n-1}{2}} d\bft}{ \int_{(0,1)^k} \prod_{j=1}^k (1-t_j)^{\frac{a+p_j}{2}-2} (1-\sum_{j=1}^k t_j R_j^2)^{-\frac{n-1}{2}} d\bft }  \\
&\leq& \frac{ \int_{(0,1)^k} t_i \prod_{j=1}^k (1-t_j)^{\frac{a+p_j}{2}-2} (1-\sum_{j \neq i} R_j^2 - t_i R_i^2)^{-\frac{n-1}{2}} d\bft}{ \int_{(0,1)^k} \prod_{j=1}^k (1-t_j)^{\frac{a+p_j}{2}-2} (1-\sum_{j \neq i} R_j^2 - t_i R_i^2)^{-\frac{n-1}{2}} d\bft } \\
&& \hspace{1 in} (\mbox{using a variation of Lemma } \ref{lem4.2}) \\
&=& \frac{ \int_0^1 t_i (1-t_i)^{\frac{a+p_j}{2}-2} (1-\sum_{j \neq i} R_j^2 - t_i R_i^2)^{-\frac{n-1}{2}} dt_i}{ \int_0^1 (1-t_i)^{\frac{a+p_j}{2}-2} (1-\sum_{j \neq i} R_j^2 - t_i R_i^2)^{-\frac{n-1}{2}} dt_i } 
\end{eqnarray*}

\begin{eqnarray*}
&=& \frac{ (1-\sum_{j \neq i} R_j^2)^{-(n-1)/2} }{ (1-\sum_{j \neq i} R_j^2)^{-(n-1)/2} } \;  \frac{ \int_0^1 t_i (1-t_i)^{\frac{a+p_j}{2}-2} (1- t_i \frac{R_i^2}{1-\sum_{j \neq i} R_j^2})^{-\frac{n-1}{2}} dt_i}{ \int_0^1 (1-t_i)^{\frac{a+p_j}{2}-2} (1 - t_i \frac{R_i^2}{1-\sum_{j \neq i} R_j^2})^{-\frac{n-1}{2}} dt_i } \\
&=& \frac{2}{a+p_i} \frac{_2F_1(\frac{n-1}{2},2;\frac{a+p_i}{2}+1;\kappa_i)}{_2F_1(\frac{n-1}{2},1;\frac{a+p_i}{2};\kappa_i)}  < 1
\end{eqnarray*}
\noindent
where $\kappa_i = \frac{R_i^2}{1-\sum_{j \neq i} R_j^2}$.

For this sequence of problems, $0<\kappa_i<1$ is fixed for all $i \neq 1$, because
\[ \kappa_i = \frac{R_i^2}{1-\sum_{j \neq i} R_j^2} =\frac{\bfy^T P_{X_i} \bfy}{(n-p-1)\hat{\sigma}^2+\bfy^T P_{X_i} \bfy} 
\]

\noindent
If $\bfy^T P_{X_i} \bfy$ is quite large compared to $\hat{\sigma}^2$, then $\kappa_i \approx 1$ and the shrinkage factor is near 1 while very small values of $\bfy^T P_{X_i} \bfy$ relative to $\hat{\sigma}^2$ implies $\kappa_i \approx 0$ and the shrinkage factor is near the lower bound $\frac{2}{a+p_i}$. For all values of $0<\kappa_i<1$, we have
\[ \lim_{N \rightarrow \infty} E \left(\frac{g_i}{1+g_i} \mid \bfy \right) < 1, \; \mbox{for} \; i\neq 1.
\]

\subsection{Detailed Proof of Theorem \ref{thm4.2}}\label{appb45}

The Bayes Factors $BF(M_i:\Mo)$ comparing the models $M_i, i=1,2$ to the null (intercept only) model are
\begin{eqnarray*}
BF(M_2:\Mo) &=& \left(\frac{a-2}{2} \right)^2 \int_{0}^1 \int_{0}^1 \prod_{i=1}^2 (1-t_i)^{\frac{a+p_i}{2}-2} (1-\sum_{i=1}^2 t_i R_i^2)^{-\frac{n-1}{2}} dt_1 dt_2 \\
\mbox{and } BF(M_1:\Mo) &=& \frac{a-2}{2} \int_{0}^1 (1-t_1)^{\frac{a+p_1}{2}-2} (1- t_1 R_1^2)^{-\frac{n-1}{2}} dt_1 
\end{eqnarray*}

\begin{eqnarray*}
\mbox{Thus, } BF(M_2:M_1) &=& \frac{BF(M_2:\Mo)}{BF(M_1:\Mo)}   \\
&=& \frac{a-2}{2} \frac{ \int_{0}^1 \int_{0}^1 \prod_{i=1}^2 (1-t_i)^{\frac{a+p_i}{2}-2} (1-\sum_{i=1}^2 t_i R_i^2)^{-\frac{n-1}{2}}  dt_1 dt_2 }{\int_{0}^1 (1-t_1)^{\frac{a+p_1}{2}-2} (1- t_1 R_1^2)^{-\frac{n-1}{2}} dt_1 } \hspace{1 cm} \\
&\geq& \frac{a-2}{2} \frac{\int_{0}^1 \int_{0}^1 \prod_{i=1}^2 \left[ (1-t_i)^{\frac{a+p_i}{2}-2} (1-t_i R_i^2)^{-\frac{n-1}{2}} \right] dt_1 dt_2}{ \int_0^1 (1-t_1)^{\frac{a+p_1}{2}-2} (1- t_1 R_1^2)^{-\frac{n-1}{2}} dt_1 } \\
&=& \frac{a-2}{2} \int_0^1 (1-t_2)^{\frac{a+p_2}{2}-2} (1- t_2 R_2^2)^{-\frac{n-1}{2}} dt_2
\end{eqnarray*}

\noindent
The above inequality comes from the following result which holds for any $k \in \mathbb{N}$
\begin{equation}\label{eq:usefulineq}
1-\sum_{i=1}^k x_i \leq \prod_{i=1}^k (1-x_i) \; \; ,  \mbox{ when } 0\leq x_i \leq 1 \; \forall \; i
\end{equation} 

As $||\beta_1|| \rightarrow \infty$, $R_1^2 \rightarrow 1$ and $R_2^2 \rightarrow 0$, and so
\[ \lim_{||\beta_1|| \rightarrow \infty} BF(M_2:M_1) \geq \frac{a-2}{2} \int_0^1 (1-t_2)^{\frac{a+p_2}{2}-2} dt_2 = \frac{a-2}{a+p_2-2}
\]

\noindent
Thus the limiting Bayes Factor is bounded away from zero. The lower bound is decreasing in $p_2$, the number of additional predictors in the superset model.

\end{document}